\documentclass[aap]{imsart}
 
 \RequirePackage{amsthm,amsmath,amsfonts,amssymb}
 \RequirePackage[numbers,sort&compress]{natbib}
 \RequirePackage[colorlinks,citecolor=blue,urlcolor=blue]{hyperref}
 \RequirePackage{graphicx}
  
 \usepackage{tikz} 
 \usepackage{mathrsfs,enumitem} 
 
 \usepackage{pgfplots}
 
 \startlocaldefs
 \theoremstyle{plain}

 \newtheorem{theorem}{Theorem}[section]
 \newtheorem{lemma}[theorem]{Lemma}
 \theoremstyle{remark}
 \newtheorem{definition}[theorem]{Definition}

 \newtheorem{corollary}[theorem]{Corollary}
 \newtheorem{proposition}[theorem]{Proposition}
 \newtheorem{remark}[theorem]{Remark}
 \newtheorem{condition}[theorem]{Condition}

  \makeatletter
 \def\namedlabel#1#2{\begingroup
 	#2%
 	\def\@currentlabel{#2}%
 	\phantomsection\label{#1}\endgroup
 }
 \makeatother

 \def\beqlb{\begin{eqnarray}}\def\eeqlb{\end{eqnarray}} 
 \def\beqnn{\begin{eqnarray*}}\def\eeqnn{\end{eqnarray*}} 
 \def\ar{&} 

 \def\proof{\noindent{\it Proof.~~}} 
 \def\qed{\hfill$\Box$\medskip}
 
 \endlocaldefs
 
\begin{document}
 
 \begin{frontmatter}
 \title{Stochastic Volterra Equations for the Local Times of Spectrally Positive Stable Processes}
 \runtitle{SVE for the Local Times of Stable Processes}

 \begin{aug}
 \author{\fnms{Wei}~\snm{Xu}\ead[label=e1]{xuwei.math@gmail.com}\orcid{0000-0002-9370-8846}}
 \address{School of Mathematics and Statistics, Beijing Institute of Technology\printead[presep={,\ }]{e1}} 
 \end{aug}
 
 \begin{abstract}
 This paper is concerned with the evolution dynamics of local times of a spectrally positive stable process in the spatial direction.
 The main results state that conditioned on the finiteness of the first time at which the local time at zero exceeds a given value, the local times at positive half line are equal in distribution to the unique solution of a stochastic Volterra equation  driven by a Poisson random measure whose intensity coincides with the L\'evy measure.
 This helps us to provide not only a simple proof for the H\"older regularity, but also a uniform upper bound for all moments of the H\"older coefficient as well as a maximal inequality for the local times.
 Moreover, based on this stochastic Volterra equation, we extend the method of duality to establish an exponential-affine representation of the Laplace functional in terms of the unique solution of a nonlinear Volterra integral equation associated with the Laplace exponent of the stable process.
 \end{abstract}

 \begin{keyword}[class=MSC2020]
 \kwd[Primary ]{60G52}
 \kwd{60J55}
 \kwd{60H20}
 \kwd[; secondary ]{60G22}
 \kwd{60F17} 
 \kwd{60G55}
 \end{keyword}

 \begin{keyword}
 \kwd{Local time}
 \kwd{stable process}
 \kwd{stochastic Volterra equation}
 \kwd{heavy tail}
 \kwd{Poisson random measure}
 \kwd{marked Hawkes point measure}
 \kwd{Laplace functional}
 \kwd{Ray-Knight theorem}
 \end{keyword}

 \end{frontmatter}
 \tableofcontents

 \section{Introduction}
  
 Local times of L\'evy processes not only have wide applications in various fields; see \cite{BarndorffMikoschResnick2001,Kyprianou2014}, but they also have been studied in depth with abundant of interesting results obtained, e.g. various constructions (see \cite{BarlowPerkinsTaylor1986,DellacherieMeyer1980,Protter2005}), Hilbert transform (see \cite{Bertoin1995,FitzsimmonsGetoor1992}), H\"older regularity (see \cite{Barlow1988,Boylan1964,FormanPalRizzoloWinkel2018}) and so on.
 We refer to \cite{Bertoin1996,GemanHorowitz1980,RevuzYor2005} for survey on local times and their applications.
 In particular, to understand thoroughly the dependence of Brownian local times in the space variable, Ray \cite{Ray1963} and Knight \cite{Knight1963} independently proved the well-known  {\it Ray-Knight theorem} that links Brownian local times to Bessel processes. It was generalized in \cite{EisenbaumKaspiMarcusRosenShi2000,SabotTarres2016} to  strongly symmetric Markov processes with finite $1$-potential densities. 
 In a recent work, A\"{i}d\'ekon et al. \cite{AidekonHuShi2024} reproved the Ray-Knight theorem for perturbed reflecting Brownian motion by means of Tanaka's formula and Walsh's stochastic integral with respect to a Gaussian white noise.

 For a general spectrally positive L\'evy process, Le Gall and Le Jan \cite{LeGallLeJan1998a,LeGallLeJan1998b} considered the reflected processes of its time-reversed processes.
 Associated to the local times at $0$, they introduced an exploration process to describe the genealogy of a \textit{continuous-state branching process} (CB-process) and generalized an analogue of the Ray-Knight theorem for a functional of local times of the L\'evy process;  see also \cite{DuquesneLeGall2002,LeGall1999} for details.

 Because of the lack of Markovianity; see \cite{EisenbaumKaspi1993}, local times of general spectrally positive L\'evy processes seem to be more untractable than their functional considered in \cite{DuquesneLeGall2002,LeGall1999}. 
 Their microstructure and evolution mechanism have received considerable attention in recent years.
 Specifically, Lambert \cite{Lambert2010} connected a compound Poisson process with unit negative drift and killed upon hitting $0$ to the jumping chronological contour processes of a splitting tree, and then showed that its local times are equal in distribution to a homogeneous, binary \textit{Crump-Mode-Jagers branching process} (CMJ-process).
 For a general spectrally positive L\'evy process, Lambert and Simatos \cite{LambertSimatos2015} explored the genealogical structure of their local times preliminarily via an approximating sequence consisting of rescaled binary CMJ-processes.
 Later, a detailed genealogical interpretation was given in \cite{LambertBravo2018} by considering the corresponding totally ordered measured tree that satisfies the splitting property.
 Meanwhile, Forman et al. \cite{FormanPalRizzoloWinkel2018} established a locally uniform approximation for the local times of a driftless spectrally positive stable process by endowing each jump with a random graph.
 Up to now,  the genealogical structure of local times of general spectrally positive L\'evy processes seems to be fairly clear.
 However, their macroevolution mechanisms, by contrast, are still incomprehensible.

 The purpose of this work is to establish stochastic equations for the local times of spectrally positive stable processes and study their evolution dynamics in the spatial direction.
 In contrast to the genealogical interpretations given in the aforementioned  literature,  stochastic equations have many advantages including
 \begin{enumerate}
	\item[$\bullet$] They provide an intuitive description for the evolution of local times in the spatial direction as well as a detailed interpretation of their perturbations caused by each jump of stable processes. \vspace{5pt}
	
	\item[$\bullet$] They allow us to study the local times by using tools and methods from the modern probability theory, e.g., stochastic integral inequalities, stochastic Fubini theorem and extreme value theory.\vspace{5pt}
	
	\item[$\bullet$] They offer a kind of novel non-Markovian models and a convenient way of numerical analysis, which will benefit greatly the related fields, e.g., processor-sharing queues and stochastic volatility models.
 \end{enumerate}

 \subsection{Overview of main results}
 Let $\xi:=\{ \xi(t):t\geq 0 \}$ be a one-dimensional spectrally positive stable process  with index $1+\alpha\in(1,2)$ and Laplace exponent
 \beqlb\label{LaplaceExp}
 {\it\Phi}(\lambda) :=  b\lambda +c\lambda^{\alpha+1} = b\lambda + \int_0^\infty \big( e^{-\lambda y}-1+\lambda y \big)\nu_\alpha(dy), \quad \lambda \geq 0,
 \eeqlb
 where $b\geq 0$, $c>0$ and $\nu_\alpha(dy)$, known as the \textit{L\'evy measure}, is a $\sigma$-finite measure on $(0,\infty)$ given by 
 \beqlb\label{LevyMeasure}
 \nu_\alpha(dy):= \frac{c\alpha(\alpha+1)}{\Gamma(1-\alpha)} \cdot y^{-\alpha-2}\cdot dy.
 \eeqlb
 It is recurrent or drifts to $-\infty$ according as $b=0$ or $>0$.
 Let $W$ be the scale function of $\xi$ and $\bar{\nu}_\alpha(x):=  \nu_\alpha([x,\infty))$ the tail function of $\nu_\alpha$.
 Let $L_\xi:=\{L_\xi(x,t):x\in \mathbb{R},t\geq 0 \}$ be the local times of $\xi$, where $L_\xi(x,t)$ is usually interpreted as the amount of time that $\xi$ spends at level $x$ up to time $t$. 
 Denote by $\tau_\xi^L(\zeta)$ the first time that the amount of local time accumulated at level $0$ exceeds a given value $\zeta>0$; more accurate definitions can be found in Section~\ref{StableProc} and \cite{Bertoin1996,Kyprianou2014}.

 Let $L^\xi_\zeta$ be the process $\{L_\xi(x,\tau_\xi^L(\zeta)):x\geq 0  \}$ conditioned on $\tau_\xi^L(\zeta)<\infty$.
 The first main result states that $L^\xi_\zeta$ is the \textit{unique weak solution}\footnote{
 A continuous process with distribution $\mathcal{P}$ is called a {\it weak solution} of (\ref{MainThm.SVE}) if there exists a stochastic basis,  a PRM  $N_0( dy,dz) $ on $(0,\infty)^2$ with intensity $  \bar{\nu}_\alpha(y)\cdot dy \cdot dz$, a PRM $N_\alpha(ds,dy,dz)$ on $(0,\infty)^3$ independent of $N_0(dy,dz)$ with intensity $ds\cdot\nu_\alpha(dy)\cdot dz$ and a continuous process $L^\xi_\zeta$ with distribution $\mathcal{P}$ such that (\ref{MainThm.SVE}) holds almost surely.
 We say the {\it weak uniqueness} holds if any two weak solutions are equal in distribution.
 \label{Def.WeakSolution}} 
 of the \textit{stochastic Volterra equation} (SVE)
 \beqlb\label{MainThm.SVE}
 L^\xi_\zeta(x)\ar=\ar   \int_0^\infty  \int_0^\zeta \nabla_y W(x)   N_0( dy,dz)  \cr
 \ar\ar + \int_0^x \int_0^\infty \int_0^{L^\xi_\zeta(s)}  \nabla_y W(x-s)   \widetilde{N}_\alpha(ds,dy,dz), \quad x\geq 0,
 \eeqlb
 where $\nabla_y W(x):= W(x)-W(x-y)$, $N_0( dy,dz) $ is a Poisson random measure (PRM) on $(0,\infty)^2$ with intensity $  \bar{\nu}_\alpha(y) \cdot dy \cdot dz$,
 $\widetilde{N}_\alpha(ds,dy,dz)$ is a compensated PRM on $(0,\infty)^3$ with intensity $ ds\cdot \nu_\alpha(dy)\cdot dz$ and independent of $N_0(dy,dz)$.
 The first stochastic integral in (\ref{MainThm.SVE})  represents the contribution of jumps up-crossing $0$ to the local time at level $x$ and the second stochastic integral, known as \textit{stochastic Volterra integral} (SVI), can be interpreted as the perturbations caused by jumps with initial positions above $0$.
 Since the convolution kernel delays the relaxation of its perturbations, the PRM $N_\alpha(ds,dy,dz)$   changes the local times continuously in the spatial variable. This stands in striking contrast to the jumps in \textit{It\^o's stochastic differential equations} (It\^o's SDEs) driven by PRM.
 Additionally, because of the joint impact of relative level $x-s$ and jump-size $y$ on the convolution kernel, the SVE (\ref{MainThm.SVE}) cannot be written into the form of SVEs in  \cite{Jaber2021,JaberCuchieroLarssonPulido2021,JaberLarssonPulido2019,PardouxProtter1990,Protter1985}.

 Based on the SVE (\ref{MainThm.SVE}), in the second main result we use stochastic integral inequalities to provide a simple proof for the H\"older continuity of $L^\xi_\zeta$ and the finiteness of all moments of the H\"older coefficient given in  \cite{Barlow1988,Boylan1964,FormanPalRizzoloWinkel2018}.
 As the novelty, we also establish a uniform upper bound for all moments of the H\"older coefficient and a maximal inequality for the local times in the spatial variable.
 With the crucial assistance from the SVE (\ref{MainThm.SVE}), in the third main result we extend the method of duality developed in \cite{JaberLarssonPulido2019} to provide an explicit representation of the Laplace functional 
 \beqnn
 \mathbf{E}  \Big[\exp  \Big\{- \lambda\cdot L^\xi_\zeta(x) -  g * L^\xi_\zeta(x)  \Big\} \Big]
 \eeqnn
 with $ \lambda \geq 0 $ and $ g\in L^\infty(\mathbb{R}_+;\mathbb{R}_+)$.
 It states that the Laplace exponent can be written as an affine functional of the initial state, in terms of the unique solution of the \textit{nonlinear Volterra integral equation} (nonlinear-VIE)
 \beqlb\label{MainThm.Volterra}
 v_\lambda^g(x) =  \lambda W'(x) + \big(g -\mathcal{V}_\alpha\circ v_\lambda^g \big)*W'(x),\quad x>0,
 \eeqlb
 where $W'$ is the derivative of $W$ and $\mathcal{V}_\alpha$ is a nonlinear operator acting on a locally integrable function $f$  by
 \beqlb\label{OperatorV}
 \mathcal{V}_\alpha\circ f(x) := \int_0^\infty
 \Big(\exp\Big\{- \int_{(x-y)^+}^{x} f(r)dr\Big\}-1  +\int_{(x-y)^+}^{x} f(r)dr\Big) \nu_\alpha( dy ),
 \eeqlb
 for $x\geq 0$. 
 Finally, we provide an alternative fractional integration and differential equation for the process $L^\xi_\zeta$ and its Laplace exponent. In contrast to the SVE (\ref{MainThm.SVE}), the  alternative equation takes it a step further and extracts the impact of drift $b$ on the local times from that of jumps.
 It also uncovers the remarkable similarity between $L^\xi_\zeta$ and CB-processes in the evolution mechanism, which, together with the genealogical
 interpretations in \cite{FormanPalRizzoloWinkel2018,Lambert2010,LambertBravo2018}, tells that the SVE (\ref{MainThm.SVE}) defines a novel non-Markovian CB-process.

 To illustrate the strength of these results,  we use the SVE (\ref{MainThm.SVE}) to establish a stochastic equation for the heavy-traffic limit of recaled queue-length processes of M/G/1 processor-sharing queues with unit service capacity, heavy-tailed service distribution and stopped upon becoming empty.
 It can be seen as a continuation of \cite{LambertSimatos2015}, where  the weak convergence of rescaled queue-length processes was proved.
 In a sense, this helps to partially answer Problem~2 stated by Zwart in \cite{Zwart2022} about the heavy-traffic limit of heavy-tailed processor-sharing queues; readers may refer to the references of Zwart and his coauthors for details.
 Enlightened by the self-exciting property observed in the SVE (\ref{MainThm.SVE}), in the forthcoming preprint \cite{HorstXuZhang2023} we use the evolution mechanism of local times of stable processes to model the sharp-raise clusters in rough volatility and introduce a novel fractional stochastic volatility model with self-excited sharp-raises.

 \subsection{Methodologies} 
 We start the construction of the SVE (\ref{MainThm.SVE}) from the result that the local times of nearly recurrent compound Poisson processes with unit negative drift, Pareto-distributed jumps are equal in distribution to a class of nearly critical binary CMJ-processes, which converge weakly to the process $L^\xi_\zeta$ after rescaling; see \cite{Lambert2010,LambertSimatos2015}.
 Enlightened by the Hawkes representation of general branching particle systems established in \cite{HorstXu2021,Xu2021}, we reconstruct the binary CMJ-processes as the intensity processes of nearly unstable \textit{marked Hawkes point measures} (MHPs) by translating the birth time, life-length and survival state of each individual into the arrival time, random mark and kernel of an event respectively.
 Furthermore, we write each intensity process into a SVE driven by an infinite-dimensional martingale in which the integrand is a functional of the \textit{resolvent function} related to the life-length distribution.
 Consequently, it suffices to prove the weak convergence of these SVEs after rescaling to the desired SVE (\ref{MainThm.SVE}).
 Unfortunately, the Pareto-distributed life-length gives raise to long-range dependence in the pre-limit SVEs, which derives a series of challenges and difficulties in the proof including
 \begin{enumerate}
	\item[$\bullet$] Along with the inseparable impact of time and life-length on the convolution kernel, the infinite-dimensional driving noises not only lead to the failure of the approximation method and the integral-derivative method developed in \cite{JaberCuchieroLarssonPulido2021,JaberLarssonPulido2019,JaissonRosenbaum2016}, but also make it hard to seek an approximation for the pre-limit SVEs.
	\vspace{5pt}
	
	\item[$\bullet$] The resolvent function fluctuates drastically and explodes around $0$ after rescaling.
	This leads to the sharp swings in the cumulative impact of infinite short-lived events on the pre-limit SVEs and also makes the uniform control on the error processes challengeable.
	\vspace{5pt}
	
	\item[$\bullet$] The resolvent function inherits long-range dependence from the life-length distribution. It prevents us from transforming the pre-limit SVEs into the form of It\^o's SDEs and obtaining the weak convergence similarly as in \cite{JaissonRosenbaum2015,Xu2021} by using the weak convergence results established in  \cite{KurtzProtter1991,KurtzProtter1996} for It\^o's SDEs. 
 \end{enumerate}
 To overcome the first two difficulties, we start by analyzing in depth the direct and indirect impact of each event on the pre-limit SVEs.
 Our analyses show that the cumulative direct impact of all events can be asymptotically ignored and a suitably rescaled version of their indirect impact asymptotically behaves as the backward difference of scale function.
 This motivates us to approximate the SVIs in the pre-limit SVEs by replacing the integrands with the backward difference of scale function.
 For the uniform control on the error processes, we first split them into several parts according to the source and then prove the finite-dimensional convergence of each part to $0$ separately.
 Based on a deep analysis about the backward difference of scale function, we prove the $C$-tightness\footnote{Readers may refer to Definition~3.25 in \cite[p.351]{JacodShiryaev2003} for  the definition of $C$-tightness.} 
 of the approximating processes, which, together with the $C$-tightness result given in \cite{LambertSimatos2015} for the local times of nearly recurrent compound Poisson processes, yields the tightness of error processes.
 To overcome the third difficulty, we establish a weak convergence result for SVIs with respect to infinite-dimensional martingales, whose tightness and finite-dimensional convergence are obtained from the foregoing tightness results and the weak convergence of the related It\^o's stochastic integrals respectively.
 More precisely, for a given finite sequence of time points, we first introduce a sequence of It\^o's stochastic integrals with respective to infinite-dimensional martingale satisfying that their finite-dimensional distributions at the given time points are equal to those of the corresponding SVIs, and then prove their weak convergence to a limit process whose finite-dimensional distribution at the given time points is equal to that of the desired limit SVI.

 In the proof of existence and uniqueness of solutions of the nonlinear-VIE (\ref{MainThm.Volterra}), the next two main difficulties steam from the nonlinear operator $\mathcal{V}_\alpha$ and the singularity of the function $W'$ at the origin
 \begin{enumerate}
	\item[$\bullet$] The interplay between the singularity of $W'$ and $\mathcal{V}_\alpha$  makes the existence of local solutions of (\ref{MainThm.Volterra}) around $0$ quite difficult.
	\vspace{5pt}
	
	\item[$\bullet$]
	Since $\mathcal{V}_\alpha$ is path-dependent and does not satisfy the Lipschitz condition, it is  difficult to identify the non-explosion of local solutions and extend them into global solutions.
 \end{enumerate}
 To bypass the first difficulty, we first prejudge the behavior of solutions near the origin with the help of an upper bound estimate of $\mathcal{V}_\alpha$ and the expansion given in \cite{CallegaroGrasselliPages2021} for solutions of fractional Riccati equations.
 In a specified closed set in some Lebesgue space, we then find a local solution of (\ref{MainThm.Volterra}) successfully by using Banach's fixed point theorem.
 To overcome the second difficulty, associated with a fractional differential equation related to $\mathcal{V}_\alpha$ we first provide an upper bound estimate for a functional of each local solution, and then, along with the comparison principle for fractional differential equations, establish a uniform control on the local solutions.

 \subsection{Related Literature}
 Let us comment on the relationship between the present work and the existing literature.
 Firstly, based on the Markov property, Brownian local times were linked to Bessel processes via their transition semigroups in \cite{Knight1963,Ray1963} or their infinitesimal generators in \cite{KawazuWatanabe1971,Li2006}.
 However, the lack of Markovianity of $L^\xi_\zeta$ makes it impossible to establish the SVE (\ref{MainThm.Volterra}) similarly as in the preceding references.
 Even if it could be established successfully, the SVE (\ref{MainThm.Volterra}) is beyond the scope of all existing literature \cite{Jaber2021,JaberCuchieroLarssonPulido2021,JaberLarssonPulido2019,PardouxProtter1990,Protter1985} and the existence of its solutions seems to be quite difficult to be proved in the standard way.
 On the other hand, the present work establishes the well-posedness of the novel SVE (\ref{MainThm.Volterra}).
 Secondly, the main results, as mentioned above, are obtained by establishing a weak convergence result for the corresponding long-range dependent MHPs.
 The first scaling limit theorem for Hawkes processes was  established by Jaisson and Rosenbaum \cite{JaissonRosenbaum2015} in the study of the asymptotic behavior of Hawkes-based price-volatility models in the context of high-frequency trading.
 Their results state that under the short-memory condition,  the rescaled intensity processes of nearly unstable Hawkes processes converge weakly to the well-known \textit{CIR-model}.  The analogous scaling limits were established for multivariate (marked) Hawkes processes in \cite{ElEuchFukasawaRosenbaum2018,Xu2021} and a jump-diffusion limit was given in \cite{HorstXu2022} for MHPs with exponential kernel.
 When the kernel is heavy-tailed, Jaisson and Rosenbaum \cite{JaissonRosenbaum2016} proved the weak convergence of the integral of rescaled intensity process to the integral of a fractional diffusion process, see also \cite{ElEuchFukasawaRosenbaum2018, RosenbaumTomas2021} for the multivariate case.
 Because of many difficulties deriving from long-range dependence, they left the weak convergence of rescaled intensity processes as an open problem, which was recently solved by \cite{HXZ2023}. 
 However, we stress that the weak convergence result in this work is established for the intensity processes of MHPs.
 As the final remark, we need to point out that different to the analogous version in \cite{LeGallLeJan1998a,LeGallLeJan1998b}, the Ray-Knight theorem in this work is established for the local times rather than their functionals.

 \smallskip
 {\it \textbf{Organization of this paper.}} \ \  In Section~\ref{preliminaries}, we first introduce general notation and properties of spectrally positive stable processes, and then formulate the main results.
 In Section~\ref{CompoundPoissonP}, we introduce some elementary results and a SVE for the local times of a compound Poisson process with negative drift by linking them to a MHP.
 Section~\ref{ProofForThm.1} is devoted to proving that $L^\xi_\zeta$ solves the SVE (\ref{MainThm.SVE}).
 Its H\"older continuity is proved in Section~\ref{ProofHolder}.
 In Section~\ref{NonLinearVolUnique}, we prove the exponential-affine transform formula as well as the existence and uniqueness of solutions of the nonlinear-VIE (\ref{MainThm.Volterra}).
 The proof for the alternative representation of $L^\xi_\zeta$ are given in Section~\ref{FIR}.
 Applications to processor-sharing queues are given in Section~\ref{Sec.PS}.
 Additional proofs and supporting results are presented in the Appendices.

 \smallskip
 {\it \textbf{Notation.}} \ \
 For any $x\in\mathbb{R}$, let $x^+:=x\vee 0$, $x^-:=x\wedge 0$ and $[x]$ be the integer part of $x$.
 For a Banach space $\mathbb{V}$ with a norm $\|\cdot\|_\mathbb{V}$, let $D([0,\infty),\mathbb{V})$ be the space of all c\'adl\'ag $\mathbb{V}$-valued functions endowed with the Skorokhod topology and $C([0,\infty),\mathbb{V})$ the space of all continuous $\mathbb{V}$-valued functions endowed with the uniform topology.
 For any $\mathcal{T}\subset[0,\infty)$ and $p\in(0,\infty]$, let $L^p(\mathcal{T}; \mathbb{V})$ be the space of $\mathbb{V}$-valued measurable functions $f$ on $\mathcal{T}$ satisfying   
 \beqnn
 \big\|f\big\|_{L^p_\mathcal{T}}^p:=  \int_\mathcal{T} \big\|f(x)\big\|_{\mathbb{V}}^p dx  <\infty.
 \eeqnn
 We also write $\|f\|_{L^p_T}$ for $\|f\|_{L^p_{[0,T]}}$ and $\|f\|_{L^{^p}}$ for $\|f\|_{L^p_\infty}$.
 We make the conventions that for $x, y\in \mathbb{R}$ with $y\geq x$,
 \beqnn
 \int_x^y=-\int_y^x= \int_{(x,y]},\quad \int_{x-}^{y-}=  \int_{[x,y)}
 \quad\mbox{and}\quad
 \int_x^\infty = \int_{(x,\infty)}.
 \eeqnn

 Denote by $f*g$  the convolution of two functions $f,g$ on $\mathbb{R}_+$.
 Let $\Delta_h$ and $\nabla_h$ be the forward and backward difference operators with step size $h>0$, i.e., 
 \beqnn
 \Delta_hf(x):= f(x+h)-f(x)
 \quad\mbox{and}\quad
 \nabla_hf(x):=f(x)-f(x-h).
 \eeqnn 
 Let $\overset{\rm u.c.}\to$,  $\overset{\rm a.s.}\to$, $\overset{\rm d}\to$  and  $\overset{\rm p}\to$ be the uniform convergence on compacts, almost sure convergence,  convergence in distribution and convergence in probability respectively.
 We also use $\overset{\rm a.s.}=$, $\overset{\rm d}=$ and $\overset{\rm p}=$ to denote almost sure equality, equality in distribution and equality in probability respectively.

 For a probability measure $\mu$ on $\mathbb{R}$, denote by $\mathbf{P}_\mu$ and $\mathbf{E}_\mu$ the law and expectation of the underlying process with initial state distributed as $\mu$.
 When $\mu$ is a Dirac measure at point $x\in\mathbb{R}$, we write $\mathbf{P}_x$ for $\mathbf{P}_\mu$ and $\mathbf{E}_x$ for $\mathbf{E}_{\mu}$.
 For simplicity, we also write $\mathbf{P}$ for  $\mathbf{P}_0$ and $\mathbf{E}$ for $\mathbf{E}_0$.
 For two $\sigma$-finite measures $\mu_1,\mu_2$ on $\mathbb{R}$, we say $\mu_1 \leq \mu_2$ if for any non-negative function $f$ on $\mathbb{R}$,
 \beqnn
 \int_\mathbb{R} f(x) \mu_1(dx) \leq \int_\mathbb{R} f(x) \mu_2(dx). 
 \eeqnn
 
 We use $C$ to denote a positive constant whose value might change from line to line.

 \section{Preliminaries and main results}
 \label{preliminaries}

 \subsection{Spectrally positive stable processes}\label{StableProc}
 Suppose that the spectrally positive stable process $\xi$ is defined on a complete probability space $(\Omega, \mathscr{F}, \mathbf{P})$ endowed with a filtration $\{ \mathscr{F}_t \}_{t\geq 0}$ satisfying the usual hypotheses.
 For every $t\geq 0$, let $\mu_{\xi,t}(dy)$ be the \textit{occupation measure} of $\xi$ on the time interval $[0,t]$ given for every non-negative and measurable function $f$ on $\mathbb{R}$ by
 \beqnn
 \int_0^t f\big(\xi(s)\big)ds \overset{\rm a.s.}= \int_\mathbb{R} f(y) \mu_{\xi,t}(dy).
 \eeqnn
 The measure $\mu_{\xi,t}$ is absolutely continuous with respect to the Lebesgue measure and the density, denoted by  $\{L_\xi(x,t):x\in \mathbb{R} \}$, is square integrable; see Theorem~1 in \cite[p.126]{Bertoin1996}.
 The quantity $L_\xi(x,t)$ is called the {\it local time} of $\xi$ at level $x$ and time $t$.
 The two-parameter process  $L_\xi:=\{L_\xi(x,t):x\in \mathbb{R}, t\geq 0 \}$ is jointly continuous and satisfies the \textit{occupation density formula}
 \beqlb\label{OccupationDensityF}
 \int_0^t f(\xi(r))dr \overset{\rm a.s.}= \int_\mathbb{R} f(x) L_\xi(x,t) dx, \quad t\geq 0,
 \eeqlb
 see Theorem~15 in \cite[p.149]{Bertoin1996}.
 Moreover, for any $(\mathscr{F}_t)$-stopping time $\tau$, it is easy to identify that
 \beqlb\label{HightMaximum}
 \inf\big\{x\geq 0:L_\xi(x,\tau) =0 \big\} \overset{\rm a.s.}= \sup\big\{  \xi(t): t\in[0,\tau]  \big\}.
 \eeqlb
 The process $\{L_\xi(0,t):t\geq 0\}$ is  continuous and non-decreasing.
 This allows us to define the \textit{inverse local time} $\tau_\xi^L:=\{\tau_\xi^L(\zeta):\zeta\geq 0\}$ at level $0$ by
 $\tau_\xi^L(\zeta)=\infty$ if $\zeta> L_\xi(0,\infty)$ and
 \beqnn
 \tau_\xi^L(\zeta):= \inf\big\{ s\geq 0:  L_\xi(0,s)\geq \zeta \big\}, \quad \mbox{if } \zeta\in[0,L_\xi(0,\infty)].
 \eeqnn
 From Proposition~4 in \cite[p.130]{Bertoin1996}, the process $\tau_\xi^L$ is a subordinator, killed at an independent exponential time if $\xi$ is transient ($b> 0$), and its Laplace transform is of the form
 \beqlb\label{eqn.500}
 \mathbf{E}\big[\exp\big\{-\lambda\cdot \tau_\xi^L(\zeta)\big\}\big] = \exp\big\{-\zeta/u^\lambda(0)\big\},\quad \lambda > 0,\, \zeta\geq 0,
 \eeqlb
 where $u^\lambda:=\{u^\lambda(y):y\in\mathbb{R}\}$ is the density of the \textit{$\lambda$-resolvent kernel} of $\xi$.
 When $b=0$, we have $L_\xi(0,\infty)\overset{\rm a.s.}=\infty$ and $ \tau_\xi^L(\zeta) <\infty$ a.s.
 When $b>0$,  the \textit{potential density} $u^0$ is well-defined as the limit case $\lambda=0$ for $u^\lambda$.
 In this case, we  consider  the limit case $\lambda\to 0+$ for (\ref{eqn.500}) to get $$\mathbf{P}\big(L_\xi(0,\infty)\geq \zeta\big)=1-\mathbf{P}\big(\tau_\xi^L(\zeta) =\infty\big)= \exp\big\{-\zeta/u^0(0)\big\},$$
 which induces that $L_\xi(0,\infty)$ is exponentially distributed with mean $u^0(0)$ and $\mathbf{P} \big(\tau_\xi^L(\zeta) =\infty\big)>0$ for any $\zeta>0$.
 For $a,\theta >0$, let $a\cdot \xi(\theta\cdot):= \{a\xi(\theta t):t\geq 0\}$.
 The equality (\ref{OccupationDensityF}), along with the change of variables, implies the following two equivalences 
 \beqlb\label{TimeSpatialChange01}
 L_{a\cdot\xi(\theta\cdot)}  \overset{\rm a.s.}= \big\{(a\theta)^{-1} \cdot L_{\xi}(x/a,\theta t):x\in\mathbb{R},t\geq 0 \big\}
 \eeqlb
 and 
 \beqlb\label{TimeSpatialChange02}
 \tau^L_{a\cdot\xi(\theta\cdot)}(\zeta) \overset{\rm a.s.}= \theta^{-1}\cdot  \tau^L_{\xi}(a\theta\zeta) .
 \eeqlb
 Actually, these two equivalences also hold for any L\'evy process. 
 



 Let $\{W(x):x\in \mathbb{R} \}$ be the \textit{scale function} of $\xi$, which is identically zero on $  (-\infty,0)$ and characterized on $[0,\infty)$ as a strictly increasing function whose Laplace transform is given by
 \beqlb\label{ScaleFunction}
 \int_0^\infty e^{-\lambda x}W(x)dx= \frac{1}{\mathit{\Phi}(\lambda )},
 \quad \lambda >0.
 \eeqlb
 The scale function $W$ is continuous on $\mathbb{R}$ and differentiable on $(0,\infty)$ with derivative denoted as $W'$; see Theorem~8 in \cite[p.194]{Bertoin1996}.
 Applying the integration by parts to (\ref{ScaleFunction}), we have
 \beqlb\label{LaplaceDW}
 \int_0^\infty e^{-\lambda x} W'(x) dx = \int_0^\infty \lambda e^{-\lambda x} W(x) dx = \frac{1}{b+c\lambda^{\alpha}},   \quad \lambda > 0.
 \eeqlb
 The  Laplace transform of {\it Mittag-Leffler function}\footnote{The  Mittag-Leffler function $E_{\alpha,\alpha}$ on $\mathbb{R}_+$ is defined by
 $$	E_{\alpha,\alpha}(x):= \sum_{n=0}^\infty \frac{x^n}{\Gamma(\alpha(n+1))};$$
 see \cite{HauboldMathaiSaxena2011} for a precise definition and a survey of  its properties, e.g.,
 	\beqnn
 	\int_0^\infty  e^{-\lambda x} a x^{\alpha-1} E_{\alpha,\alpha}(-a\cdot x^\alpha) dx = \frac{a}{a+\lambda^{\alpha}} , \quad a,\lambda \geq 0 .
 	\eeqnn
 	\label{Mittag-Leffler}}
 yields that  $W'$ has the representation
 \beqnn
 W'(x)=   c^{-1}x^{\alpha-1}  \cdot E_{\alpha,\alpha}\big(-b/c\cdot x^\alpha\big),\quad x> 0.
 \eeqnn
 The smoothness of $E_{\alpha,\alpha}$ induces that $W$ is infinitely differentiable on $(0,\infty)$.
 When $b=0$, we have $E_{\alpha,\alpha}(0)= 1/\Gamma(\alpha)$ and
 \beqlb\label{ScaleF0}
 W(x)=  \frac{x^{\alpha}}{c\cdot \Gamma( \alpha +1)},  \quad W'(x)= \frac{x^{\alpha-1}}{c\cdot \Gamma(\alpha)}, \quad W''(x)= \frac{(\alpha-1)x^{\alpha-2}}{c\cdot \Gamma(\alpha)},\quad x>0.
 \eeqlb
 When $b>0$, the function $b W'$ is a \textit{Mittag-Leffler density function} and $1-bW(x)\to0$ as $x\to\infty$.
 The properties of Mittag-Leffler distribution/density function; see \cite{HauboldMathaiSaxena2011,Mainardi2014,MathaiHaubold2008}, yield that the scale function $W$ is  H\"older continuous with index $\alpha$. Moreover, we have as $x\to 0+$,
 \beqnn
 W(x) \sim \frac{ x^{\alpha}}{c\cdot\Gamma(\alpha+1)}, \quad
 W'(x)  \sim \frac{  x^{\alpha-1}}{c\cdot\Gamma(\alpha)},
 \quad W''(x)  \sim \frac{ (\alpha-1) x^{\alpha-2}}{c\cdot \Gamma(\alpha)} 
 \eeqnn
 and as $x\to \infty $, 
 \beqnn
 W(x)  \sim \frac{1}{b}- \frac{c \cdot x^{-\alpha}}{b^2  \Gamma(1-\alpha)},
 \quad   W'(x)  \sim \frac{c\alpha \cdot x^{-\alpha-1}}{b^2  \Gamma(1-\alpha)},\quad  W''(x)  \sim -\frac{c\alpha(\alpha+1) \cdot x^{-\alpha-2}}{b^2  \Gamma(1-\alpha)}.
 \eeqnn
 A direct consequence of these asymptotic properties and (\ref{ScaleF0}) is that uniformly in $x> 0$,
 \beqlb\label{UpperboundW}
 W(x)\leq C\cdot x^\alpha,\quad 
 \big| W'(x) \big|\leq C\cdot x^{\alpha-1}
 \quad \mbox{and}\quad  
 \big|W''(x)\big|\leq C\cdot x^{\alpha-2}.
 \eeqlb
 By the mean-value theorem, it is easy to identify that uniformly in $x,y> 0$,
 \beqlb\label{UpperboundDiffW}
 \big|\nabla_yW(x)\big| =\big|\Delta_y W(x-y)\big|  \leq C\cdot \big[ x^\alpha\wedge \big(|(x-y)^+|^{\alpha-1} y \big) \big].
 \eeqlb
 
 In addition to the scale function, we will also need a \textit{Sonine pair} $(K,L_K)$ on $(0,\infty)$  defined by
 \beqlb\label{Kernel}
 K(x) := \frac{ x^{\alpha-1}}{c\cdot \Gamma(\alpha)}
 \quad \mbox{and} \quad
 L_K(x):= \frac{c \cdot x^{-\alpha}}{\Gamma(1-\alpha)},  \quad x>0,
 \eeqlb
 which satisfies the \textit{Sonine equation}, i.e.,
 \beqlb\label{eqn.Sonine}
 K*L_K= L_K*K \equiv 1.
 \eeqlb
 In the theory of Volterra equations; see \cite{GripenbergLondenStaffans1990}, the function $L_K$ is also said to be the \textit{resolvent of the first kind} related to $K$ and vice versa.
 When $b>0$, a simple calculation shows that the function $bW'$ is the \textit{resolvent of the second kind} corresponding to $bK$, which is usually introduced by means of the resolvent equation
 \beqlb\label{KenelResolvent01}
 bW' =bK  -(bK)*(bW') .
 \eeqlb
 The function $bK$ is usually referred as the \textit{resolvent} associated to $bW'$.
 Convolving both sides of (\ref{KenelResolvent01}) by $L_K$ and then dividing them by $b$, we have
 \beqlb\label{KenelResolvent02}
 L_K*W'= W'*L_K =1-bW.
 \eeqlb
 Actually, this equality also holds when $b=0$, since $W'=K$ in this case; see (\ref{ScaleF0}) and (\ref{Kernel}).
 %
 %

 \subsection{Main results}
 We now formulate the main results for the local times of $\xi$ at the stopping time $\tau_\xi^L(\zeta)$ for a given value $\zeta>0$.
 For convenience, we write $L^\xi_\zeta$ for the process $ \{L_\xi(x,\tau_\xi^L(\zeta)): x\geq 0\}$ under  $\mathbf{P}(\,\cdot\,|\tau_\xi^L(\zeta)<\infty)$.
 Since $\tau_\xi^L(\zeta)<\infty$ a.s. when $b=0$, this conditional probability law turns to be  $\mathbf{P}$.
 When $b>0$, the stopping time $\tau_\xi^L(L_\xi(0,\infty)) $ is finite almost surely and equal to the last time that  $\xi$ hits $0$. In this case, we are also interested in the process
 \beqlb \label{eqn.Linf}
 L^\xi_\infty:= \big\{L_\xi(x,\infty):x\geq 0  \big\}\overset{\rm a.s.}=
 \big\{L_\xi \big(x,\tau_\xi^L(L_\xi(0,\infty)) \big):x\geq 0 \big\} ,
 \eeqlb
 under $\mathbf{P}$.
 Let $\varrho$ be an exponential random variable with mean $u^0(0)$, independent of  $N_\alpha(ds,dy,dz)$ and $N_0(dy,dz)$.
 Our first main theorem establishes SVEs for $L^\xi_\zeta$ and $L^\xi_\infty$.
 
 \begin{theorem} \label{MainThm.01}
 	We have the following:
 	\begin{enumerate}
 		\item[(1)] For each $\zeta\geq 0$, the process $L^\xi_\zeta$ is a weak solution of (\ref{MainThm.SVE}).
 		
 		\item[(2)] If $b>0$, the process $L_\infty^\xi$ is a weak solution of  (\ref{MainThm.SVE}) with $\zeta=\varrho$.
 		\vspace{3pt}
 		
 		\item[(3)] The weak uniqueness of non-negative solutions holds for (\ref{MainThm.SVE}).
 		
 	\end{enumerate}
 \end{theorem}

 \begin{remark}\label{Remark.ItoI}
 	By the change of variables and Proposition~\ref{WeakConvergenceProp01} with $p=2$, there exists a constant $C>0$ such that for any $x\geq 0$,
 	\beqnn
 	\int_0^x ds\int_0^\infty  \big| \nabla_yW(x-s)\big|^2 \nu_\alpha(dy)
 	\ar=\ar \int_0^x ds\int_0^\infty  \big| \nabla_yW(s)\big|^2 \nu_\alpha(dy) \leq C \cdot x^\alpha.
 	\eeqnn
 	Taking expectations on both sides of (\ref{MainThm.SVE}) and then using Fubini's theorem along with (\ref{KenelResolvent02}), we have
 	\beqlb\label{BoundFirstMom}
 	\mathbf{E}\big[ L^\xi_\zeta(x) \big]
 	\ar=\ar \mathbf{E}\Big[   \int_0^\infty  \int_0^\zeta \nabla_y W(x)   N_0( dy,dz) \Big]\cr
 	\ar=\ar \zeta \int_0^\infty  \nabla_y W(x)\bar{\nu}_\alpha(y) dy\cr
 	\ar=\ar \zeta  \int_0^\infty \bar{\nu}_\alpha(y) dy\int_{(x-y)^+}^x W'(s)ds \cr
 	\ar\ar\cr
 	\ar=\ar \zeta\cdot W'*L_K(x)\cr
 	\ar\ar\cr
 	\ar=\ar \zeta \big(1-bW(x)\big) \leq \zeta ,  
 	\eeqlb
 	for any $x\geq 0$. 
 	The SVI in (\ref{MainThm.SVE}) has finite quadratic variation and is well defined as an It\^o integral; see \cite[p.59-63]{IkedaWatanabe1989}.
 \end{remark}
 
 \begin{remark}
 	By the exponential formula for PRMs; see \cite[p.8]{Bertoin1996}, we have for any $\lambda \geq 0$,
 	\beqnn
 	\lefteqn{\mathbf{E}\Big[ \exp\Big\{-\lambda  \int_0^\infty  \int_0^\zeta \nabla_y W(x)   N_0( dy,dz) \Big\} \Big] }\ar\ar\cr
 	\ar=\ar \exp\Big\{ -\zeta \int_0^\infty (1-e^{-\lambda \nabla_y W(x)}) \bar{\nu}_\alpha(y)dy \Big\}.
 	\eeqnn 
 	From (\ref{BoundFirstMom}) and the fact that $\nabla_y W(x) \to 0$ uniformly in $y$ as $x\to 0$, we have
 	\beqnn
 	\lefteqn{\mathbf{E}\Big[ \exp\Big\{-\lambda  \int_0^\infty  \int_0^\zeta \nabla_y W(x)   N_0( dy,dz) \Big\} \Big]}\ar\ar\cr
 	\ar\sim\ar  \exp\Big\{ -\zeta  \lambda   \int_0^\infty  \nabla_y W(x)\bar{\nu}_\alpha(y) dy \Big\}
 	\to e^{-\zeta  \lambda }.
 	\eeqnn
 	Thus the first term on the right side of (\ref{MainThm.SVE}) converges almost surely to $\zeta$  as $x\to 0+$.
 	We make the convention that it is equal to $\zeta$ a.s. when $x=0$, which is consistent with the fact that $L^\xi_\zeta(0)\overset{\rm a.s.}=\zeta$.
 \end{remark}

 \begin{remark}
 	By (\ref{BoundFirstMom}), the SVE  (\ref{MainThm.SVE})  can be written as
 	\beqlb\label{Equ.SVE}
 	L^\xi_\zeta(x)\ar=\ar   \zeta \big(  1-bW(x) \big) + \int_0^\infty  \int_0^\zeta \nabla_y W(x) \widetilde{N}_0(dy,dz)\cr
 	\ar\ar  + \int_0^x \int_0^\infty \int_0^{L^\xi_\zeta(s)}  \nabla_y W(x-s)   \widetilde{N}_\alpha(ds,dy,dz), \quad x\geq 0,
 	\eeqlb
 	where  $\widetilde{N}_0(dy,dz):= N_0(dy,dz) - \bar{\nu}_\alpha(y)\cdot dy\cdot dz $. Here the first term on the right side of this equality represents the average local time at level $x$.
 	The second term can be interpreted as the perturbations caused by jumps up-crossing $0$;
 	the third term can be translated into the perturbations caused by jumps with initial positions above $0$  but below $x$.
 	More precisely, the convolution kernel $\nabla_yW(x-s)$ describes the impact of a jump with initial position $s$ and size $y$ on the local time at level $x$.
 	Notice that $\nabla_yW(x-s)$ increases when $x\in[s,s+y]$ and decreases as $x\to\infty$.
 	It would be sensible to consider the jump size of each jump as its life-length/residual-life during which it perturbs the local times directly.
 	This interpretation is consistent with the genealogical interpretations in \cite{FormanPalRizzoloWinkel2018,LambertBravo2018}.
 \end{remark}

 \begin{remark}
 	Because of the delayed and smooth relaxation of its perturbations, the PRM $N_\alpha(ds,dy,dz)$ fails to make solutions of (\ref{MainThm.SVE}) jump.
 	This phenomena cannot be observed in Ito's SDEs driven by PRM, since the PRM releases its perturbations instantaneously that give raise to jumps in the solutions.
 	Consequently, the continuity of driving noises is a necessary condition for the continuity of solutions of Ito's SDEs;  see  \cite[Chapter~III-IV]{IkedaWatanabe1989} and \cite[Chapter~II-V]{Protter2005}.
 \end{remark}
 
 \begin{remark}
 	It is necessary to specify that the SVE (\ref{MainThm.SVE}) is beyond the scope of the existing literature, e.g. \cite{Jaber2021,JaberCuchieroLarssonPulido2021,JaberLarssonPulido2019,ElEuchFukasawaRosenbaum2018,JaissonRosenbaum2016,PardouxProtter1990,Protter1985}.
 	More precisely, all SVEs studied in these literature are driven by finite-dimensional semimartingale and always can be written as
 	\beqlb\label{FDSVE}
 	X(t)= H(t)+ \int_0^t\mathbf{K}(t,s, X_s) d Z(t), \quad t\geq 0,
 	\eeqlb
 	where $H$ is a given function, $\mathbf{K}$ is a $d\times k$ matrix-valued convolution kernel on $\mathbb{R}_+^2\times \mathbb{R}$ and $Z$ is a $k$-dimensional It\^o's semimartingale whose differential characteristics are functions of $X$.
 	Differently, the SVI in (\ref{MainThm.SVE}) is driven by an infinite-dimensional martingale; see \cite{KurtzProtter1996} and  Appendix \ref{AppendixB}.
 	Since the impact of time $t$ on the convolution kernel $\nabla_yW(t)$ is tightly intertwined with that of mark $y$, one cannot write (\ref{MainThm.SVE}) into the form of (\ref{FDSVE}).
 	Consequently, it is difficult to prove the existence of solutions of  (\ref{MainThm.SVE}) by using the approximation method used in \cite{Jaber2021, JaberLarssonPulido2019} or the martingale problem theory developed in \cite{JaberCuchieroLarssonPulido2021}.
 \end{remark}
 
 \begin{remark}
 	Ito's SDEs with non-negative solutions have been widely studied in \cite{Bass2004,BertoinLeGall2006,DawsonLi2006,DawsonLi2012,FuLi2010} under two key conditions: (i) when solutions hit $0$, the diffusion vanishes and the drift turns to be non-negative; (ii) solutions cannot jump into the negative half-line.
 	In particular, it is the strong Markovianity that turns the state $0$ to be a tripper or a reflecting boundary, which results in the existence of non-negative solutions.
 	However, the convolution kernel in (\ref{MainThm.SVE}) results in the lack of (strong) Markovianity of the solutions and makes the standard stopping time method fail to prove the existence of non-negative solutions.
 	Fortunately, thanks to Theorem~\ref{MainThm.01}, the existence of non-negative solutions  of (\ref{MainThm.SVE})  follows directly from the non-negativity of $L^\xi_\zeta$.
 \end{remark}
 
 \begin{remark}\label{Remark.AbsorbS}
 	The point $0$ is an absorbing state\,\footnote{Although the two terminologies {\it absorbing state} and \textit{polarity} are initially introduced for Markov processes, it is sensible to use them to describe the analogous phenomena in other stochastic processes. Precisely, a state in a process is said to be an absorbing state if once it is entered, it is impossible to leave. A set is said to be a {\it polar set} for a process if it cannot be entered in finite time.} for the process $L^\xi_\zeta$ (and also $L^\xi_\infty$), i.e., once it hits $0$, it will stay at $0$ forever.
 	Indeed, the equivalence (\ref{HightMaximum}) shows that conditioned on  $\tau_\xi^L(\zeta)<\infty$,
 	\beqnn
 	\tau_0:=\inf\big\{x\geq 0: L^\xi_\zeta(x) =0\big\} <\infty,\  a.s.
 	\quad \mbox{and}\quad
 	L^\xi_\zeta(\tau_0+x) \overset{\rm a.s.}=0 ,\quad x\geq 0.
 	\eeqnn
 	Usually, the lack of Markovianity makes it difficult to obtain this property from the SVE (\ref{MainThm.SVE}).
 	Even for the SVE (\ref{FDSVE}), the absorbing states and polarity are also unclear up to now.
 \end{remark}

 %
 
 
 The SVE (\ref{MainThm.SVE}) makes it possible to study the local times of $\xi$ by using tools and methods from stochastic analysis, e.g., stochastic integral inequalities, stochastic Fubini theorem and martingale problem theory.
 To illustrate this, the next main theorem proves the H\"older continuity of $L^\xi_\zeta$ by using the Kolmogorov continuity theorem and also provides a uniform upper bound for all moments of the H\"older coefficients by using the Garsia-Rodemich-Rumsey inequality.
 For $\kappa\in (0,1]$ and $x> 0$,  the \textit{$\kappa$-H\"older coefficient} of a H\"older continuous function $f$ on $[0,x]$ is defined by
 \beqnn
 \big\|f\big\|_{C^{0,\kappa}_x} :=\sup_{0\leq y<z\leq x} \frac{|f(y)-f(z)|}{|y-z|^\kappa}.
 \eeqnn

 \begin{theorem}[H\"older continuity] \label{Thm.Regularity}
 	For each $\zeta \geq 0$, we have the following:
 	\begin{enumerate}
 		\item[(1)] The process $L^\xi_\zeta $ is H\"older-continuous of any order strictly less than $\alpha/2$.
 		
 		\item[(2)] For each $\kappa\in(0,\alpha/2)$ and $p\geq 0$,  there exists a constant $C>0$ such that for any $x\geq 0$,
 		\beqnn
 		\mathbf{E}\Big[\big\|L^\xi_\zeta\big\|_{C^{0,\kappa}_x}^p\Big]\leq C\cdot (1+x)^{p(\alpha-\kappa)}.
 		\eeqnn
 	\end{enumerate}
 \end{theorem}
 
 As a direct consequence of this theorem, we can establish a maximal inequality for $L^\xi_\zeta$.
 In detail, for each $\kappa\in(0,\alpha/2)$ and $p\geq 0$,  we have uniformly in $x\geq 0$,
 \beqnn
 \sup_{y\in[0,x]}\big|L^\xi_\zeta(y)-L^\xi_\zeta(0)\big|^p \leq \big\|L_\zeta^\xi\big\|_{C_x^{0,\kappa}}^p\cdot x^{p\kappa},\quad a.s.
 \eeqnn
 By the power mean inequality\footnote{For any $p>0$, $k\geq 2$ and $x_1,\cdots,x_k\in\mathbb{R}$, we have $|x_1+\cdots +x_k|^p\leq k^{(p-1)^+} (|x_1|^p+ \cdots +|x_k|^p)$.  \label{Footnote.PowerMI}} 
 and Theorem~\ref{Thm.Regularity},
 \beqnn
 \mathbf{E}\Big[\sup_{y\in[0,x]}\big|L^\xi_\zeta(y)\big|^p \Big]
 \ar\leq\ar C\cdot \zeta^p +  C\cdot \mathbf{E}\Big[\sup_{y\in[0,x]}\big|L^\xi_\zeta(y)-L^\xi_\zeta(0)\big|^p\Big]\cr
 \ar\leq\ar C\cdot \zeta^p +  C\cdot \mathbf{E}\Big[\big\|L_\zeta^\xi\big\|_{C_x^{0,\kappa}}^p\Big] \cdot x^{p\kappa},
 \eeqnn
 for some constant $C$ depending only on $p$. This yields the following corollary immediately.

 \begin{corollary}[Maximal inequality]\label{Cor.Maxi}
 	For each $\zeta>0$ and $p\geq 0$, there exists a constant $C>0$ such that for any $x\geq 0$,
 	\beqnn
 	\mathbf{E}\Big[   \sup_{y\in[0,x]}  \big| L^\xi_\zeta(y) \big|^p   \Big]
 	\leq C\cdot (1+x)^{p\alpha}.
 	\eeqnn
 	
 \end{corollary}
 
 \begin{remark}
 	By applying the discretization technique to the occupation measures of stable L\'evy processes, Boylan \cite{Boylan1964} proved the uniform H\"older regularity of their local times in both the time and spatial directions.
 	Later, Barlow \cite{Barlow1988} constructed the explicit modulus of continuity and gave the optimal H\"older exponent.
 	Recently, Forman et al. \cite{FormanPalRizzoloWinkel2018} proved the finiteness of all moments of the H\"older coefficient in the case of $b=0$.
 	The novelties of our results are in the uniform upper bound for all moments of the H\"older coefficient and the maximal inequality for the local times  in the spatial direction.
 \end{remark}
 
 Because of the Markovanity and martingale property, Brownian local times are tractable and their Laplace transform can be written as an exponential affine function of the initial state, in terms of the solution of a Riccati equation.
 By contrast, the lack of Markovanity and martingale property makes the local times of stable processes quite intractable.
 As another example that illustrates the strength of the SVE (\ref{MainThm.SVE}), the third main theorem establishes an explicit representation of Laplace functionals of $L^\xi_\zeta$.
 Let $\mathcal{A}_\infty$ be the space of all functions $f$ on $(0,\infty)$ satisfying that  $$\sup_{x\in(0,T]}x^{1-\alpha} |f(x)|<\infty$$ for any $T>0$.
 A  continuous function $v_\lambda^g$ on $(0,\infty)$ is said to be a \textit{$\mathcal{A}$-global solution} of the nonlinear-VIE (\ref{MainThm.Volterra}) if $v_\lambda^g \in \mathcal{A}_\infty$ and satisfies (\ref{MainThm.Volterra}) on $(0,\infty)$.
 
 \begin{theorem}[Laplace functionals] \label{MainThm.03}
 	For each $ \lambda\geq 0 $ and $g \in L^\infty(\mathbb{R}_+;\mathbb{R}_+)$, we have for $x\geq 0$,
 	\beqlb\label{AffineRepLap}
 	\lefteqn{ \mathbf{E}\Big[\exp\big\{- \lambda\cdot L^\xi_\zeta(x) -  g * L^\xi_\zeta(x) \big\} \Big] }\ar\ar\cr
 	\ar=\ar  \exp \Big\{- \zeta \int_0^\infty \Big(1-\exp\Big\{- \int_{(x-y)^+}^x   v_\lambda^g(s) ds\Big\}\Big)\bar\nu_\alpha(y)dy \Big\},
 	\eeqlb
 	where $v_\lambda^g$ is the unique $\mathcal{A}$-global solution of (\ref{MainThm.Volterra}).
 	
 \end{theorem}
 \begin{remark}
 	By comparing (\ref{AffineRepLap}) with (2.1) in \cite{DawsonLi2006} or (2.2) in \cite{DuffieFilipovicSchachermayer2003}, we see that the process $L^\xi_\zeta$ enjoys the analogue of the affine property of affine Markov processes.
 	In addition, this property also has been observed in the rough Heston model that is defined by a SVE driven by Brownian motion; see \cite{ElEuchRosenbaum2019}.
 	Later, Abi Jaber et al. \cite{Jaber2021,JaberLarssonPulido2019} considered a class of SVEs of the form (\ref{FDSVE}) with solutions being affine and also named them {\rm affine Volterra processes}.
 	
 \end{remark}
 
 By using the resolvent equations (\ref{KenelResolvent01}) and (\ref{KenelResolvent02}), the last main result in this work provides equivalent representations of the SVE (\ref{MainThm.SVE}) and the nonlinear-VIE  (\ref{MainThm.Volterra}), which help
 a lot to clarify the similarities and differences between the process $L^\xi_\zeta$ and CB-processes; we refer to \cite{Li2020} for a review on CB-processes.

 \begin{theorem}\label{MainThm.FIR}
 	The SVE (\ref{MainThm.SVE}) is equivalent to
 	\beqlb\label{FractionalSVR}
 	L^\xi_\zeta(x)\ar=\ar \zeta- b \cdot K*L^\xi_\zeta(x) +   \int_0^\infty  \int_0^\zeta  \int_{(x-y)^+}^{x}K(r) dr \widetilde{N}_0( dy,dz) \cr
 	\ar\ar  + \int_0^x \int_0^\infty \int_0^{L^\xi_\zeta(s)} \int_{(x-s-y)^+}^{x-s}K(r) dr
 	\widetilde{N}_\alpha(ds,dy,dz),\quad x\geq 0,
 	\eeqlb
 	and the nonlinear-VIE (\ref{MainThm.Volterra}) is equivalent to
 	\beqlb\label{FractionalRiccatti}
 	v_\lambda^g(x) \ar=\ar \lambda K(x)+ g*K(x) - \big(b\cdot v_\lambda^g + \mathcal{V}_\alpha\circ v_\lambda^g \big)*K(x),\quad x> 0.
 	\eeqlb
 \end{theorem}
 
 By (\ref{Kernel}) and the Sonine equation (\ref{eqn.Sonine}), the nonlinear-VIE (\ref{FractionalRiccatti}) can be written into a fractional differential equation.
 Indeed, denote by $I^\alpha_c$ and $D^\alpha_c$ the \textit{Riemann-Liouville} fractional integral and derivative operators of order $\alpha$ modified by the constant $c>0$ respectively.
 They act on a measurable function $f$ according to
 \beqnn
 I^\alpha_c  f(x)\ar:=\ar \int_0^x\frac{(x-s)^{\alpha-1}}{c\cdot \Gamma(\alpha)}f(s)ds, \cr
 D^\alpha_c  f(x)\ar:=\ar \frac{d}{dx} \int_0^x\frac{c(x-s)^{-\alpha}}{\Gamma(1-\alpha)} f(s)ds,
 \eeqnn
 for $ x\geq 0$. 
 Convolving both sides of (\ref{FractionalRiccatti}) by $L_K$ and then using (\ref{eqn.Sonine}), we have
 \beqnn
 L_K* v_\lambda^g(x) \ar=\ar \lambda + \int_0^x g(s)ds  -\int_0^x \big(b\cdot v_\lambda^g + \mathcal{V}_\alpha\circ v_\lambda^g \big)(s)ds,\quad x\geq 0.
 \eeqnn
 Notice that $I^\alpha_c f =  K*f$ and $D^\alpha_c v_\lambda^g=   \frac{d}{dx} ( L_K* v_\lambda^g)$. A simple calculation yields the next corollary.
 \begin{corollary}
 	The nonlinear-VIE (\ref{MainThm.Volterra}) is equivalent to
 	\beqlb\label{FractionalLP}
 	D^\alpha_c v_\lambda^g  \ar=\ar   g -b\cdot v_\lambda^g - \mathcal{V}_\alpha\circ v_\lambda^g
 	\quad \mbox{with}\quad
 	I^{1-\alpha}_c v_\lambda^g(0) =   \lambda.
 	\eeqlb
 \end{corollary}
 
 \begin{remark}
 	By comparing (\ref{FractionalSVR}) with the It\^o's SDE for a CB-process\footnote{Here $a\in\mathbb{R}$, $\varsigma\geq 0$, $B_1(ds,dz)$ is a Gaussian white noise on $(0,\infty)^2$ with intensity $ds\cdot dz$,  $\widetilde{N}_1(ds,dy,dz)$ is a PRM on $(0,\infty)^3$ with intensity  $ds\cdot \nu_1(dy)\cdot dz$ and $\nu_1(dy)$ is a $\sigma$-finite measure on $(0,\infty)$ satisfying that $\int_0^\infty(1\wedge y^2)\nu_1(dy)<\infty$.}
 	\beqlb\label{CBProcess}\quad
 	Y_\zeta(t)
 	\ar=\ar \zeta -\int_0^t a Y_\zeta(s)ds+ \int_0^t \int_0^{Y_\zeta(s)}  \varsigma\, B_1(ds,dz) \cr
 	\ar\ar + \int_0^t \int_0^\infty \int_0^{Y_\zeta(s-)} y\,\widetilde{N}_1(ds,dy,dz),
 	\eeqlb
 	and then comparing (\ref{FractionalLP}) with the Riccati equation for the Laplace exponent of $Y_\zeta$
 	\beqnn
 	\frac{\partial}{\partial t}v_\lambda= -av_\lambda  - \frac{\varsigma^2}{2} |v_\lambda |^2 - \int_0^\infty \big(e^{-y v_\lambda }-1+ yv_\lambda \big) \nu_1(dy)\quad \mbox{with}\quad v_\lambda(0)=\lambda;
 	\eeqnn
 	see \cite[Section~3]{DawsonLi2012}, we see that the processes $L^\xi_\zeta$ not only owns an evolution mechanism similar to that of $Y_\zeta$, but also enjoys, formally at least, the striking analogue\footnote{For $\zeta_1,\zeta_2>0$, assume $L^\xi_{\zeta_1}$ and $L^\xi_{\zeta_2}$ are two independent solutions of (\ref{MainThm.SVE}) with $\zeta=\zeta_1$ and $\zeta=\zeta_2$ respectively, we have  $L^\xi_{\zeta_1}+L^\xi_{\zeta_2}$ is the unique weak non-negative solution of (\ref{MainThm.SVE}) with $\zeta=\zeta_1+\zeta_2$.} of the branching property.
 	Moreover, it is same to $Y_\zeta$ that the point $0$ is an absorbing state for $L^\xi_\zeta$; see Remark~\ref{Remark.AbsorbS}.
 	In conclusion, the SVE (\ref{MainThm.SVE}) defines a fully novel non-Markov CB-process whose degree of H\"older regularity is less than that of Feller branching diffusion; see Theorem~\ref{Thm.Regularity}.
 	Drawing from \cite{ElEuchRosenbaum2019,JaissonRosenbaum2016}, we may refer $L^\xi_\zeta$  as a {\rm rough CB-process}.  
 \end{remark}

 \begin{remark}
 	The CB-process defined by (\ref{CBProcess}) can be reconstructed as the scaling limit of discrete Galton-Watson branching processes.
 	The life-length of individuals uniformly tends to $0$ after rescaling and jumps in (\ref{CBProcess}) result from the simultaneous births of infinite offsprings; readers may refer to \cite{Li2020} for details.
 	In contrast, the variable $y$ in (\ref{MainThm.SVE}) is positive and can be interpreted as the life-length of each individual during which it gives birth to its children randomly.
 	This slows down the extinction of the population. In detail, the survival probability of the CB-process $Y_\zeta$ deceases to $0$ at an exponential rate when $a>0$, i.e., $\mathbf{P}(Y_\zeta(t)>0)\sim C \cdot e^{-at}$ as $t\to\infty$.
 	While, the survival probability of the rough CB-process $L^\xi_\zeta$ decreases to $0$ at a power rate when $b>0$; see the supplementary material \cite{Xu2021c}.
 \end{remark}
 
 \section{Local times of compound Poisson processes}
 \label{CompoundPoissonP}

 In this section, we first introduce some properties of local times of compound Poisson processes with negative drift and then establish a SVE for them.
 The proofs are elementary and will be merely sketched.
 
 For two constants $\beta,\gamma>0$ and a probability law $\nu(dx)$ on $(0,\infty)$ with finite mean $m_\nu$,
 let $Y:=\{Y(t):t\geq 0 \}$ be a compound Poisson process  on $(\Omega, \mathscr{F}, \mathscr{F}_t, \mathbf{P})$ with a drift $-\beta $, arrival rate $\gamma$ and jump-size distribution $\nu$.
 It is a L\'evy process with bounded variation and L\'evy measure $\gamma\cdot \nu(dx)$.
 Its Laplace exponent is of the form
 \beqnn
 \varphi (\lambda):= \beta \lambda +  \gamma \int_0^\infty (e^{-\lambda x}-1)   \nu (dx) = \lambda \Big(\beta -\gamma\int_0^\infty e^{-\lambda x}\bar\nu (x) dx \Big) ,\quad \lambda \geq 0,
 \eeqnn
 where  $\bar\nu(x):=\nu([x,\infty))$ is the tail distribution of $\nu$.
 The function $\varphi $ is zero at zero and tends to $\infty$ at infinity. Moreover, it is infinitely differentiable and strictly convex on $(0,\infty)$.
 In particular, $\varphi'(0)= -\mathbf{E}[Y(1)]= \beta - \gamma\cdot m_\nu$ and hence $\varphi$ is increasing on $[0,\infty)$ if $\varphi'(0)\geq 0$.
 The process $Y$ drifts to $-\infty$, $\infty$ or is recurrent according as $\varphi'(0)>0$, $<0$ or $=0$.
 Denote by $\tau^+_{Y}$ the first passage time of $Y$ in $[0,\infty)$, i.e.,
 $$\tau^+_{Y}:=\inf\big\{t> 0:Y(t)\geq 0\big\}.$$
 Actually, the process $Y$ always moves from the negative half line into the positive half line by jumping, i.e., $Y(\tau^+_{Y}-)<0$ and $Y(\tau^+_{Y})>0$ a.s.
 The next proposition comes from Theorem~17(ii) in \cite[p.204]{Bertoin1996}. Let $\nu^*(dx) $ be the \textit{size-biased distribution} of $\nu$ given by
 \beqnn
 \nu^*(dx) := \mathbf{1}_{\{ x> 0\}} \cdot \frac{ \bar\nu(x)}{m_\nu} \cdot dx.
 \eeqnn
 
 \begin{proposition}
 	If $\varphi'(0)\geq 0$, then $Y(\tau^+_{Y})$ under $\mathbf{P}(\,\cdot\,|\tau^+_{Y}<\infty)$ is distributed as $\nu^*$.
 \end{proposition}
 
 Let  $L_{Y}:=\{L_{Y}(x,t):x\in \mathbb{R}, t\geq 0 \}$ be the local times of $Y$  satisfying the occupation density formula (\ref{OccupationDensityF}).
 The local time $L_{Y}(x,\infty)$ is infinite almost surely for some and hence all $x\in\mathbb{R}$ if and only if $\varphi'(0)=0$.
 Denote by $\tau^L_{Y}:=\{\tau^L_{Y}(\zeta): \zeta\geq 0\}$ the right-inverse local time at level $0$.
 When $\beta =1$, the local time $L_{Y}(x,t)$ equals to the times that $Y$ hits $x$ in the time interval $(0,t]$, i.e.,
 \beqnn
 L_{Y} \overset{\rm a.s.}= \big\{\# \{s\in(0, t]: Y(s)=x \}: x\in\mathbb{R},t\geq 0 \big\}
 \eeqnn
 and $\tau^L_Y$ only jumps at positive integer points.
 The next proposition follows directly from (\ref{TimeSpatialChange01})-(\ref{TimeSpatialChange02}).
 \begin{proposition}\label{Dis.Cross0}
 	If $\beta>0$, the process $\tau^L_{Y}$ only jumps at points $\{ k/\beta: k=1,2,\cdots \}$ and
 	\beqnn
 	L_{Y} \overset{\rm a.s.}= \big\{\beta^{-1}\cdot \# \{s\in(0, t]: Y(s)=x \}: x\in\mathbb{R},t\geq 0 \big\}.
 	\eeqnn
 	
 \end{proposition}
 
 We write $ L^{Y}_{k/\beta}$ and $ L^{Y,*}_{k/\beta}$ for the process $  \{ L_{Y}(x,\tau^L_{Y}(k/\beta)):  x\geq 0\}$  under $\mathbf{P}(\,\cdot\,|\tau^L_{Y}(k/\beta)<\infty)$ and $\mathbf{P}_{\nu^*}(\,\cdot\,|\tau^L_{Y}(k/\beta)<\infty)$ respectively.
 Notice that sample paths on $[0,\tau^+_Y]$ make no contribution to the local times $  \{ L_{Y}(x,\tau^L_{Y}(k/\beta)):  x\geq 0, k\geq 1\}$.
 The next proposition can be proved immediately by using Proposition~\ref{Dis.Cross0}, the strong Markov property and independent, stationary increments of $Y$.
 
 \begin{proposition}\label{Decomposition}
 	If $\varphi'(0)\geq 0$,  for any $ k\geq 1$, we have
 	$$
 	L^{Y}_{k/\beta} \overset{\rm d}= L^{Y,*}_{k/\beta}  \overset{\rm d}=\sum_{i=1}^k L^{Y,*}_{1/\beta,i}, $$
 	where $\{L^{Y,*}_{1/\beta,i}:i= 1,2,\cdots\}$ is a sequence of  i.i.d. copies of  $  L^{Y,*}_{1/\beta}$.
 \end{proposition}

 \subsection{Branching representation}
 Lambert \cite{Lambert2010} established a one-to-one correspondence between the local times of compound Poisson processes with drift $-1$ and  homogeneous, binary CMJ-processes.
 More precisely, he observed that the jumping contour process of a homogeneous, binary CMJ-tree starting from one ancestor is a compound Poisson process with a drift $-1$; conversely, the local times on $\mathbb{R}_+$ of a compound Poisson process with drift $-1$ stopped upon hitting $0$ are equal in distribution to a homogeneous, binary  CMJ-process starting from one ancestor.
 In this section, we  extend this connection slightly to the case of compound Poisson processes with arbitrary negative drift.
 
 Recall the triplet $(\beta,\gamma, \nu)$ with $\beta > 0$. We consider a binary CMJ-process on $(\Omega,\mathscr{F},\mathscr{F}_t, \mathbf{P})$ defined by the following three properties:
 \begin{enumerate}
 	\item[\namedlabel{P1}{(P1)}]  There are $k\in\mathbb{Z}_+$ ancestors at time $0$ with residual life distributed as $\nu^*$. \vspace{5pt}
 	
 	\item[\namedlabel{P2}{(P2)}]  Offsprings have a common life-length distribution $\nu$. \vspace{5pt}
 	
 	\item[\namedlabel{P3}{(P3)}]  Each individual gives birth to its children according to a Poisson process with rate $\gamma/\beta$.
 \end{enumerate}
 
 Denote by $\mathcal{I}_k$ the collection of all individuals in the population.
 Each individual $x \in \mathcal{I}_k$ is endowed with a pair $(\sigma_x,\ell_x)$ that represents its birth time and life length.
 For convention, if $x$ is an ancestor, we assume $\sigma_x=0$ and $\ell_x$ equals to its residual life.
 For $t\geq 0$, let $\mathcal{I}_k(t):= \{x\in \mathcal{I}_k:0\leq t-\sigma_x < \ell_x \}$ and $\#\mathcal{I}_k(t)$ be the collection and number of all individuals alive at time $t$ respectively.
 
 \begin{lemma}\label{CMJRepresentation}
 	If $\varphi'(0)\geq 0$, then the process $ L^{Y}_{k/\beta}$ is equal in distribution to $  \beta^{-1}\cdot \#\mathcal{I}_k $.
 \end{lemma}
 \proof Let $Y_{\beta}:=\{ Y(t/\beta):t\geq 0 \}$, which is a compound Poisson process with a triplet $(1,\gamma/\beta, \nu)$, local times $L_{Y_\beta}$ and the right-inverse local time $\tau^L_{Y_\beta}$ at level $0$.
 Let $L^{Y_{\beta}}_k$ be the process $\{ L_{Y_\beta}(x,\tau^L_{Y_\beta}(k)):  x\geq 0\}$ under $\mathbf{P}(\,\cdot\,|\tau^L_{Y_{\beta}}(k)<\infty)$.
 Theorem 3.2 in \cite{LambertSimatosZwart2013} asserts that $L^{Y_{\beta}}_k \overset{\rm d}=  \#\mathcal{I}_k $.
 By (\ref{TimeSpatialChange01})-(\ref{TimeSpatialChange02}) and Proposition~\ref{Decomposition}, we have $ L^{Y}_{k/\beta} \overset{\rm d}= \beta^{-1}\cdot L^{Y_{\beta}}_k \overset{\rm d}= \beta^{-1}\cdot \#\mathcal{I}_k $.
 \qed

 \subsection{Stochastic Volterra representation}
 In this section, we establish a SVE for the process $L^Y_{k/\beta}$ by linking the preceding CMJ-process to a MHP; see Appendix~\ref{Appendix.HP}.
 It also can be obtained by using the connection between general branching processes and multivariate MHPs established in \cite[Section~4]{Xu2021}.
 
 As a preparation, we first clarify the genealogy of the CMJ-process.
 Denote by $\ell_{0,j}$ the residual life of the $j$-th ancestor at time $0$, which is distributed as $\nu^*$.
 We sort all offsprings by their birth times.
 Associated to the sequence $\{(\sigma_i,\ell_i):i=1,2,\cdots\}$ we define a $(\mathscr{F}_t)$-random point measure $N_\nu(ds,dy)$ on $(0,\infty)^2$ by
 \beqnn
 N_\nu(ds,dy) := \sum_{i=1}^\infty \mathbf{1}_{\{\sigma_i\in ds,\ell_i\in dy\}}.
 \eeqnn
 By the branching mechanism, a child will be born at time $t$ at the rate $\gamma/\beta\cdot \#\mathcal{I}_k(t-)$.
 Thus the random point measure $N_\nu(ds,dy)$ has a $(\mathscr{F}_t)$-intensity $\gamma/\beta \cdot \#\mathcal{I}_k(s-)\cdot ds \cdot \nu(dy) $.
 Notice that $\#\mathcal{I}_k$ has the representation
 \beqnn
 \#\mathcal{I}_k(t) = \sum_{j=1}^{k} \mathbf{1}_{\{ \ell_{0,j}>t \}} +   \sum_{\sigma_i \leq t}\mathbf{1}_{\{ \ell_{i}>t-\sigma_i \}},\quad t\geq 0.
 \eeqnn
 Here the two sums on the right side of this equality represent the numbers of ancestors and offsprings alive at time $t$ respectively.
 Similarly as in  Appendix~\ref{Appendix.HP}, we can write the foregoing equation into
 \beqlb\label{PreHR}
 \#\mathcal{I}_k(t)  = \sum_{j=1}^{k} \mathbf{1}_{\{ \ell_{0,j}>t \}} + \int_0^t \int_0^\infty  \mathbf{1}_{\{y>t-s\}} N_\nu(ds,dy),\quad t\geq 0.
 \eeqlb
 Hence  $N_\nu(ds,du)$  is a MHP on $(0,\infty)^2$.
 Before applying Proposition~\ref{SVR} to establish a SVE for the process $L^Y_{k/\beta}$, we need to introduce two important quantities related to the triplet $(\beta,\gamma,\nu)$.
 Let $R_\nu$ be the \textit{resolvent} associated to the function $\gamma/\beta \cdot \bar\nu$ defined by the unique solution of
 \beqlb\label{Resolvent.01}
 R_\nu =  \frac{\gamma}{\beta} \cdot \bar\nu +  \frac{\gamma}{\beta}\cdot \bar\nu * R_\nu .
 \eeqlb
 Actually, the function $R_\nu$ can be interpreted as the mean  reproduction rate of an individual and its descents.
 In addition, we introduce a two-parameter function on $\mathbb{R}_+^2$
 \beqlb\label{Resolvent.02}
 R(t,y)\ar=\ar \mathbf{1}_{\{y>t\}} + \int_0^t R_\nu (t-s)\cdot \mathbf{1}_{\{y>s\}}ds,\quad t,y\geq 0, \label{Resolvent02}
 \eeqlb
 to describe the mean reproduction rate of an individual with life-length $y$ and its descents; we refer to \cite[Section~2]{Xu2021} for more detailed explanations.

 \begin{theorem}\label{Thm.SVR}
 	If $\varphi'(0)\geq 0$, we have for $ k\geq 1$, the process $L^Y_{k/\beta}$ is equal in distribution to the unique solution of the SVE
 	\beqlb\label{HawkesRepresentation}
 	Z_k(t)\ar=\ar \frac{1}{\beta}\sum_{j=1}^{k} \mathbf{1}_{\{ \ell_{0,j}> t \}} + \frac{1}{\beta}\int_0^t R_\Lambda  (t-s) \sum_{j=1}^{k} \mathbf{1}_{\{ \ell_{0,j}> s \}} ds\cr
 	\ar\ar + \int_0^t \int_0^\infty \int_0^{Z_k (s-)} \frac{1}{\beta}R  (t-s,y ) \widetilde{N} (ds,dy,dz), \quad t\geq 0,
 	\eeqlb
 	where $\widetilde{N}(ds,dy,dz)$ is a compensated PRM on $(0,\infty)^3$ with intensity $ \gamma \cdot ds\cdot \nu(dy)\cdot dz$. 
 \end{theorem}
 \proof From (\ref{PreHR}) and Proposition~\ref{SVR}, the population process $\#\mathcal{I}_k(\cdot)$ is the unique solution of
 \beqnn
 \#\mathcal{I}_k(t) 
 \ar=\ar \sum_{j=1}^{k} \mathbf{1}_{\{ \ell_{0,j}> t \}} + \int_0^t R_\nu (t-s) \sum_{j=1}^{k} \mathbf{1}_{\{ \ell_{0,j}>   s \}} ds \cr
 \ar\ar + \int_0^t \int_0^\infty \int_0^{\#\mathcal{I}_k(s-) } R (t-s,y) \widetilde{N}_\nu (ds,dy,dz),
 \eeqnn
 where $\widetilde{N}_\nu(ds,dy,dz)$ is compensated PRM on $(0,\infty)^3$ with intensity $ \gamma/\beta \cdot ds \cdot \nu(dy)\cdot dz$.
 By the change of variables, the process $\beta^{-1}\cdot \#\mathcal{I}_k $ is the unique solution to (\ref{HawkesRepresentation}) with $\widetilde{N}(ds,dy,dz)=\widetilde{N}_\nu( ds,dy,\beta\cdot dz)$ and the desired result follows  from Lemma~\ref{CMJRepresentation}.
 \qed

 \section{Stochastic Volterra equation for $L^{\xi}_\zeta$}\label{ProofForThm.1}

 In this section, we prove the first two claims in Theorem~\ref{MainThm.01}.
 By the self-similarity of stable processes and their local times; see (\ref{TimeSpatialChange01})-(\ref{TimeSpatialChange02}), the next lemma shows that it suffices to prove them for a subclass of spectrally positive stable processes.
 It will help a lot to simplify our proofs.
 \begin{lemma}\label{lemma.01}
 	Let $c_0>0$. If Theorem~\ref{MainThm.01}(1) holds for any $\xi$ with $(\alpha,b,c)\in (0,1)\times [0,\infty)\times \{ c_0 \}$, then it holds for all $\xi$ with $(\alpha,b,c)\in (0,1)\times [0,\infty)\times (0,\infty)$.
 \end{lemma}
 \proof For any $\xi$ with parameters $(\alpha,b,c)\in (0,1)\times [0,\infty)\times (0,\infty)$, its self-similarity induces that $ a\xi  $ with $a=(c_0/c)^{1/(1+\alpha)} $ is a spectrally positive stable process with parameters $(\alpha,ab,c_0)$  and L\'evy measure $\nu_\alpha(a^{-1}\cdot dy)=c_0/c\cdot \nu_\alpha(dy)$.
 For $\zeta >0$, by Theorem~\ref{MainThm.01}(1) the process $L^{a\xi}_\zeta$ is a weak solution to
 \beqnn
 L^{a\xi}_\zeta(x)\ar=\ar  \int_0^\infty  \int_0^\zeta \nabla_y W_a(x)N_{0,a}(dy,dz) \cr
 \ar\ar + \int_0^x \int_0^\infty \int_0^{L^{a\xi}_\zeta(s)} \nabla_y W_a(x-s)  \widetilde{N}_a(ds,dy,dz), \quad x\geq 0,
 \eeqnn
 where $W_a$ is the scale function of $a\xi$, $N_{0,a}(dy,dz) $ is PRM on $(0,\infty)^2$ with intensity $   c_0/c\cdot \bar\nu_\alpha(y)\cdot dy\cdot dz $  and $\widetilde{N}_a(ds,dy,dz)$ is a compensated PRM on $(0,\infty)^3$ with intensity $ c_0/c\cdot ds \cdot  \nu_\alpha(dy)\cdot dz $.
 By (\ref{ScaleFunction}) as well as the continuity of $W$ and $W_a$, we have $W(x)=aW_a(ax)$ for any $x\in\mathbb{R}$.
 Application of (\ref{TimeSpatialChange01})-(\ref{TimeSpatialChange02}) gives that $ L^{\xi}_\zeta \overset{\rm a.s.}= aL_{\zeta/a}^{a\xi}(a\cdot)$.
 Taking these back into the preceding SVE and then using the change of variables,
 \beqnn
 L^{\xi}_\zeta(x)
 \ar=\ar  \int_0^\infty  \int_0^\zeta \nabla_{y} W(x)N_{0,a}(a\cdot dy,a^{-1}\cdot dz) \cr
 \ar\ar +  \int_0^{x} \int_0^\infty \int_0^{ L^{\xi}_\zeta(s)} \nabla_y W(x-s) \widetilde{N}_a(a\cdot ds,a\cdot dy,a^{-1}\cdot dz),
 \quad x\geq 0.
 \eeqnn
 Notice that $\nu_\alpha(a\cdot dy)= a^{-\alpha-1} \nu_\alpha(dy)$; see (\ref{LevyMeasure}).
 The PRM $N_{0,a}(a\cdot dy,a^{-1}\cdot dz)$ has intensity $\bar\nu_\alpha(y)dydz$ and
 the compensated PRM $\widetilde{N}_a(a\cdot ds,a\cdot dy,a^{-1}\cdot dz)$ has intensity $ds\cdot  \nu_\alpha(dy)\cdot dz$. Thus $L^{\xi}_\zeta$ is a weak solution of (\ref{MainThm.SVE}) and Theorem~\ref{MainThm.01}(1) holds for $\xi$.
 \qed

 %
 \subsection{Compound Poisson approximation}\label{Sec.CP}
 Based on Lemma~\ref{lemma.01}, we start to prove Theorem~\ref{MainThm.01}(1) for any $\xi$ with
 \beqnn
 \alpha\in(0,1),\quad b\geq 0, \quad c \equiv \Gamma(1-\alpha),\quad \nu_\alpha(dy)=\alpha (\alpha+1) \frac{dy}{y^{\alpha+2}}, \quad 
 \bar\nu_\alpha(y)=\frac{\alpha}{y^{\alpha+1}}.
 \eeqnn
 Let us consider a sequence of  compound Poisson processes with negative drift and positive jumps, whose local times give a good approximation for the process $L^\xi_\zeta$.
 Denote by $\Lambda $ the Pareto-II distribution on $\mathbb{R}_+$ with location $0$ and shape $\alpha+1$, i.e.,
 \beqlb\label{ParetoDis}
 \Lambda (dx):= \frac{\alpha+1}{(1+x )^{\alpha+2}} dx
 \quad \mbox{with tail-distribution}\quad
 \bar{\Lambda}(x)= (1+x )^{-\alpha-1},\quad x\geq 0.
 \eeqlb
 For $n\geq 1$, let $\{Y^{(n)}(t):t\geq 0\}$ be a compound Poisson process with a drift $-1$, arrival rate $\gamma_n\in(0,\alpha)$, jump-size distribution $\Lambda$, Laplace exponent $\varphi^{(n)}$ and local times $L_{Y^{(n)}}$.
 We now consider the behavior of $Y^{(n)}$ and $L_{Y^{(n)}}$ at a large time scale under the following condition.
 \begin{condition}\label{Main.Condition}
 	Assume that $n^\alpha(1-\gamma_n/\alpha)\to b  $ as $n\to\infty$.
 \end{condition}
 
 A routine computation, along with this condition, induces that the rescaled Laplace exponent $ {\it\Phi}^{(n)}:=\{ n^{1+\alpha} \varphi^{(n)}(\lambda/n):\lambda \geq 0 \}$ converges locally uniformly to ${\it\Phi} $ on $\mathbb{R}_+$  as $n\to\infty$.
 By Corollary~4.3 in \cite[p.440]{JacodShiryaev2003}, the rescaled compound Poisson process
 $$  \xi^{(n)}:=\big\{n^{-1}\cdot Y^{(n)}(n^{1+\alpha}t):t\geq 0\big\} $$
 converges weakly to $\xi$ in $D([0,\infty);\mathbb{R})$ as $n\to\infty$.
 Let $L_{\xi^{(n)}}$ be the local times of $\xi^{(n)}$ satisfying (\ref{OccupationDensityF}) and $\tau^L_{\xi^{(n)}}$ the right-inverse local time at level $0$.
 Similarly, we also write $L^{\xi^{(n)}}_\zeta$ for the process $\{L_{ \xi^{(n)}}(x,\tau^L_{\xi^{(n)}}(\zeta)):x\geq 0\} $ under  $\mathbf{P}(\,\cdot\, |\tau^L_{\xi^{(n)}}(\zeta)<\infty)$.
 The following  convergence result comes from Theorem~2.4 in \cite{LambertSimatos2015}.

 \begin{lemma}\label{ConvergenceLocalTime}
 	For $\zeta>0$, we have $L^{\xi^{(n)}}_\zeta \to L^\xi_\zeta$ weakly in $D([0,\infty);\mathbb{R}_+)$ as $n\to\infty$.
 \end{lemma}
 
 For $n\geq 1$, let  $Z^{(n)}_{[n^\alpha\zeta]}$ be the unique solution of (\ref{HawkesRepresentation}) with $k=[n^\alpha\zeta]$, $\beta =1$, $\gamma=\gamma_n$ and $\nu=\Lambda$.
 By (\ref{TimeSpatialChange01})-(\ref{TimeSpatialChange02}) and Theorem~\ref{Thm.SVR},
 \beqlb\label{eqn.28}
 L^{\xi^{(n)}}_\zeta \overset{\rm a.s.}= \big\{ n^{-\alpha}\cdot L^{Y^{(n)}}_{[n^\alpha\zeta]}(nx): x\geq 0  \big\} \overset{\rm d}= \big\{ n^{-\alpha}\cdot Z^{(n)}_{[n^\alpha\zeta]}(nt): t\geq 0 \big\} =: X^{(n)}_\zeta
 \eeqlb
 and hence $X^{(n)}_\zeta\to L^\xi_\zeta$ weakly in $D([0,\infty);\mathbb{R}_+)$ as $n\to\infty$; see Lemma~\ref{ConvergenceLocalTime}.
 Clearly, Theorem~\ref{MainThm.01}(1) can be proved by characterizing the cluster points of  $\{X^{(n)}_\zeta \}_{n\geq 1}$ as weak solutions of (\ref{MainThm.SVE}).
 
 Let $R^{(n)}_\Lambda$ and $R^{(n)}$ be the resolvent associated to the tail-distribution $\bar\Lambda$ defined as in (\ref{Resolvent.01})-(\ref{Resolvent.02}) with $\beta =1$, $\gamma=\gamma_n$ and $\nu=\Lambda$, i.e., $R^{(n)}_\Lambda(t)= R^{(n)}(t,y) =0$ if $t\leq 0$ or $y\leq 0$ and  for $t,y\geq 0$,
 \beqlb
 R^{(n)}_\Lambda(t)\ar=\ar  \gamma_n \bar\Lambda(t) +  \gamma_n  \bar\Lambda * R_\Lambda^{(n)}(t), \label{Resolvent.n01}\\
 R^{(n)}(t,y) \ar=\ar \mathbf{1}_{\{y>t\}} + \int_0^t R_\Lambda^{(n)} (t-s)\cdot \mathbf{1}_{\{y>s\}}ds. \label{Resolvent.n02}
 \eeqlb
 By (\ref{eqn.28}), Theorem~\ref{Thm.SVR} and the change of variables, the process $X^{(n)}_\zeta$ satisfies the SVE
 \beqlb\label{Eqn.HR}
 X_{\zeta}^{(n)}(t)
 \ar=\ar \frac{1}{n^{\alpha}}\sum_{i=1}^{[n^\alpha \zeta ] } \mathbf{1}_{\{\ell_{0,i}>nt\}} + \int_0^{nt} R_\Lambda^{(n)}(nt-s) \cdot  \frac{1}{n^{\alpha}}\sum_{i=1}^{[n^\alpha \zeta ]} \mathbf{1}_{\{\ell_{0,i}>s\}} ds\cr
 \ar\ar + \int_0^t \int_0^\infty \int_0^{X_{\zeta}^{(n)}(s-)} \frac{1}{n^{\alpha}} R^{(n)}\big(n(t-s),ny\big) \widetilde{N}^{(n)}(n\cdot ds,n\cdot dy,n^{\alpha}\cdot dz), 
 \eeqlb
 where $\widetilde{N}^{(n)}(n\cdot ds,n\cdot dy,n^{\alpha}\cdot dz)$ is a compensated PRM on $(0,\infty)^3$ with intensity $n^{1+\alpha} \gamma_n \cdot ds\cdot \Lambda(n\cdot dy)\cdot dz$.
 Here $\ell_{0,i}$ is distributed as the size-biased distribution of $\Lambda$
 \beqlb\label{SizeBiasedLambda}
 \Lambda^*(dx):=  \frac{\alpha\cdot dx}{(1+x)^{\alpha+1}} 
 \quad \mbox{with tail-distribution}\quad
 \overline{\Lambda^*}(x):= (1+x)^{-\alpha},\quad
 x\geq 0.
 \eeqlb

 \subsection{Proofs for Theorem~\ref{MainThm.01}(1)-(2)}
 In order to make the proofs much easier to be understood, we offer an intuitive description on how to derive the convergence of the SVE (\ref{Eqn.HR}) to the SVE (\ref{MainThm.SVE}).
 The following auxiliary lemmas will be proved in Section~\ref{ProofAuxilliaryLemma}.

 We start the asymptotic analysis from the first two terms on the right side of (\ref{Eqn.HR}).
 Denote by $X_{\zeta,0}^{(n)}(t)$ their sum.
 Let $M(\mathbb{R}_+)$ be the space of all $\sigma$-finite measures on $\mathbb{R}_+$ endowed with the topology of weak convergence and $M_p(\mathbb{R}_+)\subset M(\mathbb{R}_+)$ the space of all Radon point measures on $\mathbb{R}_+$.
 Let $N^{(n)}_\zeta(dy)$ be a finite point measure on $(0,\infty)$ defined by
 \beqnn
 N^{(n)}_\zeta(dy): = \sum_{i=1}^{[n^\alpha \zeta ]} \mathbf{1}_{\{\ell_{0,i}/n \in dy\}} .
 \eeqnn
 From this and the change of variables, we can write $X_{\zeta,0}^{(n)}$ into
 \beqnn
 X_{\zeta,0}^{(n)}(t)
 =  \int_0^\infty \Big(\frac{ \mathbf{1}_{\{y>t\}}}{n^{\alpha}} + \int_0^{t}  n^{1-\alpha} R_\Lambda^{(n)}(n(t-s)) \cdot  \mathbf{1}_{\{y>s\}} ds \Big) N^{(n)}_\zeta(dy), \quad t\geq 0.
 \eeqnn
 Obviously, the convergence of $X_{\zeta,0}^{(n)}$ relies on the asymptotic behavior of $\{n^{1-\alpha}R_\Lambda^{(n)}(nt):t\geq 0\}$ and $N^{(n)}_\zeta$.
 \begin{lemma}\label{ResolventConvergence}
 	We have $\int_0^\cdot n^{1-\alpha}R_\Lambda^{(n)}(ns) ds \overset{\rm u.c.}\to W $ as $n\to\infty $.
 \end{lemma}
 Notice that
 $ n^\alpha\mathbf{P}(\ell_{0,i}/n \in dy ) \to  \alpha y^{-\alpha-1}dy=\bar\nu_\alpha(y)dy$ vaguely in $M(\mathbb{R}_+)$.
 The next proposition follows directly from the basic convergence theorem of empirical measures; see Theorem~5.3 in \cite[p.138]{Resnick2007}.
 \begin{proposition}
 	We have $N^{(n)}_\zeta(dy)\to N_\zeta(dy):= N_0(dy,[0,\zeta])$ weakly in $M_p(\mathbb{R}_+)$ as $n\to\infty $.
 \end{proposition}

 We start to consider the convergence of the SVI in (\ref{Eqn.HR}).
 By (\ref{Resolvent.n02}),
 \beqnn
 n^{-\alpha}R^{(n)}(nt,ny) = n^{-\alpha} \mathbf{1}_{\{y>t\}} + \int_0^{t} n^{1-\alpha}R^{(n)}_\Lambda(n(t-s)) \mathbf{1}_{\{y>s\}} ds,\quad t,y\geq 0.
 \eeqnn
 Clearly, the first term on the right side of this equality vanishes uniformly in $t,y\in\mathbb{R}$ as $n\to\infty$.
 Additionally, using the change of variable to the second term,
 \beqnn
 \int_0^{t} n^{1-\alpha}R^{(n)}_\Lambda\big(n(t-s)\big) \mathbf{1}_{\{y>s\}} ds = \int_{(t-y)^+}^{t}  n^{1-\alpha}R^{(n)}_\Lambda(ns)  ds .
 \eeqnn
 By Lemma~\ref{ResolventConvergence} and the fact that $W(x)=0$ if $x<0$, it can be well approximated by $W(t)-W(t-y)=\nabla_yW(t)$.
 Thus the SVI in (\ref{Eqn.HR}) can be well approximated by
 \beqnn
 M^{(n)}(t) :=  \int_0^t \int_0^\infty \int_0^{ X_{\zeta}^{(n)}(s-)} \nabla_y W(t-s)  \widetilde{N}^{(n)}(n\cdot ds,n\cdot dy,n^{\alpha}\cdot dz)
 \eeqnn
 with the error denoted as $ \varepsilon^{(n)}(t)$.
 By the preceding analysis and notation,  we can write the SVE (\ref{Eqn.HR})  into
 \beqlb\label{ApproxHR}
 X_{\zeta}^{(n)}(t) = \varepsilon^{(n)}(t)+ X_{\zeta,0}^{(n)}(t)+ M^{(n)}(t),\quad t\geq 0.
 \eeqlb
 As the last preparation,  the next two lemmas establish some convergence results for $\{ M^{(n)} \}_{n\geq1} $ and $\{ \varepsilon^{(n)} \}_{n\geq 1}$.
 
 \begin{lemma}\label{Lemma.TightM}
 	The sequence $\{ M^{(n)} \}_{n\geq1} $ is $C$-tight. 
 \end{lemma}
 
 \begin{lemma}\label{Lemma.Error02}
 	We have $ \varepsilon^{(n)} + X_{\zeta,0}^{(n)} \to X_{\zeta,0} $ weakly in $D([0,\infty),\mathbb{R})$ as $n\to\infty$, where
 	\beqlb\label{MartingaleX}
 	X_{\zeta,0}(t):= \int_0^\infty \nabla_yW(t)N_\zeta(dy) =\int_0^\infty \int_0^\zeta \nabla_yW(t)N_0(dy,dz)  ,\quad t\geq 0.
 	\eeqlb
 \end{lemma}

 \textsc{ Proof for Theorem~\ref{MainThm.01}(1).}
 Corollary~3.33(b) in \cite[p.353]{JacodShiryaev2003}, together with  Lemma~\ref{ConvergenceLocalTime}, \ref{Lemma.TightM} and \ref{Lemma.Error02}, yields the $C$-tightness of the sequence $\{ (X_{\zeta}^{(n)}, \varepsilon^{(n)}+X_{\zeta,0}^{(n)}, M^{(n)})  \}_{n\geq1}  $.
 Additionally, if
 \beqlb\label{ConvergeFFD}
 \big(X_{\zeta}^{(n)}, \varepsilon^{(n)}+X_{\zeta,0}^{(n)}, M^{(n)}\big) \to \big( X_\zeta,  X_{\zeta,0}, M  \big),
 \eeqlb
 as $n\to\infty$ in the sense of finite-dimensional distributions, where
 \beqlb\label{MartingaleM}
 M(t):= \int_0^t \int_0^\infty \int_0^{X_{\zeta}(s)} \nabla_y W(t-s) \widetilde{N}_\alpha(ds,dy,dz),\quad t\geq 0,
 \eeqlb
 then Theorem~13.1 in \cite[p.139]{Billingsley1999} implies  that $ (X_{\zeta}^{(n)}, \varepsilon^{(n)}+ X_{\zeta,0}^{(n)}, M^{(n)}) $ converges weakly to $ ( X_\zeta,  X_{\zeta,0}, M  ) $ in $D([0,\infty);\mathbb{R}^3)$.
 By Proposition~2.4 in \cite[p.339]{JacodShiryaev2003} and the continuous mapping theorem,
 \beqnn
 \sup_{t\in[0,T]} \big| X_\zeta(t)- X_{\zeta,0}(t)- M(t) \big|   \ar \overset{\rm d}=\ar \lim_{n\to\infty} \sup_{t\in[0,T]} \big|  X_{\zeta}^{(n)}(t)-  \varepsilon^{(n)}(t) -X_{\zeta,0}^{(n)}(t) -M^{(n)}(t) \big| ,
 \eeqnn
 which equal almost surely to $0$ for any $T\geq 0$.
 In conclusion, we have $X_\zeta\overset{\rm a.s.}=X_{\zeta,0}+ M$ and hence $X_\zeta$ is a weak solution of the SVE (\ref{MainThm.SVE}).
 
 To finish this proof, it remains to prove the finite-dimensional convergence (\ref{ConvergeFFD}).
 It follows if for any $d\geq 1$ and $0\leq t_1<\cdots<t_d$, the $3d$-dimensional random variable
 \beqnn
 Z_d^{(n)}:=\big(X_{\zeta}^{(n)}(t_i),  \varepsilon^{(n)} (t_i) + X_{\zeta,0}^{(n)}(t_1),M^{(n)}(t_i)\big)_{i=1,\cdots,d}
 \eeqnn
 converges in distribution to
 \beqnn
 Z_d:=\big( X_\zeta(t_i), X_{\zeta,0}(t_i),M(t_i)  \big)_{i=1,\cdots,d}. 
 \eeqnn 
 For $i\in\{1,\dots, d\} $ and $t\geq 0$, let
 \beqlb
 M_i^{(n)}(t)\ar:=\ar \int_0^t \int_0^\infty \int_0^{ X_{\zeta}^{(n)}(s-)} \nabla_y W(t_i-s) \widetilde{N}^{(n)}(n\cdot ds,n\cdot dy,n^{\alpha}\cdot dz), \label{eqn.Mi01}\\
 M_i (t)\ar:=\ar \int_0^t \int_0^\infty \int_0^{ X_{\zeta} (s)} \nabla_y W(t_i-s)  \widetilde{N}_\alpha( ds, dy, dz) . \label{eqn.Mi02}
 \eeqlb
 Similarly as in Remark~\ref{Remark.ItoI}, we can identify that both $ M_i^{(n)}$ and $ M_i$ are well-defined and local martingales.
 Notice that $M_i^{(n)}(t_i)\overset{\rm a.s.}=M^{(n)}(t_i)$ and $M_i(t_i)\overset{\rm a.s.}=M(t_i)$ for $i=1,\cdots,d$. Then
 \beqnn
 Z_d^{(n)}\ar \overset{\rm d}= \ar  \big(X_{\zeta}^{(n)}(t_i),  \varepsilon^{(n)} (t_i) + X_{\zeta,0}^{(n)}(t_i),M^{(n)}_1(t_i)\big)_{i=1,\cdots,d},\cr
 \ar\ar\cr
 Z_d \ar \overset{\rm d}= \ar \big( X_\zeta(t_i),X_{\zeta,0}(t_i),M_1(t_i)  \big)_{i=1,\cdots,d} .
 \eeqnn
 By the continuity of $W$ and the fact that $N_\alpha(\{t\},\mathbb{R}_+^2) \overset{\rm a.s.}=0$ for any $t>0$, the set 
 \beqnn
 \{ t\geq 0: \mathbf{P}(| Z_d(t)-Z_d(t-)|\neq 0)>0\}
 \eeqnn
 is null.
 From Proposition~3.14 in \cite[p.349]{JacodShiryaev2003}, we have  $Z_d^{(n)} \overset{\rm d}\to Z_d$ if
 \beqlb\label{Convergence.01}
  \big(X_{\zeta}^{(n)}, \varepsilon^{(n)}+X_{\zeta,0}^{(n)}, M_1^{(n)},\cdots, M_d^{(n)}\big) \to \big(X_{\zeta},X_{\zeta,0}, M_1 ,\cdots, M_d \big),
 \eeqlb
 weakly in $D([0,\infty),\mathbb{R}^{2+d})$ as $n\to\infty$.
 
 We start to prove (\ref{Convergence.01}) by using the weak convergence results established in \cite{KurtzProtter1996} for It\^o's stochastic integrals driven by infinite-dimensional semimartingale; see Appendix~\ref{AppendixB}.
 Let $\hat\nu_\alpha(dy,dz):= \nu_{\alpha}(dy)dz$ be a $\sigma$-finite measure on $\mathbb{R}_+^2$ and $L^2(\hat\nu_\alpha)$ the collection of all functions on $\mathbb{R}_+^2$ that are square-integrable with respect to $\hat\nu_\alpha$.
 For $n\geq 1$ and $t>0$, let 
 \beqnn
 \widetilde{\bf N}^{(n)}(t):=\widetilde{N}^{(n)}((0,nt],n\cdot dy,n^{\alpha}\cdot dz)
 \quad \mbox{and}\quad
 \widetilde{\bf N}_\alpha(t):=\widetilde{N}_\alpha((0,t],dy, dz),
 \eeqnn 
 which are two standard $L^2(\hat\nu_\alpha)^\#$-martingales.
 We can write (\ref{eqn.Mi01}) and (\ref{eqn.Mi02}) into 
 \beqnn
 (M_1^{(n)}(t),\cdots,M_d^{(n)}(t))\ar=\ar F( X_\zeta^{(n)},-) \cdot  \widetilde{\bf N}^{(n)}(t),\cr
 (M_1(t),\cdots,M_d(t))\ar=\ar F(X_{ \zeta},-) \cdot  \widetilde{\bf N}_\alpha (t),
 \eeqnn
 where the function $F:D([0,\infty);\mathbb{R}_+)\times\mathbb{R}_+\mapsto (L^2(\hat\nu))^d$ is defined by
 \beqnn
 F(x,s):= \big(\nabla_y W(t_1-s) , \cdots ,\nabla_y W(t_d-s) \big)\cdot \mathbf{1}_{\{0<z\leq x(s-)\}}.
 \eeqnn
 From Condition~\ref{Main.Condition} and Theorem~2.7 in \cite{KurtzProtter1996}, we have $\gamma_n n^{\alpha+1}\Lambda(n\cdot dy)\to \nu_\alpha(dy)$ vaguely and hence $\widetilde{\bf N}^{(n)}  \Rightarrow \widetilde{\bf N}_\alpha$ as $n\to\infty$.
 Similarly as in the proof of Lemma~4.9 in \cite{HorstXu2022}, we can prove that the sequence of $L^2(\hat\nu)^\#$-martingales $\{\widetilde{\bf N}^{(n)} \}_{n\geq1} $ is uniformly tight.
 By Example~5.3 in \cite{KurtzProtter1991}, the function $F$ satisfies Condition~C.2 and also Condition~C.1 in \cite[p.248-249]{KurtzProtter1996}, which induces that
 \beqnn
 \big(X_\zeta^{(n)},  \varepsilon^{(n)} +X_{\zeta,0}^{(n)}, F(X_\zeta^{(n)},-), \widetilde{\bf N}^{(n)} \big) \Rightarrow \big(X_\zeta,X_{\zeta,0}, F(X_\zeta,-), \widetilde{\bf N}_\alpha\big).
 \eeqnn
 By Theorem~4.2 or 5.5 in \cite{KurtzProtter1996}, there exists a filtration such that  $(X_\zeta,\tilde{\bf N}_\alpha)$ is adapted and
 \beqnn
 \ar\ar \big(X_\zeta^{(n)},  \varepsilon^{(n)} +X_{\zeta,0}^{(n)}, F(X_\zeta^{(n)},-), \widetilde{\bf N}^{(n)}, F(X_\zeta^{(n)},-)\cdot \widetilde{ \bf N}^{(n)} \big)  \cr 
 \ar\ar\cr
 \ar\ar \Rightarrow   \big(X_\zeta,X_{\zeta,0}, F(X_\zeta,-), \widetilde{\bf N}_\alpha, F(X_\zeta,-)\cdot \widetilde{\bf N}_\alpha \big),
 \eeqnn
 as $n\to\infty$, which implies (\ref{Convergence.01}) immediately.
 The whole proof has been finished.
 \qed
 
 {\sc Proof for Theorem~\ref{MainThm.01}(2).}
 From (\ref{eqn.500}), the inverse local time $\tau_\zeta^L$ is a subordinator killed at an independent exponential time $\varrho$. For any $\zeta>0$, we have $\tau_\xi^L(\zeta)<\infty$ if and only if $\varrho>\zeta$. Moreover,
 the independent increments of $\xi$ and the memorylessness of $\varrho$ yield that conditioned on $\tau_\xi^L(\zeta)<\infty$, the two random elements $L_\xi(\cdot ,\tau_\xi^L(\zeta) )$ and $\varrho$ are independent.
 Hence
 \beqlb\label{eqn.300}
 \lefteqn{\mathbf{P}\big(L_\xi(\cdot ,\tau_\xi^L(\zeta) )\in A,  \varrho \in d\zeta | \tau_\xi^L(\zeta)<\infty \big)}\ar\ar\cr
 \ar\ar\cr
 \ar=\ar  \mathbf{P}\big(L_\xi(\cdot ,\tau_\xi^L(\zeta) )\in A | \tau_\xi^L(\zeta)<\infty )\mathbf{P}\big(\varrho \in d\zeta| \tau_\xi^L(\zeta)<\infty \big)\cr
 \ar\ar\cr
 \ar=\ar   \mathbf{P}\big(L^\xi_\zeta \in A \big) \mathbf{P}\big(\varrho \in d\zeta | \varrho>\zeta \big).
 \eeqlb
 Notice that the event $\{L_\xi(0,\infty) \in d\zeta\}$ occurs if and only if $\{ \tau_\xi^L(\zeta)<\infty, \varrho \in d\zeta  \}$.
 By (\ref{eqn.Linf}) and (\ref{eqn.300}),  we have for any Borel set $A$ of $ C([0,\infty),\mathbb{R}_+) $,
 \beqnn
 \mathbf{P}\big(L_\infty^\xi\in A\big)
 \ar=\ar \mathbf{P}\big(L_\xi(\cdot ,\tau_\xi^L(L_\xi(0,\infty)))\in A \big)\cr
 \ar\ar\cr
 \ar=\ar \int_0^\infty \mathbf{P}\big(L_\xi(\cdot ,\tau_\xi^L(\zeta) )\in A, L_\xi(0,\infty) \in d\zeta \big)\cr
 \ar=\ar \int_0^\infty \mathbf{P}\big(L_\xi(\cdot ,\tau_\xi^L(\zeta) )\in A, \tau_\xi^L(\zeta)<\infty ,\varrho \in d\zeta \big)\cr
 \ar=\ar \int_0^\infty  \mathbf{P}\big( \tau_\xi^L(\zeta)<\infty  \big) \mathbf{P}\big(L_\xi(\cdot ,\tau_\xi^L(\zeta) )\in A,  \varrho \in d\zeta | \tau_\xi^L(\zeta)<\infty \big)\cr
 \ar=\ar  \int_0^\infty  \mathbf{P}\big(  \varrho>\zeta  \big) \mathbf{P}\big(L^\xi_\zeta \in A \big) \mathbf{P}\big(\varrho \in d\zeta | \varrho>\zeta \big) \cr
 \ar=\ar \int_0^\infty   \mathbf{P}\big(L^\xi_\zeta \in A \big) \mathbf{P}\big(\varrho \in d\zeta \big).
 \eeqnn
 Hence the process $L_\infty^\xi$ is a weak solution of (\ref{MainThm.SVE}) with $\zeta =\varrho$.
 \qed
 
 
 \subsection{Proofs for auxiliary lemmas}\label{ProofAuxilliaryLemma}
 
 In this section we provide the detailed proofs for the auxiliary lemmas in the last section.
 Denote by $\mathcal{L}_{\bar\Lambda}$, $\mathcal{L}_{\Lambda^*}$ and $\mathcal{L}_{R_\Lambda^{(n)}}$  the Laplace transforms of $\bar\Lambda$, $\Lambda^*$ and $R_\Lambda^{(n)}$ respectively, i.e., for $\lambda \geq 0$, 
 \beqnn
 \mathcal{L}_{\bar\Lambda}(\lambda)\ar:=\ar \int_0^\infty e^{-\lambda x} \bar\Lambda(x)dx, \cr
 \quad  \mathcal{L}_{\Lambda^*}(\lambda)\ar:=\ar \int_0^\infty e^{-\lambda x}  \Lambda^*(dx),\cr
 \quad \mathcal{L}_{R_\Lambda^{(n)}}(\lambda)\ar:=\ar \int_0^\infty e^{-\lambda x}   R^{(n)}_\Lambda(x)dx.
 \eeqnn

 \subsubsection{Upper bounds for resolvents} As a preparation, we first give some upper bound estimates for the resolvent $R^{(n)}_\Lambda$ and $R^{(n)}$, which will play an important role in the following proofs and analyses.
 Notice that $\bar\Lambda$ is a complete monotone function.
 Theorem 5.3.1 in \cite[p.148]{GripenbergLondenStaffans1990} tells that $R^{(n)}_\Lambda $ is also completely monotone and $\|R^{(n)}_\Lambda\|_{L^\infty}<\infty$ .
 The next proposition gives a uniform upper bound estimate for $\{R^{(n)}_\Lambda\}_{n\geq1}$.

 \begin{proposition}\label{UpperBoundResolvent}
 	There exists a constant $C>0$ such that for any  $t\geq 0$ and $n\geq 1$,
 	\beqlb\label{eqn.UpperBoundResolvent}
 	R^{(n)}_\Lambda (t) \leq   C \cdot \big(1+t\big)^{\alpha-1} .
 	\eeqlb	
 \end{proposition}
 \proof
 Let $R_\Lambda$ be the resolvent associated to $\alpha\bar\Lambda$, i.e., 
 $$ R_\Lambda(t)=\alpha\bar\Lambda(t)+\alpha\bar\Lambda* R_\Lambda(t),
 \quad t\geq 0 . $$
 It is bounded and completely monotone with $R_\Lambda(0)=\alpha$.
 It is easy to identify that  $R^{(n)}_\Lambda $ and $ R_\Lambda$ have the representations
 \beqnn
 R^{(n)}_\Lambda  = \sum_{k=1}^\infty (\gamma_n \cdot \bar\Lambda)^{*k}
 \quad\mbox{and}\quad
 R_\Lambda = \sum_{k=1}^\infty (\alpha\cdot \bar\Lambda)^{*k}.
 \eeqnn
 Here $f^{*k}$ is the $k$-th convolution of $f$.
 Since $\gamma_n\leq \alpha$, we have $R^{(n)}_\Lambda(t)\leq   R_\Lambda(t)$ uniformly in $n\geq 1$ and $t\geq 0$.
 Let $ \bar{R}_\Lambda:= 1 - R_\Lambda/\alpha$ for $t\geq 0$, which is a distribution function on $\mathbb{R}_+$.
 Using the integration by parts,
 \beqlb\label{eqn.4.11}
 \int_0^\infty e^{-\lambda t} d \bar{R}_\Lambda(t)
 \ar=\ar \lambda\int_0^\infty e^{-\lambda t}\bar{R}_\Lambda(t)dt \cr
 \ar=\ar\lambda\int_0^\infty e^{-\lambda t} \big( 1 - R_\Lambda(t)/\alpha\big)dt\cr
 \ar=\ar 1 - \frac{\lambda}{\alpha} \mathcal{L}_{R_\Lambda}(\lambda),
 \eeqlb
 where $ \mathcal{L}_{R_\Lambda}$ denotes the Laplace transform of $R_\Lambda$.
 Taking Laplace transforms of both sides of $R_\Lambda=\alpha \bar{\Lambda} + \alpha \bar{\Lambda} *R_\Lambda$, we have $\mathcal{L}_{R_\Lambda}(\lambda)= \alpha\mathcal{L}_{\bar\Lambda}(\lambda)\big[1+ \mathcal{L}_{R_\Lambda}(\lambda)\big]$ and
 \beqnn
 \mathcal{L}_{R_\Lambda}(\lambda) \ar=\ar \frac{ \alpha \mathcal{L}_{\bar\Lambda}(\lambda)}{1-\alpha \mathcal{L}_{\bar\Lambda}(\lambda)},\quad \lambda>0.
 \eeqnn
 It is obvious that the numerator goes to $1$ as $\lambda\to 0+$.
 Moreover, a simple calculation with (\ref{ParetoDis}) induces that $$\int_t^\infty \alpha \bar\Lambda(s)ds= (1+t)^{-\alpha},\quad t\geq 0.$$  
 By Karamata's Tauberian theorem\footnote{Let $F$ be a probability distribution function on $\mathbb{R}_+$ with Laplace transform $\mathcal{L}_F$. If $(1-F(tx))/(1-F(t))\to x^{-\alpha}$ or $(1-\mathcal{L}_F(1/(t\lambda)))/(1-\mathcal{L}_F(1/t))\to x^{-\alpha}$ as $t\to\infty$ for some $\alpha\in(0,1)$, then $1-\mathcal{L}_F(1/t) \sim \Gamma(1-\alpha) (1-F(t))$ as $t\to \infty$. \label{FootnoteKaraTaub}}; see Theorem~8.1.6 in \cite[p.333]{BinghamGoldieTeugels1987}, we have as $\lambda \to 0+$,
 \beqnn
 1-\alpha \mathcal{L}_{\bar\Lambda}(\lambda) \sim  \Gamma(1-\alpha) \lambda^{\alpha}
 \quad \mbox{and}\quad 
 \mathcal{L}_{R_\Lambda}(\lambda)  \sim \lambda^{-\alpha}/\Gamma(1-\alpha).
 \eeqnn 
 Taking these back into (\ref{eqn.4.11}),
 \beqnn
 1-\int_0^\infty e^{-\lambda t} d \bar{R}_\Lambda(t)
 \sim  \frac{\lambda^{1-\alpha}}{\alpha\Gamma(1-\alpha) }.
 \eeqnn
 Using Karamata's Tauberian theorem again, we have as $t \to \infty$,
 \beqnn
 1- \bar{R}_\Lambda(t) \sim  \frac{ t^{\alpha-1}}{\alpha \Gamma(\alpha)\Gamma(1-\alpha) }
 \quad \mbox{and}\quad
 R_\Lambda(t) =\alpha(1- \bar{R}_\Lambda(t))  \sim  \frac{ t^{\alpha-1}}{\Gamma(\alpha)\Gamma(1-\alpha) }.
 \eeqnn
 Consequently, there exists a constant $C>0$ such that $R^{(n)}_\Lambda(t) \leq R_\Lambda(t)  \leq C(1+t)^{\alpha-1}$ for any $t\geq 0$.
 \qed

 \begin{proposition}\label{UpperBoundResolvent02}
 	For any $p> 1+\alpha$, there exists a constant $C>0$ such that for any $t\geq 0$ and $n\geq 1$,
 	\beqlb\label{eqn.4.16}
 	n^{\alpha+1}\int_0^t  ds \int_0^\infty \big|n^{-\alpha }R^{(n)}(ns,ny)\big|^p \Lambda (n\cdot dy)\leq C\cdot (1+t)^{(p-1)\alpha}.
 	\eeqlb
 \end{proposition}
 \proof  
 By (\ref{Resolvent.n02}) and  the power mean inequality,
 the term on the left side of (\ref{eqn.4.16}) can be bounded by the sum of the following four terms
 \beqnn
 I_1^{(n)}(t)\ar:=\ar C\int_0^t ds\int_0^\infty \big|n^{-\alpha } \mathbf{1}_{\{y>s\}}\big|^p  n^{\alpha+1}\Lambda (n\cdot dy),\cr
 I_2^{(n)}(t)\ar:=\ar C\int_0^t\Big| \int_{0}^{s}n^{1-\alpha} R^{(n)}_\Lambda(nr)dr  \Big|^p n^{\alpha+1}\bar\Lambda(ns) ds, \cr
 I_3^{(n)}(t)\ar:=\ar C\int_0^t ds\int_0^{s/2} \Big| \int_{s-y}^{s}n^{1-\alpha} R^{(n)}_\Lambda(nr)dr  \Big|^p n^{\alpha+1}\Lambda(n\cdot dy),\cr
 I_4^{(n)}(t)\ar:=\ar C\int_0^t ds\int_{s/2}^s \Big| \int_{s-y}^{s}n^{1-\alpha} R^{(n)}_\Lambda(nr)dr  \Big|^p n^{\alpha+1}\Lambda(n\cdot dy),
 \eeqnn
 for some constant $C$ depending only on $p$.
 Notice that the inner integral in $I^{(n)}_1(t)$ equals to $n^{(1-p)\alpha +1}\bar\Lambda (ns)$.
 By the change of variables and the fact that $\|\bar\Lambda\|_{L^1}<\infty$,
 \beqnn
 \sup_{t\geq 0}|I_1^{(n)}(t)| =C\cdot n^{(1-p)\alpha  }\cdot \sup_{t\geq 0}\int_0^{nt}  \bar\Lambda (s)ds = C\cdot \|\bar\Lambda\|_{L^1}\cdot n^{(1-p)\alpha }  ,
 \eeqnn
 which vanishes as $n\to\infty$, since $p>1$.
 Integrating (\ref{eqn.UpperBoundResolvent}) induces that
 \beqnn
  \int_0^{s} R^{(n)}_\Lambda(r)dr \leq C\cdot (1+s)^\alpha,
  \eeqnn
 uniformly in $n\geq 1$ and $s\geq 0$.
 Using the change of variables to $I^{(n)}_2$,
 \beqnn
 I_2^{(n)}(t) \ar=\ar \frac{C}{n^{(p-1)\alpha}}\int_0^{nt}\Big| \int_{0}^{s}  R^{(n)}_\Lambda(r)dr  \Big|^p \bar\Lambda(s) ds \cr
 \ar\leq\ar   \frac{C}{n^{(p-1)\alpha}}\int_0^{nt} (1+s)^{(p-1)\alpha -1} ds
 \leq  C \cdot (1+t)^{(p-1)\alpha } .
 \eeqnn
 Here the constant $C>0$ is independent of $n$ and $t$.
 We now turn to consider $I^{(n)}_3$.
 Using (\ref{eqn.UpperBoundResolvent}) again, we have  $n^{1-\alpha} R^{(n)}_\Lambda(nr)\leq C\cdot r^{\alpha-1}$ uniformly in $n\geq 1$ and $r\geq0$.
 Hence uniformly in $s,y\geq 0$
 \beqlb\label{eqn.310}
 \int_{(s-y)^+}^{s}n^{1-\alpha} R^{(n)}_\Lambda(nr)dr \leq C\cdot \big( s^{\alpha}\wedge (|(s-y)^+|^{\alpha-1}\cdot y) \big).
 \eeqlb
 Plugging this back into $I^{(n)}_3$ and then using the fact that  $n^{\alpha+1}\Lambda(n\cdot dy) \leq  \nu_\alpha(dy)$, we have uniformly in $n\geq 1$ and $t\geq 0$,
 \beqnn
 I^{(n)}_3(t)
 \ar\leq\ar C \int_0^tds \int_0^{s/2}  |s-y|^{p(\alpha-1)} \cdot y^{p -\alpha-2}dy \cr 
 \ar\leq\ar   C \int_0^t s^{p(\alpha-1)}  ds \int_0^{s/2}  y^{p-\alpha-2}dy
 \leq     C\cdot t^{(p-1)\alpha}.
 \eeqnn
 Similarly, we also have $\int_{s-y}^{s}n^{1-\alpha} R^{(n)}_\Lambda(nr)dr  \leq C\int_0^{s}r^{\alpha-1}dr \leq C s^\alpha $ uniformly in $s>y\geq 0$. Then
 \beqnn
 I^{(n)}_4(t)\ar\leq\ar C\int_0^t s^{p\alpha} ds \int_{s/2}^s y^{-\alpha-2}dy
 \leq    C \int_0^t s^{(p-1)\alpha-1}ds
 \leq C\cdot t^{(p-1)\alpha},
 \eeqnn
 uniformly in $n\geq 1$ and $t\geq 0$.
 The desired result follows by putting all estimates above together.
 \qed
 
 \subsubsection{Moment estimates}
 In this section, we provide some uniform upper bounds for all moments of the sequence $\{ X_\zeta^{(n)} \}_{n\geq 1}$ and the impact of each ancestor on the population system.

 \begin{proposition}\label{MomentEstimate01}
 	For each $p\geq 0$, there exists a constant $C>0$ such that for any $t\geq 0$ and $n\geq 1$,
 	\beqnn
 	\mathbf{E}\Big[\Big| \int_0^{t} R_\Lambda^{(n)}(t-s) \cdot   \mathbf{1}_{\{\ell_{0,1}>s\}} ds \Big|^p\Big]\leq C \cdot t^{\alpha\cdot (p-1)^+}.
 	\eeqnn
 \end{proposition}
 \proof  By Jensen's inequality, (\ref{SizeBiasedLambda}) and  (\ref{eqn.UpperBoundResolvent}), there exists a constant $C>0$ such that for any $n\geq 1$,  $p\in[0,1]$ and $t\geq 0$,
 \beqnn
 \mathbf{E}\Big[\Big| \int_0^{t} R_\Lambda^{(n)}(t-s) \cdot   \mathbf{1}_{\{\ell_{0,1}>s\}} ds \Big|^p\Big]
 \ar\leq \ar \Big( \int_0^{t} R_\Lambda^{(n)}(t-s)     \mathbf{E} [\mathbf{1}_{\{\ell_{0,1}>s\}} ] ds\Big)^p \cr
 \ar\leq\ar  C\Big( \int_0^{t}  (t-s)^{\alpha-1} \frac{ds}{s^{\alpha}}\Big)^p \leq C.
 \eeqnn
 For $p>1$, by H\"older's inequality\footnote{\label{Footnote.Holder} For two integrable functions $g,h$ on $[0,T]$ and $p>1$, we have
 \beqnn
 \Big|\int_0^T |g(s) h(s)| ds\Big|^p= \Big|\int_0^T |g(s)||h(s)|^{1/p}\cdot |h(s)|^{1-1/p} ds\Big|^p\leq \int_0^T|g(s)|^p|h(s)| ds\cdot \Big|\int_0^T |h(s)| ds\Big|^{p-1}.
 \eeqnn
 } 
 and the previous upper bound,
 \beqnn
 \mathbf{E}\Big[\Big| \int_0^{t} R_\Lambda^{(n)}(t-s) \cdot   \mathbf{1}_{\{\ell_{0,1}>s\}} ds \Big|^p\Big] 
 \ar\leq\ar \|R_\Lambda^{(n)}\|_{L^1_t}^{p-1}\cdot\mathbf{E}\Big[  \int_0^{t} R_\Lambda^{(n)}(t-s) \cdot   \mathbf{1}_{\{\ell_{0,1}>s\}} ds  \Big] ,
 \eeqnn
 which can be bounded by $C\cdot  \|R_\Lambda^{(n)}\|_{L^1_t}^{p-1} $. 
 By (\ref{eqn.UpperBoundResolvent}), we have $\sup_{n\geq 1}\|R_\Lambda^{(n)}\|_{L^1_t}   \leq C  \cdot t^\alpha$ uniformly in $t\geq 0$ and the desired result follows.
 \qed

 \begin{proposition}\label{MomentEstimate}
 For each $p\geq 0$ and $T\geq0$, we have 
 \beqnn
 \sup_{n\geq 1}\sup_{t\in[0,T]}\mathbf{E} \big[\big|X_\zeta^{(n)}(t)\big|^p \big] <\infty.
 \eeqnn
 \end{proposition}
 \proof It suffices to prove the case of $p\in \mathbb{Z}_+$.
 Taking expectations on both sides of (\ref{Eqn.HR}),
 \beqnn
 \mathbf{E}[X_\zeta^{(n)}(t)]\ar\leq \ar \zeta +
 n^{-\alpha} \sum_{k=1}^{[\zeta n^\alpha]} \mathbf{E}\Big[ \int_0^{nt} R_\Lambda^{(n)}(nt-s) \mathbf{1}_{\{\ell_{0,k}>s\}} ds \Big].
 \eeqnn
 From Proposition~\ref{MomentEstimate01}, we have $ \mathbf{E}[X_\zeta^{(n)}(t)]\leq C$ uniformly in $t\geq 0$ and $n\geq 1$.
 By mathematical induction, it suffices to prove that for any $p\geq 1$, the desired $2p$-order moment estimate holds under the assumption that $\sup_{n\geq 1}\sup_{t\in[0,T]}\mathbf{E} [|X_\zeta^{(n)}(t)|^p ] <\infty$.
 By the power mean inequality, 
 \beqlb\label{MomentEstimate.Eqn01}
 \lefteqn{\mathbf{E}[|X_\zeta^{(n)}(t)|^{2p}] \leq  C\cdot \zeta^{2p} + C \mathbf{E}\Big[\Big| \int_0^{nt} R_\Lambda^{(n)}(nt-s)  n^{-\alpha}\sum_{k=1}^{[\zeta n^\alpha]} \mathbf{1}_{\{\ell_{0,k}>s\}} ds  \Big|^{2p} \Big]}\ar\ar\cr
 \ar\ar + C\mathbf{E}\Big[\Big| \int_0^t \int_0^\infty \int_0^{X_\zeta^{(n)}(s-)} n^{-\alpha}R^{(n)}(n(t-s),ny) \widetilde{N}^{(n)} (n\cdot ds,n\cdot dy,n^{\alpha}\cdot dz) \Big|^{2p} \Big],
 \eeqlb
 for some constant $C>0$ independent of $n$ and $t$.
 By (\ref{BDG}), the last expectation can be bounded by
 \beqnn
 \ar\ar C  \Big|  \int_0^t ds\int_0^\infty \big|n^{-\alpha }R^{(n)}(ns,ny)\big|^{2}  n^{ \alpha+1} \Lambda (n\cdot dy)  \Big|^{p} \cr
 \ar\ar +C  \int_0^t ds\int_0^\infty \big|n^{-\alpha }R^{(n)}(ns,ny)\big|^{2p}  n^{ \alpha+1} \Lambda (n\cdot dy),
 \eeqnn
 uniformly in $n\geq 1$ and $t\in[0,T]$.
 A simple calculation, together with (\ref{eqn.4.16}), implies that these two terms can be bounded  uniformly in $t\geq 0$ and $n\geq 1$ by $C\cdot (1+t)^{\alpha p} $ and $C\cdot  (1+t)^{\alpha(2p-1)}$ respectively.
 We now consider the first expectation on the right side of (\ref{MomentEstimate.Eqn01}).
 Notice that
 \beqnn
 \Big|\sum_{k=1}^{[\zeta n^\alpha]} \int_0^{nt} R_\Lambda^{(n)}(nt-s) \mathbf{1}_{\{\ell_{0,k}>s\}} ds  \Big|^{2p}
 \ar=\ar \sum_{|\mathbf{k}^{(n)}|=2p}\prod_{i=1}^{[\zeta n^\alpha]}\Big|  \int_0^{nt} R_\Lambda^{(n)}(nt-s) \mathbf{1}_{\{\ell_{0,i}>s\}} ds  \Big|^{k_i}.
 \eeqnn
 Here the sum on the right side of this inequality is over all  $\mathbf{k}^{(n)}:=(k_1,\cdots,k_{[\zeta n^\alpha]})\in\mathbb{N}^{[\zeta n^\alpha]}$  with $|\mathbf{k}^{(n)}|:= \sum_{i=1}^{[\zeta n^\alpha]}k_i=2p$.
 By Proposition~\ref{MomentEstimate01},  we have  for some constant $C>0$ depending only on $p$,
 \beqnn
 \mathbf{E}\Big[\Big| \int_0^{nt} R_H^{(n)}(nt-s)  n^{-\alpha}\sum_{k=1}^{[\zeta n^\alpha]} \mathbf{1}_{\{\ell_{0,k}>s\}} ds  \Big|^{2p} \Big]
 \leq \frac{C}{n^{2p\alpha}}  \cdot \sum_{|\mathbf{k}^{(n)}|=2p}(nt)^{\alpha \sum_{i=1}^{[\zeta n^\alpha]} (k_i-1)^+} .
 \eeqnn
 Using the multinomial distribution and then the combination formula to the last sum,
 \beqnn
 \sum_{|\mathbf{k}^{(n)}|=2p}(nt)^{\alpha \sum_{i=1}^{[\zeta n^\alpha]} (k_i-1)^+}  =\sum_{j=1}^{[ n^\alpha\zeta ]\wedge (2p)} \,_{[n^\alpha\zeta ]}\mathrm{C}_j \cdot
 (nt)^{\alpha(2p-j) },
 \eeqnn
 Here $_{[n^\alpha\zeta ]}\mathrm{C}_j:=[n^\alpha\zeta ]!/(j!([n^\alpha\zeta ]-j)!) \leq \zeta^j \cdot n^{\alpha j}/j!$ for any $j=1,\cdots,[n^\alpha\zeta ]\wedge (2p)$. Thus there exists a constant $C>0$ such that for any $n\geq 1$ and $t\in[0,T]$,
 \beqnn
 \mathbf{E}\Big[\Big| \int_0^{nt} R_H^{(n)}(nt-s)  n^{-\alpha}\sum_{k=1}^{[\zeta n^\alpha]} \mathbf{1}_{\{\ell_{0,k}>s\}} ds  \Big|^{2p} \Big]
 \ar\leq\ar C \sum_{j=1}^{[n^\alpha\zeta ]\wedge (2p)} \zeta^{j}\cdot  \frac{t^{\alpha(2p-j) }}{j!} \cr
 \ar\ar\cr
 \ar\leq\ar C
 \cdot \zeta^{2p}\cdot (1+ t)^{\alpha(2p-1)}.
 \eeqnn
 The desired result follows by putting all estimates above together.
 \qed
 
 \subsubsection{Proof for Lemma~\ref{ResolventConvergence}}
 The desired convergence is obtained by showing that the Laplace transform
 of the measure with density $n^{1-\alpha}R_\Lambda^{(n)}(n\cdot)$ converges toward the Laplace transform of the measure with density $W'$.
 Taking the Laplace transforms of both sides of (\ref{Resolvent.n01}), we have  
 \beqnn
 \mathcal{L}_{R_\Lambda^{(n)}}(\lambda)=\gamma_n \mathcal{L}_{\bar\Lambda}(\lambda) \big[ 1+ \mathcal{L}_{R_\Lambda^{(n)}}(\lambda)\big]
 \quad\mbox{and hence}\quad 
 \mathcal{L}_{R_\Lambda^{(n)}}(\lambda)= \frac{\gamma_n \mathcal{L}_{\bar\Lambda}(\lambda)}{ 1-\gamma_n \mathcal{L}_{\bar\Lambda}(\lambda)}. 
 \eeqnn 
 By the change of variables,
 \beqnn
 \int_0^\infty e^{-\lambda t} n^{1-\alpha}R_\Lambda^{(n)}(nt) dt
 \ar=\ar \frac{\mathcal{L}_{R_\Lambda^{(n)}}(\lambda/n)}{n^{\alpha}}\cr
 \ar=\ar\frac{\gamma_n \mathcal{L}_{\bar\Lambda}(\lambda/n)}{n^{\alpha}(1-\gamma_n\mathcal{L}_{\bar\Lambda}(\lambda/n))} \cr
 \ar=\ar \frac{\gamma_n \mathcal{L}_{\bar\Lambda}(\lambda/n)}{n^{\alpha}(1-\frac{\gamma_n}{\alpha}) + \frac{\gamma_n}{\alpha}\cdot n^{\alpha}(1-\mathcal{L}_{\Lambda^*}(\lambda/n))}.
 \eeqnn
 A simple calculation, along with Condition~\ref{Main.Condition}, shows that $\gamma_n \mathcal{L}_{\bar\Lambda}(\lambda/n)\to 1$ as $n\to\infty$.
 By (\ref{SizeBiasedLambda}) and Karamata's Tauberian theorem; see footnote~\ref{FootnoteKaraTaub}, we have $n^{\alpha}\big(1-\mathcal{L}_{\Lambda^*}(\lambda/n)\big) \to \Gamma(1-\alpha) \lambda^{\alpha}$ as $n\to\infty$
 and
 \beqnn
 \int_0^\infty e^{-\lambda t} n^{1-\alpha}R_\Lambda^{(n)}(nt) dt
 \to  \frac{1}{b+c\lambda^\alpha}.
 \eeqnn
 By (\ref{LaplaceDW}), the function whose Laplace transform is equal to the last quantity is $W'$.
 
 \subsubsection{Proof for Lemma~\ref{Lemma.TightM}}
 As a direct consequence of the continuity of $\nabla_yW$, the process $M^{(n)}$ is also continuous.
 Let $p>2/\alpha$.
 By the Kolmogorov tightness criterion; see Theorem 13.5 in \cite{Billingsley1999}, the sequence $\{ M^{(n)} \}_{n\geq 1}$ is $C$-tight if for any $T>0$, there exists a constant $C\geq 0$ such that for any $h\in[0,1]$,
 \beqlb\label{eqn.302}
 \sup_{n\geq 1}\sup_{t\in[0,T]} \mathbf{E} \Big[ \big| \Delta_hM^{(n)}(t)\big|^{2p}\Big] \leq  C \cdot  h^2.
 \eeqlb
 We start to prove (\ref{eqn.302}) with the help of the technical results about $W$ in Appendix~\ref{TechResW}.
 We first split  $\Delta_h M^{(n)}(t)$ into the following five parts:
 \beqnn
 M_{1}^{(n)}(t,h)\ar:=\ar  \int_t^{t+h} \int_0^{t+h-s} \int_0^{X_\zeta^{(n)}(s-)} \nabla_y W(t+h-s)  \widetilde{N}^{(n)}(n\cdot ds,n\cdot dy,n^{\alpha}\cdot dz),\cr
 M_{2}^{(n)}(t,h)\ar:=\ar  \int_t^{t+h} \int_{t+h-s}^\infty \int_0^{X_\zeta^{(n)}(s-)} W(t+h-s) \widetilde{N}^{(n)}(n\cdot ds,n\cdot dy,n^{\alpha}\cdot dz),\cr
 M_{3}^{(n)}(t,h)\ar:=\ar \int_0^t \int_0^{t-s} \int_0^{X_\zeta^{(n)}(s-)} \nabla_y\Delta_h W(t-s)  \widetilde{N}^{(n)}(n\cdot ds,n\cdot dy,n^{\alpha}\cdot dz),\cr
 M_{4}^{(n)}(t,h)\ar:=\ar \int_0^t \int_{t-s}^{t+h-s} \int_0^{X_\zeta^{(n)}(s-)} \nabla_y\Delta_h W(t-s) \widetilde{N}^{(n)}(n\cdot ds,n\cdot dy,n^{\alpha}\cdot dz) ,\cr
 M_{5}^{(n)}(t,h)\ar:=\ar \int_0^t \int_{t+h-s}^\infty \int_0^{X_\zeta^{(n)}(s-)} \Delta_h W(t-s) \widetilde{N}^{(n)}(n\cdot ds,n\cdot dy,n^{\alpha}\cdot dz).
 \eeqnn
 Applying (\ref{BDG}), Proposition~\ref{MomentEstimate} and~\ref{WeakConvergenceProp01} to $\mathbf{E}[| M^{(n)}_{1}(t,h)|^{2p} ]$,
 \beqnn
 \mathbf{E}\Big[ \big| M^{(n)}_{1}(t,h) \big|^{2p}  \Big]\ar\leq\ar C \Big|  \int_t^{t+h}  ds  \int_0^{t+h-s} \frac{|\nabla_y W(t+h-s)|^{2}}{y^{\alpha+2}} dy  \Big|^{p} \cr
 \ar\ar + C \int_t^{t+h}  ds  \int_0^{t+h-s} \frac{|\nabla_y W(t+h-s)|^{2p}}{y^{\alpha+2}}dy\cr
 \ar\leq\ar C \Big|  \int_0^{h}  ds  \int_0^\infty \frac{|\nabla_y W(s) |^{2}}{y^{\alpha+2}} dy  \Big|^{p} \cr
 \ar\ar  + C \int_0^{h}  ds  \int_0^\infty \frac{|\nabla_y W(s)|^{2p}}{y^{\alpha+2} } dy,
 \eeqnn
 which can be bounded by $C h^{p\alpha } \leq C h^{2} $ uniformly in $n\geq1$, $ t\in[0,T]$ and $h\in[0,1]$.
 With the help of Proposition~\ref{WeakConvergenceProp02} and \ref{WeakConvergenceProp03}, we also can prove the similar results for other terms. 
 The inequality (\ref{eqn.302}) follows by putting them together and the power mean inequality. 
 \qed

 \subsubsection{Proof for Lemma~\ref{Lemma.Error02}}
 From (\ref{ApproxHR}), we have $\varepsilon^{(n)} + X^{(n)}_{\zeta,0}=X_\zeta^{(n)}-M^{(n)}$.
 Notice that $X^{(n)}_\zeta\to L^\xi_\zeta$ weakly in $D([0,\infty);\mathbb{R}_+)$; see Lemma~\ref{ConvergenceLocalTime} and (\ref{eqn.28}).
 By the continuity of $L^\xi_\zeta$, the sequence $\{ X^{(n)}_\zeta\}_{n\geq 1}$ is $C$-tight.
 Together with this and Lemma~\ref{Lemma.TightM}, Corollary~3.33 in \cite[p.353]{JacodShiryaev2003} yields the $C$-tightness of the sequence $\{\varepsilon^{(n)}+X^{(n)}_{\zeta,0}\}_{n\geq 1}$ immediately.
 Hence it remains to prove $ \varepsilon^{(n)}+X^{(n)}_{\zeta,0} \to X_{\zeta,0}$ in the sense of finite-dimensional distributions, which follows directly from the next two propositions.
 
 
 \begin{proposition}
 	We have $\varepsilon^{(n)}\to 0$ in the sense of finite-dimensional distributions as $n\to\infty$.
 \end{proposition}
 \proof
 It suffices to prove $|\varepsilon^{(n)}(t)|\overset{\rm p}\to0$ as $n\to\infty$ for any $t\geq 0$.
 By (\ref{Resolvent.n02}), we split $\varepsilon^{(n)}(t)$ into the following two parts:
 \beqnn
 \varepsilon_{1}^{(n)}(t) \ar:=\ar \int_0^t \int_{0}^\infty \int_0^{X_\zeta^{(n)}(s-)} n^{-\alpha}\cdot \mathbf{1}_{\{ y>t-s  \}} \widetilde{N}^{(n)}(n\cdot ds,n\cdot dy,n^{\alpha}\cdot dz),\cr
 \varepsilon_{2}^{(n)}(t)\ar:=\ar \int_0^t \int_0^\infty \int_0^{X_\zeta^{(n)}(s-)} \Big( \int_{(t-s-y)^+}^{t-s}n^{1-\alpha} R^{(n)}_\Lambda(nr)dr \cr
 \ar\ar \qquad \qquad \qquad\qquad  -  \nabla_yW(t-s) \Big) \widetilde{N}^{(n)}(n\cdot ds,n\cdot dy,n^{\alpha}\cdot dz).
 \eeqnn
 Applying (\ref{BDG}) and Proposition~\ref{MomentEstimate} to $\mathbf{E}[| \varepsilon_{1}^{(n)}(t)|^2]$, we have
 \beqnn
 \sup_{t\geq 0} \mathbf{E}\big[|\varepsilon_{1}^{(n)}(t)|^2\big]
 \ar \leq \ar  C \sup_{t\geq 0}\int_0^t   n^{1-\alpha} \bar\Lambda(n(t-s))ds
 = C n^{-\alpha}\int_0^\infty\bar\Lambda(s)ds,
 \eeqnn
 which goes to $0$  as $n\to\infty$.
 Similarly, we also can bound $\sup_{t\in[0,T]}\mathbf{E}\big[|\varepsilon_{2}^{(n)}(t)|^2\big]$ by 
 \beqnn
  C \int_0^T ds \int_0^\infty \Big| \int_{(s-y)^+}^{s}n^{1-\alpha} R^{(n)}_\Lambda(nr)dr -  \nabla_yW(s)\Big|^2 n^{\alpha+1}\Lambda(ndy) .
 \eeqnn
 For $\vartheta\in (\frac{1+\alpha}{2},\frac{1/2}{1-\alpha}\wedge 1)$, the preceding integral can be bounded by the multiplication of the next two terms
 \beqnn
 \varepsilon_{21}^{(n)}(T)\ar:=\ar  \sup_{x\in[0,T]} \Big| \int_0^{x}n^{1-\alpha} R^{(n)}_\Lambda (nr)dr -  W(x)  \Big|^{2(1-\vartheta)}, \cr
 \varepsilon_{22}^{(n)}(T)\ar:=\ar  \int_0^T ds \int_0^\infty \Big| \int_{(s-y)^+}^{s}n^{1-\alpha} R^{(n)}_\Lambda (nr)dr - \nabla_yW(s)\Big|^{2\vartheta} n^{\alpha+1}\Lambda(n\cdot dy).
 \eeqnn
 From Lemma~\ref{ResolventConvergence}, we have $\varepsilon_{21}^{(n)}(T)\to 0$ as $n\to\infty$.
 It remains to prove that $\varepsilon_{22}^{(n)}(T)$ is bounded.
 By the power mean inequality and the fact that $n^{\alpha+1}\Lambda(n\cdot dy) \leq  \nu_\alpha(dy)$,
 \beqlb\label{eqn.4.17}
  \varepsilon_{22}^{(n)}(T)\ar\leq \ar   C\int_0^T ds \int_0^\infty \frac{| \nabla_yW(s) |^{2\vartheta}} {y^{\alpha+2} } dy \cr
 \ar\ar + C\int_0^T ds \int_0^\infty \Big| \int_{(s-y)^+}^{s}n^{1-\alpha} R^{(n)}_\Lambda(nr)dr \Big|^{2\vartheta} \frac{dy}{y^{\alpha+2}},
 \eeqlb
 uniformly in $n\geq 1$.
 By  Proposition~\ref{WeakConvergenceProp01}, the first integral on the right side of (\ref{eqn.4.17}) can be bounded by $C\cdot T^{\alpha (2\vartheta-1) }$.
 Plugging (\ref{eqn.310}) into the second integral on the right side of (\ref{eqn.4.17}), it can be bounded by
 \beqnn
 C\int_0^T ds \int_s^\infty  \frac{s^{2\vartheta\alpha}}{y^{\alpha+2}}dy + C\int_0^T ds \int_0^s  (s-y)^{2\vartheta(\alpha-1)} y^{2\vartheta-\alpha-2}dy\leq C\cdot T^{\alpha(2\vartheta-1)}.
 \eeqnn
 Hence $ \varepsilon_{22}^{(n)}(T)<\infty$ and $\sup_{t\in[0,T]}\mathbf{E}[|\varepsilon_{2}^{(n)}(t)|^2]\to 0$ as $n\to\infty$.
 The desired result follows by putting these estimates together.
 \qed

 \begin{proposition}
 	We have $X^{(n)}_{\zeta,0}\to X_{\zeta,0}$ in the sense of finite-dimensional distributions as $n\to\infty$.
 \end{proposition}
 \proof Notice that $X^{(n)}_{\zeta,0}(0)=[n^\alpha \zeta]/n^\alpha$ is deterministic and converges to $ \zeta=X_{\zeta,0}(0)$ as $n\to\infty$. For any $T>0$, $d\in\mathbb{Z}_+$,  $0 <t_1<\cdots<t_d\leq T$ and $ \lambda_1,\cdots,\lambda_d\geq 0$, let
 \beqnn
 Y_d^{(n)}(\lambda,y) \ar:=\ar  \sum_{i=1}^d \lambda_i \cdot \Big(\frac{ \mathbf{1}_{\{y>t_i\}}}{n^{\alpha}} + \int_0^{t_i}  n^{1-\alpha} R_\Lambda^{(n)}(n(t_i-s)) \cdot  \mathbf{1}_{\{y>s\}} ds \Big),\cr
 Y_d(\lambda,y)\ar:=\ar  \sum_{i=1}^d \lambda_i \int_0^{t_i} W'(t_i-s)\cdot  \mathbf{1}_{\{y>s\}} ds
 = \sum_{i=1}^d \lambda_i \cdot \nabla_y W(t_i).
 \eeqnn
 It suffices to prove that 
 \beqnn
  \mathbf{E}\Big[\exp\Big\{ - \int_0^\infty Y_d^{(n)}(\lambda, y) N_\zeta^{(n)}(dy) \Big\}\Big] \to  \mathbf{E}\Big[\exp\Big\{ - \int_0^\infty Y_d (\lambda, y) N_\zeta (dy) \Big\}\Big] .
  \eeqnn 
 By the definition of $N^{(n)}_\zeta$ we have
 \beqnn
 \lefteqn{\mathbf{E}\Big[\exp\Big\{ - \int_0^\infty Y_d^{(n)}(\lambda, y) N_\zeta^{(n)}(dy) \Big\}\Big]}\ar\ar\cr
 \ar=\ar \Big(  \mathbf{E} \Big[e^{ - Y_d^{(n)}(\lambda, \ell_{0,1}) }\Big] \Big)^{[n^\alpha\zeta]} \cr
 \ar=\ar \Big( 1-  \frac{1}{n^\alpha}  \int_0^\infty \big(1-e^{ - Y_d^{(n)}(\lambda,y) }\big)\cdot n^{\alpha}\cdot  \Lambda^*(n\cdot dy) \Big)^{[n^\alpha\zeta]},
 \eeqnn
 and this converges as $n\to\infty$ to
 \beqnn
 \exp\Big\{-\zeta \int_0^\infty (1-e^{- Y_d(\lambda,y)}) \bar\nu_\alpha(y)dy  \Big\}
 =\mathbf{E} \Big[\exp \Big\{ - \int_0^\infty Y_d (\lambda, y) N_\zeta (dy) \Big\} \Big] ,
 \eeqnn
 if and only if
 \beqlb\label{eqn.200}
 \int_0^\infty \big(1- e^{ - Y_d^{(n)}(\lambda,y)}\big)\cdot n^{\alpha}\cdot  \Lambda^*(n\cdot dy) - \int_0^\infty (1-e^{- Y_d(\lambda,y)}) \bar\nu_\alpha(y)dy   \to0.
 \eeqlb
 For any $\epsilon \in (0,t_1/2)$, we can write the preceding subtraction into the sum of the following four terms
 \beqnn
 \varepsilon_{3}^{(n)}  \ar:=\ar  \int_0^\infty \big(1- e^{ - Y_d(\lambda,y) }\big) \cdot  \big(n^{\alpha}\cdot\Lambda^*(n\cdot dy) -  \bar\nu(y)dy\big) ,\cr
 \varepsilon_{4}^{(n)}  \ar:=\ar  \int_\epsilon^\infty \big( e^{ - Y_d (\lambda,y) } - e^{ - Y_d^{(n)}(\lambda,y) } \big)\cdot  n^{\alpha}\cdot \Lambda^*(n\cdot dy) ,\cr
 \varepsilon_{5}^{(n)}  \ar:=\ar \int_0^\epsilon  \big( 1 -e^{ - Y_d^{(n)}(\lambda,y) } \big)\cdot  n^{\alpha}\cdot \Lambda^*(n\cdot dy),\cr
 \varepsilon_{6}^{(n)}  \ar:=\ar \int_0^\epsilon  \big(1- e^{ - Y_d (\lambda,y) }  \big)\cdot  n^{\alpha}\cdot \Lambda^*(n\cdot dy) .
 \eeqnn
 The vague convergence of $n^{\alpha}\cdot\Lambda^*(n\cdot dy) $ to $ \bar\nu_\alpha(y)dy  $ implies that $|\varepsilon_{3}^{(n)}|\to 0$ as $n\to\infty$.
 By Lemma~\ref{ResolventConvergence}, we have $ Y_d^{(n)}(\lambda,y) \to  Y_d(\lambda,y)$ uniformly in  $y>0$, which yields that as $n\to\infty$,
 \beqnn
 \big|\varepsilon_{4}^{(n)}\big| 
 \ar\leq\ar    \sup_{y>0}\big|e^{ - Y_d (\lambda,y) } - e^{ - Y_d^{(n)}(\lambda,y) } \big| \cdot n^\alpha \overline{\Lambda^*}(n\epsilon)\cr
 \ar\leq\ar  \alpha \cdot \epsilon^{-\alpha} \cdot   \sup_{y>0}\big|e^{ - Y_d (\lambda,y) } - e^{ - Y_d^{(n)}(\lambda,y) } \big| \to 0.
 \eeqnn
 Notice that for any $y\in(0,\epsilon)$,
 \beqnn
 Y_d^{(n)}(\lambda,y) \ar=\ar \sum_{i=1}^d \lambda_i  \int_0^{y}  n^{1-\alpha} R_\Lambda^{(n)}(n(t_i-s)) ds,\cr
 Y_d(\lambda,y)  \ar=\ar \sum_{i=1}^d \lambda_i \cdot \nabla_y W(t_i).
 \eeqnn 
 By the change of variables and (\ref{eqn.310}),
 \beqnn
 \sup_{n\geq 1} \big| 1- e^{ - Y_d^{(n)} (\lambda,y) }  \big|
 \ar\leq\ar   \sup_{n\geq 1}\big| Y_d^{(n)} (\lambda,y)\big|\cr
 \ar=\ar  \sup_{n\geq 1} \sum_{i=1}^d \lambda_i  \int_{t_i-y}^{t_i}  n^{1-\alpha} R_\Lambda^{(n)}(ns) ds  \leq C\cdot  |t_1|^{\alpha-1} \cdot y .
 \eeqnn
 Plugging this back into $\varepsilon_{5}^{(n)}$ and then using the fact that $ n^{\alpha}\cdot \Lambda^*(n\cdot dy) \leq \bar\nu(y)dy$,
 \beqnn
 \sup_{n\geq 1} \big| \varepsilon_{5}^{(n)}\big| \leq  C \int_0^\epsilon y\bar\nu(y)dy \leq C \int_0^\epsilon y^{-\alpha} dy \leq C \epsilon^{1-\alpha},
 \eeqnn
 which goes to $0$ as $\epsilon \to 0 +$. Similarly, we also can prove that $ \sup_{n\geq 1} | \varepsilon_{6}^{(n)}| \to 0$ as   $\epsilon \to 0+$.
 The convergence  (\ref{eqn.200})  follows directly by putting these estimates together.
 \qed
 
 \section{H\"older continuity}\label{ProofHolder}

 In this section, we prove the H\"older continuity  of $L^\xi_\zeta$; see Theorem~\ref{Thm.Regularity}.
 For simplicity on exposition, we also denote by $X_{\zeta,0}$ and $M$  the
 two terms on the right side of (\ref{MainThm.SVE}) respectively without ambiguity, i.e.,
 the SVE  (\ref{MainThm.SVE}) can be written into
 $ L^\xi_\zeta(x)=   X_{\zeta,0}(x) + M(x)$ for any $x\geq 0$. 
 
 \begin{lemma}\label{Lemma.501}
 	For each $p\geq 0$ and $\zeta>0$, the exists a constant $C>0$ such that for any $x\geq 0$,
 	\beqlb\label{MomentEsti01}
 	\mathbf{E}\Big[\big|X_{\zeta,0}(x)\big|^{p}\Big] \leq  C \cdot (1+x)^{\alpha(p-1)^+ }.
 	\eeqlb
 \end{lemma}
 \proof By (\ref{Equ.SVE}), we have $X_{\zeta,0}(x)\geq 0 $ a.s. and
 \beqlb\label{NewXzeta}
 X_{\zeta,0}(x)= \zeta \big(1-bW(x)\big)  + \int_0^\infty  \int_0^\zeta \nabla_y W(x) \widetilde{N}_0( dy,dz),\quad x\geq 0.
 \eeqlb
 When $p\leq 1$, by Jensen's inequality we have 
 $$ \mathbf{E}\Big[ \big|X_{\zeta,0}(x)\big|^p\Big]\leq |\mathbf{E}[X_{\zeta,0}(x)]|^p \leq \zeta^p .$$
 For $p>1$,
 By the power mean inequality and Theorem~\ref{Thm.BDG}, there exits a constant $C>0$ such that for any $x>0$,
 \beqlb\label{eqn.400}
 \mathbf{E}\Big[\big|X_{\zeta,0}(x)\big|^{p}\Big]
 \ar\leq \ar C\zeta^p+ C \Big| \int_0^\infty \big|\nabla_y W(x)\big|^2 \bar\nu_\alpha(y)dy  \Big|^{p/2}\cr
 \ar\ar + C \int_0^\infty \big|\nabla_y W(x)\big|^p \bar\nu_\alpha(y)dy .
 \eeqlb
 By (\ref{UpperboundW}) and (\ref{BoundFirstMom}),
 \beqnn
 \int_0^\infty \big|\nabla_y W(x)\big|^p \bar\nu_\alpha(y)dy \leq  \int_0^\infty \nabla_y W(x)  \bar\nu_\alpha(y)dy \cdot \big|W(x)\big|^{p-1} \leq C x^{(p-1)\alpha}.
 \eeqnn
 The desired result follows by plugging this back into (\ref{eqn.400}).
 \qed
 
 \begin{lemma}\label{MomentIncreX}
 	For each $p\geq1$, there exists a constant $C>0$ such that for any $ x_1,x_2\geq 0$,
 	\beqnn
 	\mathbf{E}\Big[\big|X_{\zeta,0}(x_2) -X_{\zeta,0}(x_1)\big|^{2p}\Big] \leq C\cdot  \big(|x_2-x_1|^{2p\alpha}+|x_2-x_1|^{p\alpha}\big).
 	\eeqnn
 \end{lemma}
 \proof
 By the power mean inequality and (\ref{NewXzeta}), there exists a constant $C>0$ such that for any $ x_1,x_2\geq 0$, the expectation $ \mathbf{E} [ |X_{\zeta,0}(x_2) -X_{\zeta,0}(x_1) |^{2p} ]$ can be bounded by
 \beqnn
 C   \big|W(x_2)-W(x_1)\big|^{2p}  + C\mathbf{E}\Big[\Big| \int_0^\infty \int_0^\zeta  \big(\nabla_y W(x_2)-\nabla_y W(x_1)\big) \widetilde{N}_0(dy,dz) \Big|^{2p}  \Big].
 \eeqnn
 By the uniform $\alpha$-H\"older continuity of $W$ on $\mathbb{R}_+$; see Section~\ref{StableProc},  the first term can be bounded by $C\cdot |x_2-x_1|^{2p\alpha}$  uniformly in  $ x_1,x_2\geq 0$.
 Applying (\ref{BDG}) to the second term, it can be bounded by
 \beqlb\label{eqn.4002}
 \ar\ar C \Big| \int_0^\infty  |\nabla_y W(x_2)-\nabla_y W(x_1)|^2 \bar\nu(y)dy\Big|^p \cr
 \ar\ar +C \int_0^\infty  |\nabla_y W(x_2)-\nabla_y W(x_1)|^{2p} \bar\nu(y)dy,
 \eeqlb
 for some constant $C>0$ independent of $x_1,x_2$.
 Using the uniform $\alpha$-H\"older continuity of $W$ on $(0,\infty)$ again, we have $ |\nabla_y W(x_2)-\nabla_y W(x_1)|\leq C\big(|x_2-x_2|^\alpha\wedge y^\alpha\big)$ uniformly in $x_1,x_2,y\geq 0$.
 Plugging this into (\ref{eqn.4002}), it can be bounded by
 \beqnn 
 \ar\ar C  \Big| \int_0^\infty  \big(|x_2-x_2|^\alpha\wedge y^\alpha\big)^2 \bar\nu(y)dy\Big|^p   +C \int_0^\infty  \big(|x_2-x_2|^\alpha\wedge y^\alpha\big)^{2p} \bar\nu(y)dy ,
 \eeqnn
 which bounded by $C\cdot |x_2-x_1|^{p\alpha}$ for some constant $C>0$ independent of $x_1,x_2$.
 The desired result follows by putting these estimates together.
 \qed

 \begin{lemma}\label{MomentM}
 	For each $p\geq 0$ and $\zeta>0$, there exists a constant $C>0$ such that for any $x\geq 0$,
 	\beqlb\label{MomentEsti}
 	\mathbf{E}\big[|L^\xi_\zeta(x)|^{p}\big] \leq  C \cdot  (1 + x )^{p\alpha }.
 	\eeqlb
 \end{lemma}
 \proof Here we just prove this lemma with $p=2^k$ and $k\in\mathbb{N}$.
 The general case can be proved in the same way.
 When $k=0$, by (\ref{BoundFirstMom}) we have $\mathbf{E}[L^\xi_\zeta(x)]\leq \zeta$.
 For $k\geq 1$, by mathematical induction it suffices to prove (\ref{MomentEsti}) for $p=2^k$ under the assumption that it holds for $p=2^{i}$ with $i=0,1,\cdots, k-1$.
 Applying the Burkholder-Davis-Gundy inequality and the power mean inequality to $\mathbf{E}[|M(x)|^{2^k}]$,
 \beqlb\label{eqn.305}
 \mathbf{E}\Big[|M(x)|^{2^k}\Big]\ar\leq\ar  C\cdot \mathbf{E}\Big[ \Big|\int_0^x\int_0^\infty \int_0^{L^\xi_\zeta(s)} \big|\nabla_y W(x-s) \big|^2 N_\alpha(ds,dy,dz)\Big|^{2^{k-1}}\Big] \cr
 \ar\leq\ar C\cdot\mathbf{E}\Big[ \Big|\int_0^x\int_0^\infty \int_0^{L^\xi_\zeta(s)} \big|\nabla_y W(x-s) \big|^2 \widetilde{N}_\alpha(ds,dy,dz)\Big|^{2^{k-1}}\Big]\cr
 \ar\ar + C\cdot\mathbf{E}\Big[ \Big|\int_0^x ds\int_0^\infty  L^\xi_\zeta(s) \big|\nabla_y W(x-s) \big|^2 \nu_\alpha(dy)\Big|^{2^{k-1}}\Big]\cr
 \ar\leq\ar C\cdot \sum_{i=1}^{k}\mathbf{E}\Big[\Big| \int_0^x ds\int_0^\infty  L^\xi_\zeta(s) \big|\nabla_y W(x-s) \big|^{2^i} \nu_\alpha(dy)\Big|^{2^{k-i}}\Big],
 \eeqlb
 with $C>0$ depending only on $k$.
 Applying H\"older's inequality to each term in the last sum; see footnote~\ref{Footnote.Holder},
 \beqnn
 \lefteqn{\mathbf{E}\Big[\Big| \int_0^x ds \int_0^\infty  L^\xi_\zeta(s) \big|\nabla_y W(x-s) \big|^{2^i} \nu_\alpha(dy)\Big|^{2^{k-i}}\Big]}\ar\ar\cr
 \ar\leq\ar  \int_0^x ds \int_0^\infty \mathbf{E}\big[| L^\xi_\zeta(s)|^{2^{k-i}} \big] \big|\nabla_y W(x-s) \big|^{2^i} \nu_\alpha(dy) \cr
 \ar\ar\quad \times \Big| \int_0^x ds\int_0^\infty\big|\nabla_z W(x-s) \big|^{2^i} \nu_\alpha(dz) \Big|^{2^{k-i}-1} \cr
 \ar\leq\ar \sup_{t\in[0,x]} \mathbf{E}\big[| L^\xi_\zeta(t)|^{2^{k-i}} \big] \cdot \Big| \int_0^x ds \int_0^\infty\big|\nabla_yW(s) \big|^{2^i} \nu_\alpha(dy) \Big|^{2^{k-i}}\cr
 \ar\leq\ar C\cdot \big(1+x\big)^{2^{k-i}\alpha}\cdot \Big| \int_0^x ds \int_0^\infty\big|\nabla_yW(s) \big|^{2^i} \nu_\alpha(dy) \Big|^{2^{k-i}}.
 \eeqnn
 By Proposition~\ref{WeakConvergenceProp01} with $p=2^i$,  the right hand of the last inequality can be bounded by $C\cdot (1+x)^{p\alpha } $ uniformly in $x\geq0$.
 Taking these estimates back into (\ref{eqn.305}), we have 
 $$ \mathbf{E}\Big[\big|M(x)\big|^{p}\Big]\leq C\cdot (1+x)^{p\alpha},$$ 
 uniformly in $x\geq 0$.
 Hence the desired result follows from this, (\ref{MomentEsti01}) and the power mean inequality.
 \qed

 \begin{lemma}\label{MomentIncreM}
 	For each $p\geq1$, there exists a constant $C>0$ such that for any $x\geq 0$ and $ x_1,x_2\in[0,x]$,
 	\beqnn
 	\mathbf{E}\Big[\big|M(x_1)-M(x_2)\big|^{2p}\Big] \leq C\cdot (1+x)^{p\alpha}\cdot |x_2-x_1|^{p\alpha}.
 	\eeqnn
 \end{lemma}
 \proof Without loss of generality, we assume $0\leq x_1< x_2\leq x$.
 Similarly as in the proof of Lemma~\ref{Lemma.TightM},  by the power mean inequality we have $$\mathbf{E}\Big[\big|M(x_1)-M(x_2)\big|^{2p}\Big] \leq 5^{2p}\cdot \sum_{i=1}^5\mathbf{E}\Big[\big| M_i(x_1,x_2)\big|^{2p}\Big]$$ 
 with
 \beqnn
 M_1(x_1,x_2)\ar:=\ar \int_{x_1}^{x_2} \int_0^{x_2-s} \int_0^{L^\xi_\zeta(s)} \nabla_yW(x_2-s)\widetilde{N}_\alpha(ds,dy,dz), \cr
 M_2(x_1,x_2)\ar:=\ar \int_{x_1}^{x_2} \int_{x_2-s}^\infty \int_0^{L^\xi_\zeta(s)} W(x_2-s)\widetilde{N}_\alpha(ds,dy,dz), \cr
 M_3(x_1,x_2)\ar:=\ar \int_0^{x_1}  \int_{x_1-s}^\infty \int_0^{L^\xi_\zeta(s)} \nabla_y \Delta_{x_2-x_1}W(x_1-s)\widetilde{N}_\alpha(ds,dy,dz), \cr
 M_4(x_1,x_2)\ar:=\ar \int_0^{x_1}  \int_{x_1-s}^{x_2-s} \int_0^{L^\xi_\zeta(s)} \nabla_y \Delta_{x_2-x_1}W(x_1-s)\widetilde{N}_\alpha(ds,dy,dz), \cr
 M_5(x_1,x_2)\ar:=\ar \int_0^{x_1}  \int_{x_2-s}^\infty \int_0^{L^\xi_\zeta(s)} \Delta_{x_2-x_1}W(x_1-s)\widetilde{N}_\alpha(ds,dy,dz).
 \eeqnn
 Applying (\ref{BDG}) to $\mathbf{E}[| M_1(x_1,x_2)|^{2p} ]$ and  then using the change of variables, there exists a constant $C>0$ depending only on $p$ such that
 \beqnn
 \mathbf{E}\Big[\big| M_1(x_1,x_2)\big|^{2p} \Big] \ar\leq\ar C \sup_{t\in[0,x]}\mathbf{E}\Big[\big|L^\xi_\zeta(t)\big|^p\Big]\cdot \Big|\int_{x_1}^{x_2} \int_0^{x_2-s} \big|\nabla_yW(x_2-s)\big|^2\nu_\alpha(dy)ds\Big|^p\cr
 \ar\ar + C \sup_{t\in[0,x]}\mathbf{E}\Big[\big|L^\xi_\zeta(t)\big|\Big]\cdot  \int_{x_1}^{x_2} \int_0^{x_2-s} \big|\nabla_yW(x_2-s)\big|^{2p}\nu_\alpha(dy)ds \cr
 \ar=\ar C \sup_{t\in[0,x]}\mathbf{E}\Big[\big|L^\xi_\zeta(t)\big|^p\Big]\cdot \Big|\int_{0}^{x_2-x_1} ds \int_0^{s} \big|\nabla_yW(s)\big|^2\nu_\alpha(dy)\Big|^p\cr
 \ar\ar + C \sup_{t\in[0,x]}\mathbf{E}\Big[\big|L^\xi_\zeta(t)\big|\Big]\cdot  \int_{0}^{x_2-x_1} ds \int_0^{s} \big|\nabla_yW(s)\big|^{2p}\nu_\alpha(dy).
 \eeqnn
 By Proposition~\ref{WeakConvergenceProp01} and Lemma~\ref{MomentM},
 the foregoing quantities can be bounded by
 \beqnn
  C \sup_{t\in[0,x]}\mathbf{E}\Big[\big|L^\xi_\zeta(t)\big|^p\Big] \cdot |x_2-x_1|^{p\alpha}   +   C \sup_{t\in[0,x]}\mathbf{E}\Big[\big|L^\xi_\zeta(t)\big|\Big] \cdot |x_2-x_1|^{\alpha(2p-1)},
 \eeqnn
 which is bounded by $ C\cdot (1+x)^{p\alpha}\cdot |x_2-x_1|^{p\alpha}$
 uniformly  in $x\geq 0$ and $x_1,x_2\in[0,x]$.
 Similarly, for $i\in \{2,3,4,5\}$, by Proposition~\ref{WeakConvergenceProp02} and \ref{WeakConvergenceProp03} we also have 
 $$\mathbf{E} \Big[\big|M_i(x_1,x_2) \big|^{2p} \Big] \leq C\cdot (1+x)^{p\alpha}\cdot |x_2-x_1|^{p\alpha},$$ 
 uniformly in  $ x\geq 0$ and $x_1,x_2\in[0,x]$.
 The desired inequality follows by putting these estimates together.
 \qed

 {\sc Proof for Theorem~\ref{Thm.Regularity}.}
 By the Kolmogorov continuity theorem along with Lemma~\ref{MomentIncreX} and \ref{MomentIncreM}, the two processes $X_{\zeta,0}$ and $M$ are locally H\"older continuous of any order strictly less than $\alpha/2$. Then claim (1) holds.
 For the second claim, it suffices to consider the case of $p>1$.
 By the Garsia-Rodemich-Rumsey inequality; see Lemma~1.1 in \cite{GarsiaRodemichRumseyRosenblatt1970} with $\psi(u)=|u|^p$ and $p(u)=|u|^{q+1/p}$ for $q>1/p$, there exists a constant $C_{p,q}>0$ such that for any $x_2>x_1\geq 0$,
 \beqnn
 \big|M(x_2)-M(x_1)\big|^p\leq C_{p,q}\cdot |x_2-x_1|^{pq-1} \int_{x_1}^{x_2} dv \int_{x_1}^{x_2} \frac{|M(u)-M(v)|^p}{|u-v|^{pq+1}}du,\quad a.s.
 \eeqnn
 In particular, choosing $p>(\alpha/2-\kappa)^{-1}$ and $q=1/p+\kappa$ we have
 \beqnn
 \big\|M\big\|_{C^{0,\kappa}_x}^p
 \ar=\ar \sup_{0\leq x_1<x_2\leq x}\frac{|M(x_2)-M(x_1)|^p}{|x_2-x_1|^{p\kappa}}\cr
 \ar\leq\ar C_{p,q} \int_{0}^{x} dv \int_{0}^{x} \frac{|M(u)-M(v)|^p}{|u-v|^{p\kappa+2}}du,\quad a.s.
 \eeqnn
 From Lemma~\ref{MomentIncreM}, there exists a constant $C>0$ such that for any $x\geq 0$,
 \beqnn
 \mathbf{E}\Big[\big\|M\big\|_{C^{0,\kappa}_x}^p \Big]
 \ar\leq\ar  C   (1+x)^{p\alpha/2} \int_0^x dv \int_0^x |u-v|^{p \alpha/2 -p\kappa-2} du \leq C   \cdot   (1+x)^{p(\alpha-\kappa)}.
 \eeqnn
 Similarly, by Lemma~\ref{MomentIncreX} we also have $\mathbf{E}\big[\|X_{\zeta,0}\|_{C^{0,\kappa}_x}^p \big] \leq C (1+x)^{p(\alpha-\kappa)}$ uniformly in $x\geq 0$.
 Then claim (2) follows from these two upper bound estimates.  
 \qed

 \section{Laplace functionals and weak uniqueness}
 \label{NonLinearVolUnique}

 In this section, we firstly prove the affine representation of the Laplace functionals of $L^\xi_\zeta$ in two steps: (i)  on some finite interval, local solutions of (\ref{MainThm.Volterra}) exist uniquely and (\ref{AffineRepLap}) holds; (ii) this finite interval can be extended successfully to the whole positive real line. Then we prove the weak uniqueness holds for the SVE (\ref{MainThm.SVE}).
 As a preparation, we introduce several function spaces that will be used in the following proofs.
 For $T,J>0$,
 \begin{enumerate}
 	\item[$\bullet$] $\mathcal{A}_{T,J}$: the space of functions $f$ on $(0,\infty)$ satisfying $\sup_{t\in(0,T]}t^{1-\alpha}|f(t)|\leq J $;
 	\smallskip
 	
 	\item[$\bullet$] $\mathcal{B}_{T,J}$: the space of  functions $f$ on $\mathbb{R}_+$ satisfying $\|f\|_{L^\infty_T}\leq J$.
 \end{enumerate}
 Notice that  $\mathcal{A}_{T,J_1}\subset \mathcal{A}_{T,J_2}$ for any $J_2\geq J_1>0$.
 Let $\mathcal{A}_T:= \cup_{J>0} \mathcal{A}_{T,J}$, which collects all functions on $(0,\infty)$ satisfying that   $\sup_{x\in(0,T]}x^{1-\alpha}|f(x)|<\infty  $.
 It is obvious that $\mathcal{A}_{T}$ is decreasing in $T$ and $\mathcal{A}_\infty= \cap_{T>0} \mathcal{A}_T$.

 \subsection{Nonlinear Volterra equation}\label{NonlinearVolterra}
 We say a pair $(v_\lambda^g,T) \in \mathcal{A}_T \times (0,\infty)$ is a \textit{$\mathcal{A}$-local solution} of (\ref{MainThm.Volterra}) if $v_\lambda^g $ is continuous and satisfies (\ref{MainThm.Volterra}) on $(0,T]$.
 Moreover, for a pair $(v_\lambda^g,T_\lambda^g)$ with $T_\lambda^g\in (0,\infty]$ and a function $v_\lambda^g$ on $[0,T_\lambda^g)$, we say it is a \textit{$\mathcal{A}$-noncontinuable solution}\footnote{The terminology ``noncontinuable solution" comes from the theory of Volterra equations; see Chapter 12 in \cite{GripenbergLondenStaffans1990}. } of (\ref{MainThm.Volterra}) if for any  $T\in(0,T_\lambda^g)$, the pair $(v_\lambda^g,T)$ is a $\mathcal{A}$-local solution of (\ref{MainThm.Volterra}) and $\limsup_{x\to T_\lambda^g- }|v_\lambda^g(x)|= \infty$ if $T_\lambda^g-<\infty$.
 In particular,  if $T_\lambda^g=\infty$ the function $v_\lambda^g$ turns to be a  $\mathcal{A}$-global solution of (\ref{MainThm.Volterra}).
 In this section, we prove the existence and uniqueness of $\mathcal{A}$-noncontinuable solutions of (\ref{MainThm.Volterra}) with the help of the following technical estimates for the nonlinear operator $\mathcal{V}_\alpha$.
 The next useful inequality can be proved immediately by using the mean-value theorem
 \beqlb\label{ExponentialIneq}
 \big|(e^{-x}-1+x)-(e^{-z}-1+z) \big| \leq  e^{|x|\vee|z|} \cdot(|x|\vee|z|)\cdot |x-z|,\quad x,z\in\mathbb{R}.
 \eeqlb
 
 \begin{proposition}\label{UpBoundVf01}
 	There exists a constant $C> 0$ such that for any $T,J>0$, $f\in \mathcal{A}_{T,J}$ and $t\in(0,T]$,
 	\beqlb\label{eqn.UpBoundVf01}
 	\sup_{y\geq 0} \Big| \int_{(t-y)^+}^t f(r)dr \Big|  \leq  \frac{J}{\alpha}\cdot t^\alpha,
 	\eeqlb
 	and 
 	\beqlb \label{eqn.UpBoundVf02}
 	\int_0^\infty \Big|\int_{(t-y)^+}^t f(r) dr\Big|^2 \nu_\alpha(dy) \leq C \cdot J^2 \cdot t^{\alpha-1}.
 	\eeqlb
 \end{proposition}
 \proof The first desired inequality follows directly from a simple calculation. For the second one, we have
 \beqlb\label{eqn.101}
 \int_0^\infty \Big|\int_{(t-y)^+}^t f(r) dr\Big|^2 \nu_\alpha(dy)
 \ar\leq\ar
 \|f\|_{L^1_t}^2 \cdot \bar\nu_\alpha(t/2) +   \int_0^{t/2} \|f\|_{L^1_{[t-y,t]}}^2 \nu_\alpha(dy).
 \eeqlb
 It is easy to see that $\|f\|_{L^1_t}^2 \cdot \bar\nu_\alpha(t/2) \leq C \cdot J^2\cdot  t^{\alpha-1}$ uniformly in $t\in(0,T]$ and $T>0$.
 Since 
 $$\|f\|_{L^1_{[t-y,t]}} \leq J \cdot (t-y)^{\alpha-1}\cdot y \leq  2J\cdot t^{\alpha-1}\cdot y$$ 
 for any $ 0<y<t/2<t\leq T$, we have
 \beqlb\label{eqn.102}
 \int_0^{t/2}\|f\|_{L^1_{[t-y,t]}}^2  \nu_\alpha(dy)
 \ar\leq\ar 4J^2 t^{2\alpha-2} \int_0^{t/2} y^2 \nu_\alpha(dy)
 \leq C\cdot J^2   \cdot t^{\alpha-1},
 \eeqlb
 with $C>0$ depending only on $\alpha$ and $c$.
 The second desired inequality follows.
 \qed
 
 \begin{proposition}\label{UpBoundVf}
 	There exists a constant $C> 0$ such that for any $T,J>0$, $f\in \mathcal{A}_{T,J}$ and $t\in(0,T]$,
 	\beqnn
 	\big|\mathcal{V}_\alpha \circ f (t) \big| \leq C \cdot J^2   e^{\frac{J}{\alpha}\cdot t^\alpha}\cdot t^{\alpha-1}
 	\quad\mbox{and}\quad
 	\big|(\mathcal{V}_\alpha \circ f )* W'(t)\big| \leq C \cdot J^2  e^{\frac{J}{\alpha}\cdot t^\alpha}\cdot t^{2\alpha-1}.
 	\eeqnn
 \end{proposition}
 \proof
 By (\ref{eqn.UpBoundVf01}) and (\ref{ExponentialIneq}) with $z=0$,
 \beqnn
 \Big| \exp\Big\{-\int_{(t-y)^+}^t f(r)dr  \Big\}- 1+ \int_{(t-y)^+}^t f(r)dr  \Big|
 \leq  e^{ \frac{J}{\alpha}\cdot t^\alpha } \cdot \Big|\int_{(t-y)^+}^t f(r) dr\Big|^2.
 \eeqnn
 Plugging this back into (\ref{OperatorV}) and then using  (\ref{eqn.UpBoundVf02}), we have
 \beqnn
 \big|\mathcal{V}_\alpha\circ f(t)\big|
 \leq  e^{ \frac{J}{\alpha}\cdot t^\alpha } 	 \int_0^\infty \Big|\int_{(t-y)^+}^t f(r) dr\Big|^2 \nu_\alpha(dy) \leq C \cdot J^2   e^{\frac{J}{\alpha}\cdot t^\alpha}\cdot t^{\alpha-1}.
 \eeqnn
 From this and (\ref{UpperboundW}), we have  uniformly in $T,J\geq 0$, $f\in\mathcal{A}_{T,J}$ and $t\in[0,T]$, 
 \beqnn
 \big|(\mathcal{V}_\alpha  \circ f) * W'(t)\big| \ar\leq\ar C\cdot J^2 \cdot e^{\frac{J}{\alpha}\cdot t^\alpha }\cdot \int_0^t s^{\alpha-1} (t-s)^{\alpha-1}ds
 \leq C \cdot J^2 \cdot e^{\frac{J}{\alpha}\cdot t^\alpha }\cdot t^{2\alpha-1} . 
 \eeqnn
 \qed
 
 \begin{proposition}\label{ContractiveMap}
 	For $\theta\in(1,\frac{1}{1-\alpha} )$, there exists a constant $C> 0$ such that for any $T,J>0$ and $f_1,f_2\in \mathcal{A}_{T,J}$,
 	\beqnn
 	\big\|(\mathcal{V}_\alpha  \circ f_1-\mathcal{V}_\alpha  \circ f_2 )* W'\big\|_{L_T^\theta} \leq  C\cdot J T^{\alpha} e^{\frac{J}{\alpha}\cdot T^\alpha}\cdot\big\|f_1-f_2\big\|_{L^\theta_T}.
 	\eeqnn
 \end{proposition}
 \proof Let $\bar{f}:=f_1-f_2$. By (\ref{OperatorV}), (\ref{ExponentialIneq}) and  (\ref{eqn.UpBoundVf01}), 
 we have for any $t\in[0,T]$,
 \beqnn
 \lefteqn{\big|\mathcal{V}_\alpha  \circ f_1(t)-\mathcal{V}_\alpha  \circ f_2(t) \big|}\ar\ar\cr 
 \ar\ar\cr
 \ar\leq\ar e^{\frac{J}{\alpha}\cdot t^\alpha} \int_0^\infty\big(  \|f_1 \|_{L^1_{[(t-y)^+,t]}} \vee \| f_2\|_{L^1_{[(t-y)^+,t]}}\big) \|\bar f\|_{L^1_{[(t-y)^+,t]}}  \nu_\alpha(dy)
 \eeqnn
 and hence 
 \beqnn
 \big|(\mathcal{V}_\alpha \circ f_1-\mathcal{V}_\alpha \circ f_2 )* W' (t)\big| 
 \leq e^{\frac{J}{\alpha}\cdot t^\alpha}\cdot \big(I_1(t)+I_2(t)\big),
 \eeqnn
 where 
 \beqnn
 I_1(t) \ar:=\ar \int_0^t W'(t-s)ds \int_{s/2}^\infty \big(  \|f_1 \|_{L^1_{[(s-y)^+,s]}}\vee \| f_2\|_{L^1_{[(s-y)^+,s]}} \big) \|\bar f\|_{L^1_{[(s-y)^+,s]}}  \nu_\alpha(dy),\cr
 I_2(t) \ar:=\ar \int_0^t W'(t-s)ds \int_0^{s/2}  \big(  \|f_1 \|_{L^1_{[s-y,s]}}\vee \| f_2\|_{L^1_{[s-y,s]}} \big) \|\bar f\|_{L^1_{[s-y,s]}}  \nu_\alpha(dy).
 \eeqnn
 By Minkowski's inequality, 
 \beqnn
  \big\|(\mathcal{V}_\alpha \circ  f_1-\mathcal{V}_\alpha \circ  f_2 )* W' \big\|_{L_T^\theta}
  \leq  e^{\frac{J}{\alpha}T^\alpha}\cdot \big(\|I_1 \|_{L_T^\theta}+\|I_2 \|_{L_T^\theta} \big).
 \eeqnn
 Notice that
 $$I_1(t) \leq  \int_0^t W'(t-s) \cdot \big(  \|f_1 \|_{L^1_s}\vee \| f_2\|_{L^1_s} \big) \cdot \|\bar f\|_{L^1_s}  \cdot \bar\nu_\alpha(s/2) ds.$$
 By H\"older's inequality, we have
 $\| \bar{f} \|_{L^1_s}\leq \| \bar{f} \|_{L^\theta_s} \cdot s^{1-1/\theta}$.
 By (\ref{UpperboundW}) and (\ref{eqn.UpBoundVf01}), there exits a constant $C>0$ such that for any $T,J> 0$ and $t\in(0,T]$,
 \beqnn
 I_1(t) \ar\leq\ar C\cdot J\cdot\int_0^t (t-s)^{\alpha-1}s^{-1/\theta} \| \bar{f} \|_{L^\theta_s} ds
 \leq C \cdot J \cdot t^{\alpha-1/\theta}\cdot \| \bar{f} \|_{L^\theta_t}
 \eeqnn
 and
 \beqnn
 \|I_1\big\|_{L_T^\theta} \leq C\cdot J\cdot T^\alpha \cdot  \| \bar{f} \|_{L^\theta_T}.
 \eeqnn
 Noting that $\|f_1 \|_{L^1_{(s-y,s]}}\vee \| f_2\|_{L^1_{(s-y,s]}} \leq J\cdot(s-y)^{\alpha-1}\cdot y$ for any $y\in(0,s/2)$, we have
 \beqnn
 \lefteqn{\int_0^{s/2}  \big(  \|f_1 \|_{L^1_{[s-y,s]}}\vee \| f_2\|_{L^1_{[s-y,s]}} \big)\cdot \|\bar f\|_{L^1_{[s-y,s]}}  \nu_\alpha(dy)}\ar\ar\cr
 \ar\leq\ar  C\cdot J\int_0^{s/2} (s-y)^{\alpha-1} y^{-\alpha-1} \|\bar f\|_{L^1_{[s-y,s]}}   dy\cr
 \ar=\ar C\cdot J\int_0^{s/2} (s-y)^{\alpha-1} y^{-\alpha-1} \int_0^y |\bar{f}(s-r)|dr   dy\cr
 \ar=\ar  C\cdot J\int_0^{s/2} |\bar{f}(s-y)|dy \int_{y}^{s/2} (s-r)^{\alpha-1} r^{-\alpha-1} dr\cr
 \ar\leq\ar C\cdot J\cdot s^{\alpha-1} \int_0^{s/2} y^{-\alpha}|\bar{f}(s-y)|dy\cr
 \ar\leq\ar C\cdot J\cdot s^{\alpha-1} \int_0^{s} (s-y)^{-\alpha}|\bar{f}(y)|dy,
 \eeqnn
 for some constant $C$ independent of $s$, $J$, $T$ and $\bar{f}$.
 Here the two equalities follow from the change of variables and Fubini's theorem respectively.
 Plugging this into $I_2(t)$ and then using (\ref{UpperboundW}),
 \beqnn
 I_2(t) \ar\leq\ar C\cdot J \cdot \int_0^t (t-s)^{\alpha-1} s^{\alpha-1} ds \int_0^{s} (s-y)^{-\alpha}|\bar{f}(y)|dy\cr
 \ar=\ar C\cdot J \cdot \int_0^t  |\bar{f}(s)|ds \int_0^{t-s} (t-s-y)^{\alpha-1} (y+s)^{\alpha-1} y^{-\alpha} dy.
 \eeqnn
 This equality comes from Fubini's theorem and the change of variables.
 Let $\eta\in(0\vee \frac{1-\theta\alpha}{\theta-\theta\alpha},1)$.
 Since $(y+s)^{\alpha-1} \leq s^{(1-\eta)(\alpha-1)}\cdot y^{\eta(\alpha-1)}$ for any $y,s>0$, we have
 \beqnn
 I_2(t)\ar\leq\ar C\cdot J \cdot \int_0^t s^{(1-\eta)(\alpha-1)}  |\bar{f}(s)|ds \int_0^{t-s} (t-s-y)^{\alpha-1}  y^{\eta(\alpha-1)-\alpha} dy \cr
 \ar\leq\ar C\cdot J \cdot\int_0^t (t-s)^{\eta(\alpha-1) }s^{(1-\eta)(\alpha-1)}  |\bar{f}(s)|ds .
 \eeqnn
 By Young's convolution inequality and then H\"older's inequality,
 \beqnn
 \big\|I_2\big\|_{L_T^\theta} \ar\leq\ar
 C \cdot J\cdot \Big(\int_0^T t^{\theta\eta(\alpha-1) }dt \Big)^{1/\theta}  \cdot \int_0^T s^{(1-\eta)(\alpha-1)}  |\bar{f}(s)|ds \cr
 \ar\leq\ar C \cdot J\cdot T^{\eta(\alpha-1)+1/\theta} \cdot \Big(\int_0^T s^{(1-\eta)(\alpha-1)\frac{\theta}{\theta-1}} ds \Big)^{\frac{\theta-1}{\theta}} \|\bar{f}\|_{L^\theta_T}\cr
 \ar\ar\cr
 \ar\leq\ar C \cdot J\cdot T^\alpha \cdot \|\bar{f}\|_{L^\theta_T}
 \eeqnn
 with $C>0$ independent of $T,J$ and $\bar{f}$.
 The desired result follows by putting these estimates together.
 \qed
 
 The next two propositions can be proved similarly and their detailed proofs are omitted.
 \begin{proposition}\label{Prop.VL1}
 	There exists a constant $C> 0$ such that for any $T,J>0$, $f \in \mathcal{B}_{T,J}$ and $t\in(0,T]$,
 	\beqnn
 	\big|\mathcal{V}_\alpha\circ f(t)\big|   \leq C\cdot J^2 e^{J t}\cdot t^{1-\alpha}
 	\quad\mbox{and}\quad
 	\big|(\mathcal{V}_\alpha\circ f)*W'(t)\big| \leq C\cdot J^2 e^{Jt}\cdot t.
 	\eeqnn
 \end{proposition}

 \begin{proposition}\label{Prop.DVL1}
 	There exists a constant $C> 0$ such that for any $T,J>0$ and $f_1,f_2 \in \mathcal{B}_{T,J}$,
 	\beqnn
 	\big\|(\mathcal{V}_\alpha  \circ f_1-\mathcal{V}_\alpha  \circ f_2 )* W'\big\|_{L^1_T} \leq C\cdot JTe^{JT}\cdot \|f_1-f_2\|_{L^1_T}.
 	\eeqnn
 \end{proposition}

 \begin{lemma}\label{Lemma.NonlinearVE}
 	For each $\lambda \geq 0$ and $g\in L^\infty(\mathbb{R}_+;\mathbb{R}_+)$, the nonlinear Volterre equation (\ref{MainThm.Volterra}) has a unique $\mathcal{A}$-noncontinuable solution.
 \end{lemma}
 \proof By using Banach's fixed point theorem, we prove this lemma in the following three steps.
 
 {\it Step 1.} We first prove the existence of $\mathcal{A}$-local solutions near $0$.
 We consider the map $\mathcal{R}_0$ that acts on a locally integrable function $f$ on $\mathbb{R}_+$ according to
 $$
 \mathcal{R}_0 \circ f := \lambda W' +   (g - \mathcal{V}_\alpha \circ  f)*W'.
 $$
 Recall the constant  $\theta\in (1,1/(1-\alpha) )$.
 By (\ref{UpperboundW}), Proposition~\ref{UpBoundVf} and \ref{ContractiveMap}, there exists a constant $C_0>0$ such that for any $T,J>0$, $f_1,f_2\in \mathcal{A}_{T,J}$ and $t\in(0,T]$,
 \beqnn
 \big|\mathcal{R}_0\circ f_1(t) \big| \leq C_0\big( |\lambda | +\|g\|_{L^\infty}+ J^2\cdot T^\alpha  e^{\frac{J}{\alpha}T^\alpha} \big) \cdot t^{\alpha-1}
 \eeqnn
 and
 \beqnn
 \big\|\mathcal{R}_0\circ f_1-\mathcal{R}_0\circ  f_2  \big\|_{L_T^\theta} 
 \ar=\ar \big\|(\mathcal{V}_\alpha\circ  f_1-\mathcal{V}_\alpha\circ  f_2 )* W'\big\|_{L_T^\theta} \cr
 \ar\ar\cr
 \ar\leq\ar  C_0 J\cdot T^{\alpha}  e^{\frac{J}{\alpha}T^\alpha } \big\|f_1-f_2\big\|_{L^\theta_T} .
 \eeqnn
 Choosing $J_0>2C ( |\lambda| +\|g\|_{L^\infty} )$ and $T_0\in(0,1)$ such that $T_0^\alpha \cdot C_0\cdot J_0 e^{J_0/\alpha}< 1/2$, we have for any $f_1,f_2\in \mathcal{A}_{T_0,K_0}$ and $t\in (0,T_0]$,
 \beqnn
 |\mathcal{R}_0\circ  f_1(t)| \leq J_0\cdot t^{\alpha-1}
 \quad \mbox{and}\quad
 \big\|\mathcal{R}_0 \circ  f_1-\mathcal{R}_0\circ  f_2 \big\|_{L_{T_0}^\theta} < \big\|f_1-f_2 \big\|_{L^\theta_{T_0}}.
 \eeqnn
 Thus $\mathcal{R}_0$ is a contractive map from $\mathcal{A}_{T_0,K_0}$ to itself.
 It can be easily identify that $\mathcal{A}_{T_0,J_0}$ is a closed, bounded and convex subset in $L^\theta((0,T_0];\mathbb{R})$.
 By Banach's fixed point theorem, there exists a unique function $v_0 \in \mathcal{A}_{T_0,J_0}$  satisfying (\ref{MainThm.Volterra}), i.e., $v_0=\mathcal{R}_0\circ v_0$ almost everywhere on $(0,T_0]$.
 By the properties of convolution, the function $v_\lambda^g:= \mathcal{R}_0\circ  v_0 $ is continuous and equal to $v_0$ almost everywhere.
 Hence $\mathcal{R}_0 \circ v_\lambda^g = \mathcal{R}_0\circ  v_0 = v_\lambda^g$ pointwisely on $(0,T_0]$ and $(v_\lambda^g,T_0)$ is a $\mathcal{A}$-local solution of (\ref{MainThm.Volterra}).
 
 {\it Step 2.} We now extend the preceding $\mathcal{A}$-local solution into a $\mathcal{A}$-noncontinuable solution.
 Denote by $\mathcal{T}$ the collection of all $T>0$ such that (\ref{MainThm.Volterra}) has a $\mathcal{A}$-local solution on $(0,T]$.
 We assert that $\mathcal{T}$ is an open interval containing $(0,T_0]$.
 Indeed, for any $t_0\in \mathcal{T}$ and some $j_0>0$, assume that $v_\lambda^g \in \mathcal{A}_{t_0,j_0}$ is a $\mathcal{A}$-local solution of (\ref{MainThm.Volterra}).
 For $t\geq 0$, let
 \beqlb\label{FunH}
 H_1(t):= \lambda W'(t_0+t)+ W'*g(t_0+t)-\int_0^{t_0} \mathcal{V}_\alpha\circ v_\lambda^g(s) W'(t_0+t-s)ds.
 \eeqlb
 By (\ref{UpperboundW}), we have uniformly in $t\geq 0$,
 \beqnn
 \big|\lambda W'(t_0+t)\big|\leq C\cdot t_0^{\alpha-1}
 \quad \mbox{and}\quad
 \big|W'*g(t_0+t)\big|\leq C\cdot \|g\|_{L^\infty}\cdot (t_0+t)^{\alpha}.
 \eeqnn 
 From Proposition~\ref{UpBoundVf}, we also have $|\mathcal{V}_\alpha\circ v_\lambda^g(s) | \leq C \cdot  s^{\alpha-1}$ uniformly in $s\in (0,t_0]$ and hence 
 $$\int_0^{t_0} \mathcal{V}_\alpha\circ v_\lambda^g(s) W'(t_0+t-s)ds\leq C,$$ 
 uniformly in $t\geq 0$.
 Putting these  estimates together, there exists a constant $C_{H_1}>0$ such that $|H_1(t)| \leq C_{H_1}$ for any $t\in[0,1]$.
 We consider the map $\mathcal{R}_1$ acting on functions $f\in L^\infty(\mathbb{R}_+;\mathbb{R})$ by 
 $$\mathcal{R}_1\circ f:=  H_1 - W'* (\mathcal{V}_\alpha \circ f) .$$
 Recall $\mathcal{B}_{T,J}$ for $T,J>0$.
 From Proposition~\ref{Prop.VL1} and \ref{Prop.DVL1}, there exists a constant $C_1>0$ such that for any $T\in[0,1]$, $J>0$, $f_1,f_2\in \mathcal{B}_{T,J}$ and $t\in[0,T]$,
 \beqnn
 \big|\mathcal{R}_1\circ f_1(t) \big| \leq C_1 (C_{H_1}+ J^2e^{J\cdot T}\cdot T) 
 \eeqnn
 and
 \beqnn
 \big\|\mathcal{R}_1\circ f_1-\mathcal{R}_1\circ f_2 \big\|_{L^1_T} \leq C_1\cdot J T\cdot e^{2JT}\cdot \|f_1-f_2\|_{L^1_T}.
 \eeqnn
 Choosing $J_1>2C_1\cdot C_{H_1}$ and $T_1\in [0,1]$ such that $T_1\cdot C_1 J_1e^{2J_1}<1/2$, we have for any $f_1,f_2\in \mathcal{B}_{T_1,J_1}$,
 \beqnn
 \sup_{t\in[0,T_1]} \big|\mathcal{R}_1\circ f_1(t)\big| \leq J_1
 \quad \mbox{and}\quad
 \big\|\mathcal{R}_1\circ f_1-\mathcal{R}_1\circ f_2\big\|_{L^1_{T_1}}< \big\|f_1-f_2\big\|_{L^1_{T_1}}.
 \eeqnn
 Thus $\mathcal{R}_1$ is a contractive map from $\mathcal{B}_{T_1,J_1}$ to itself.
 Notice that $\mathcal{B}_{T_1,J_1}$ is a closed, bounded and convex subset of $L^1([0,T_1];\mathbb{R})$.
 Applying Banach's fixed point theorem again, there exists a unique function $v_1 \in \mathcal{B}_{T_1,J_1}$ satisfying that $v_1= \mathcal{R}_1\circ v_1$ almost everywhere. 
 Let 
 $$v_\lambda^g(t_0+t):= \mathcal{R}_1\circ v_1(t),\quad t\in[0,T_1].$$  
 One can verifies that $(v_\lambda^g,t_0+T_1)$ is a $\mathcal{A}$-local solution of (\ref{MainThm.Volterra}) and hence $t_0$ is an interior point of $\mathcal{T}$.
 
 Let $T_\lambda^g:= \sup \mathcal{T}$ and $v_\lambda^g $ be a continuous function on $(0,T_\lambda^g)$ satisfying that $(v_\lambda^g,T)$ is a $\mathcal{A}$-local solution of (\ref{MainThm.Volterra}) for any $T\in (0,T_\lambda^g)$.
 To assert that $(v_\lambda^g,T_\lambda^g)$ is a $\mathcal{A}$-noncontinuable solution, it remains to identify that
 \beqnn
 \limsup_{t\to T_\lambda^g-}\big|v_\lambda^g(t)\big|=\infty
 \quad\mbox{if}\quad
 T_\lambda^g<\infty.
 \eeqnn 
 Actually, if $\sup_{t\in[T_0, T_\lambda^g-)}|v_\lambda^g(t)|\leq J_2$ for some constant $J_2>0$, then 
 $$\sup_{t\in(0,T_\lambda^g)} t^{1-\alpha}\big|v_\lambda^g (t)\big|<\infty.$$
 Let $H_2$ be a function on $\mathbb{R}_+$ defined as in (\ref{FunH}) with $t_0$ replaced by $T_\lambda^g$.
 Let $\mathcal{R}_2$ be a map acting on functions 
 $f\in L^\infty(\mathbb{R}_+;\mathbb{R})$ by 
 $$\mathcal{R}_2\circ f:=  H_2  + W'* (\mathcal{V}_\alpha \circ f).$$ 
 Similarly as in the previous paragraph, there exist constants $C_{H_2}>0$, $T_2\in[0,1]$ and $J_2>0$ such that $|H_2(t)|\leq C_{H_2}$ for any $t\in [0,1]$ and $\mathcal{R}_2$ is a contractive map from $\mathcal{B}_{T_2,J_2}$ to itself.
 By Banach's fixed point theorem again, there exists a unique function $v_2 \in \mathcal{B}_{T_2,J_2}$ satisfying that $v_2= \mathcal{R}_2\circ v_2$ almost everywhere.   
 Let 
 $$v_\lambda^g(T_\lambda^g+t):= \mathcal{R}_1\circ v_2(t),\quad t\in[0,T_2].$$
 Then $(v_\lambda^g,T_\lambda^g+T_2)$ is a $\mathcal{A}$-local solution of (\ref{MainThm.Volterra}) and $T_\lambda^g+T_2\in \mathcal{T}$, which contradicts the assumption that $T_\lambda^g= \sup\mathcal{T}$.
 Consequently,   $(v_\lambda^g,T_\lambda^g)$ is a $\mathcal{A}$-noncontinuable solution of (\ref{MainThm.Volterra}).

 {\it Step 3.} We prove the uniqueness.
 Assume that $(v_\lambda^g,T_\lambda^g)$ and $(\hat{v}_\lambda^g,\hat{T}_\lambda^g)$ are two $\mathcal{A}$-noncontinuable solutions of (\ref{MainThm.Volterra}) with $T_\lambda^g  \leq \hat{T}_\lambda^g$. 
 Similarly as in Step~1, there exist two constants $T_0\in(0,1)$ and $J_0>0$ such that
 $v_\lambda^g,\hat{v}_\lambda^g \in \mathcal{A}_{T_0,J_0}$  and the map $\mathcal{R}_0$ is  contractive from $\mathcal{A}_{T_0,J_0}$ to itself. 
 Then  Banach's fixed point theorem induces that $$\big\|v_\lambda^g-\hat{v}_\lambda^g \big\|_{L^\theta_{T_0}}=0.$$ 
 Their continuity yields that $v_\lambda^g=\hat{v}_\lambda^g$ on $(0,T_0]$.
 Similarly as in Step~2, let $t_0:=\inf\{ t>0: v_\lambda^g(t)\neq\hat{v}_\lambda^g(t) \}$, $H_1$ be the function defined by (\ref{FunH}), $v_{\lambda,1}^g(t)= v_\lambda^g(t_0+t)$ and $\hat{v}_{\lambda,1}^g(t)= \hat{v}_\lambda^g(t_0+t)$ for $t\in[0, T_\lambda^g-t_0)$.
 We also can find some constants $T_1\in(0,1)\cap[0, T_\lambda^g-t_0)$ and $J_1>0$ such that $v_{\lambda,1}^g, \hat{v}_{\lambda,1}^g \in \mathcal{B}_{T_1,J_1}$ and the map $\mathcal{R}_1$ is contractive from $\mathcal{B}_{T_1,J_1}$ to itself. Again,  Banach's fixed point theorem induces that $\| v_{\lambda,1}^g- \hat{v}_{\lambda,1}^g \|_{L^1_{T_1}}=0$. 
 Their continuity yields that $v_\lambda^g=\hat{v}_\lambda^g$ on $(0,t_0+T_1]$, which contracts the definition of $t_0$. Hence the uniqueness holds.
 \qed

 \subsection{Laplace functionals}
 For convenience, we assume that the process $L^\xi_\zeta$, the PRM $\widetilde{N}_\alpha(ds,dy,dz)$ are defined on a filtrated probability basis $(\Omega, \mathscr{G},\mathscr{G}_r,\mathbf{P})$ satisfying the general hypotheses and $N_0(dy,dz)$ is $\mathscr{G}_0$-measurable\footnote{Otherwise, by the definition of weak solutions of (\ref{MainThm.SVE}); see footnote~\ref{Def.WeakSolution}, a realization of $L^\xi_L$, $N_0(dy,dz)$ and $\widetilde{N}_\alpha(ds,dy,dz)$ can be found on some filtrated probability basis such that (\ref{MainThm.SVE}) holds. }.
 For each $\lambda \geq 0$ and $g\in L^\infty(\mathbb{R}_+;\mathbb{R}_+)$, we assume that $(v_\lambda^g,T)$ is a $\mathcal{A}$-local solution of (\ref{MainThm.Volterra}), i.e., $T>0$ and $v_\lambda^g \in \mathcal{A}_{T,J}$ for some $J>0$.
 For $x,r\geq 0$, conditioned on  $\mathscr{G}_r$ we take expectations on both sides of (\ref{MainThm.SVE}) and get
 \beqlb\label{eqn.ConExp}
 \mathbf{E} \big[L^\xi_\zeta(x)\,\big|\,\mathscr{G}_r\big] 
 \ar=\ar \int_0^\infty  \int_0^\zeta \nabla_y W(x) N_0( dy,dz)  \cr
 \ar\ar + \int_0^{x \wedge r} \int_0^\infty \int_0^{L^\xi_\zeta(s)}  \nabla_y W(x-s)\widetilde{N}_\alpha(ds,dy,dz).
 \eeqlb

 \begin{proposition} \label{MartingaleReProp}
 	For any $x\in[0,T]$, the following hold:
 	\begin{enumerate}
 		\item[(1)] The random variable $Y_x(x):= \lambda L^\xi_\zeta(x)+ (g -\mathcal{V}_\alpha \circ v_\lambda^g) *L^\xi_\zeta(x)$ is integrable, i.e., $\mathbf{E}\big[|Y_x(x)|\big]<\infty$.
 		
 		\item[(2)] The Doob martingale $\big\{Y_x(t):= \mathbf{E}\big[Y_x(x)|\mathscr{G}_t\big]:t\in[0,x]\big\}$ has the representation
 		\beqnn
 		Y_x(t) \ar=\ar \int_0^\infty  \int_0^\zeta \int_{(x-y)^+}^x  v_\lambda^g(s) ds N_0( dy,dz) \cr
 		\ar\ar +  \int_0^{t} \int_0^\infty \int_0^{L^\xi_\zeta(s)} \int_{(x-s-y)^+}^{x-s} v_\lambda^g(r)dr \widetilde{N}_\alpha(ds,dy,dz).
 		\eeqnn
 	\end{enumerate}
 \end{proposition}
 \proof
 By Proposition~\ref{UpBoundVf}, we have $\big\| \mathcal{V}_\alpha \circ v_\lambda^g\big\|_{L^1_x} \leq C \cdot x^\alpha$.
 Moreover, by (\ref{BoundFirstMom}) we have
 $$
 \mathbf{E}\big[|Y_x(x)|\big] \leq   \zeta  \cdot \big( \lambda + \|g\|_{L^1_x} +\big\| \mathcal{V}_\alpha \circ v_\lambda^g\big\|_{L^1_x} \big) <\infty
 $$
 and
 \beqlb\label{eqn.401}
 Y_x(0)
 \ar=\ar \int_0^\infty  \int_0^\zeta \lambda \nabla_y W(x)   N_0( dy,dz) \cr
 \ar\ar + \int_0^\infty  \int_0^\zeta \nabla_y W(\cdot)   N_0( dy,dz) * (g - \mathcal{V}_\alpha \circ v_\lambda^g)  (x).
 \eeqlb
 By using Proposition~\ref{UpBoundVf} and (\ref{BoundFirstMom}) again, we have for any $\epsilon>0$, 
 \beqnn
 \lefteqn{ \Big|\int_0^x  (g - \mathcal{V}_\alpha \circ v_\lambda^g)  (x-s) \int_0^\epsilon \nabla_y W(s) \bar\nu_\alpha(y)dy ds \Big|}\ar\ar\cr
 \ar\leq\ar \big( \|g\|_{L^1_x} +\big\| \mathcal{V}_\alpha \circ v_\lambda^g\big\|_{L^1_x} \big) \sup_{s\in[0,x]}  \int_0^\epsilon \nabla_y W(s) \bar\nu_\alpha(y)dy  ,
 \eeqnn
 which is finite and goes to $0$ as $\epsilon\to 0 + $.
 Moreover, by Fubini's theorem,
 \beqnn
 \int_0^x  (g - \mathcal{V}_\alpha \circ v_\lambda^g)(x-s) \nabla_y W(s) ds
 \ar=\ar \int_0^x  (g - \mathcal{V}_\alpha \circ v_\lambda^g)(x-s) \int_{(s-y)^+}^s W'(r)dr ds\cr
 \ar=\ar \int_{(x-y)^+}^x   (g - \mathcal{V}_\alpha \circ v_\lambda^g) *W'(s)ds ,
 \eeqnn
 which can be bounded by $C(1\wedge y)$ uniformly in $y>0$; see Proposition~\ref{UpBoundVf}.
 Hence for any $\epsilon>0$, 
 \beqnn
 \int_0^\epsilon \Big|\int_{(x-y)^+}^x   (g - \mathcal{V}_\alpha \circ v_\lambda^g) *W'(s)ds\Big| \bar\nu_\alpha(y)dy \leq C \cdot \epsilon^{\alpha} ,
 \eeqnn
 which vanishes as $\epsilon\to0 + $. 
 By these and the stochastic Fubini theorem; see Theorem~\ref{StoFubiniThm},
 the stochastic integral
 \beqnn
 \int_0^\infty  \int_0^\zeta \int_{(x-y)^+}^x   (g - \mathcal{V}_\alpha \circ v_\lambda^g) *W'(s)ds  N_0( dy,dz)
 \eeqnn
 is well defined and equal almost surely to the second stochastic integral on the right side of (\ref{eqn.401}).
 By (\ref{MainThm.Volterra}),
 \beqnn
 \lefteqn{ \lambda \nabla_y W(x) + \int_{(x-y)^+}^x   (g - \mathcal{V}_\alpha \circ v_\lambda^g) *W'(s)ds}\ar\ar\cr
 \ar = \ar \int_{(x-y)^+}^x  \big(\lambda W'(s) + (g - \mathcal{V}_\alpha \circ v_\lambda^g) *W'(s) \big) ds =  \int_{(x-y)^+}^x v_\lambda^g(s) ds .
 \eeqnn
 Plugging this back into (\ref{eqn.401}), we have
 \beqlb\label{eqn.402}
 Y_x(0)  \ar=\ar \int_0^\infty  \int_0^\zeta \int_{(x-y)^+}^x   v_\lambda^g(s) ds N_0( dy,dz).
 \eeqlb
 By (\ref{eqn.ConExp}), we have for $t\in[0,x]$,
 \beqlb\label{MartingaleRep.01}
 Y_x(t) \ar=\ar \int_0^\infty  \int_0^\zeta \int_{(x-y)^+}^x   v_\lambda^g(s) ds N_0( dy,dz) \cr
 \ar\ar + \int_0^{t} \int_0^\infty \int_0^{L^\xi_\zeta(s)} \lambda \nabla_y W(x-s)  \widetilde{N}_\alpha(ds,dy,dz) \cr
 \ar\ar  + \int_0^x \big( g - \mathcal{V}_\alpha \circ v_\lambda^g \big)(x-r)dr \int_0^{r\wedge t} \int_0^\infty \int_0^{L^\xi_\zeta(s)} \nabla_y W(r-s)  \widetilde{N}_\alpha(ds,dy,dz).
 \eeqlb
 Proposition~\ref{UpBoundVf}, together with the assumption that $v_\lambda^g \in \mathcal{A}_T$, implies that $g -\mathcal{V}_\alpha \circ v_\lambda^g \in \mathcal{A}_T$. Then there exists a constant $C>0$ such that $|g(t) -\mathcal{V}_\alpha \circ v_\lambda^g(t)|\leq C\cdot K(t)$ for any $t\in(0,T]$.
 By Proposition~\ref{Prop.D4},
 \beqlb \label{eqn.107}
 \lefteqn{\int_0^x  \big( g - \mathcal{V}_\alpha \circ v_\lambda^g \big)(x-r) dr \Big| \int_0^r ds\int_0^\epsilon |\nabla_y W(r-s) |^{2}\nu_\alpha(dy) \Big|^{1/2}}\qquad   \ar\ar\cr
 \ar\leq\ar C\int_0^x K(r)dr\cdot  \Big| \int_0^x ds\int_0^\epsilon |\nabla_y W(s) |^{2}\nu_\alpha(dy) \Big|^{1/2} \cr
 \ar\leq\ar C\cdot x^{\alpha}\cdot  \Big| \int_0^x ds\int_0^\epsilon |\nabla_y W(s) |^{2}\nu_\alpha(dy) \Big|^{1/2}  ,
 \eeqlb
 which goes to $0$ as $\epsilon\to 0+$. Moreover, by the change of variables,
 \beqnn
 \lefteqn{\int_0^x \big( g-\mathcal{V}_\alpha \circ v_\lambda^g \big)(x-r) \cdot \mathbf{1}_{\{0\leq   s<r\wedge t \}}\cdot \nabla_y W(r-s) dr}\ar\ar\cr
 \ar=\ar \int_s^x \big( g - \mathcal{V}_\alpha v_z^g \big)(x-r)  \cdot \nabla_y W(r-s) dr \cdot \mathbf{1}_{\{0\leq s< t\}}\cr
 \ar\ar\cr
 \ar=\ar\big( g -\mathcal{V}_\alpha \circ v_\lambda^g \big)* \nabla_y W(x-s)   \cdot \mathbf{1}_{\{0\leq s< t\}} .
 \eeqnn
 By Proposition~\ref{Prop.B5}, we have for $\epsilon >0$,
 \beqnn
 \lefteqn{\int_0^x ds \int_0^\epsilon | \big( g -\mathcal{V}_\alpha \circ v_\lambda^g \big)* \nabla_y W(x-s)|^2 \nu_\alpha(dy)}\ar\ar\cr
 \ar \leq\ar C \int_0^x ds \int_0^\epsilon \frac{|K* \nabla_y W(x-s)|^2}{y^{\alpha+2}}dy,
 \eeqnn
 which is finite and goes to $0$ as $\epsilon\to 0+$.
 From this, (\ref{eqn.107}) and the stochastic Fubini theorem; see Theorem~\ref{StoFubiniThm},
 the stochastic integral
 \beqnn
 \int_0^t \int_0^\infty \int_0^{L^\xi_\zeta(s)} \big( g -\mathcal{V}_\alpha \circ v_\lambda^g \big)* \nabla_y W(x-s)   \widetilde{N}_\alpha(ds,dy,dz)
 \eeqnn
 is well defined and equal almost surely to the last stochastic integral on the right side of (\ref{MartingaleRep.01}).
 Moreover, by Fubini's theorem we have for any $s\in[0,T]$ and $y>0$,
 \beqnn
 \lambda\nabla_yW(s) +\big( g- \mathcal{V}_\alpha\circ v_\lambda^g \big)*  \nabla_y W(s) 
 \ar=\ar \int_{(s-y)^+}^s \big( \lambda W'(r)+   ( g - \mathcal{V}_\alpha\circ v_\lambda^g )*  W'(r)\big)dr \cr
 \ar=\ar \int_{ (s-y)^+ }^s  v_\lambda^g(r)dr.
 \eeqnn
 Consequently, the sum of the last two stochastic integrals on the right side of (\ref{MartingaleRep.01}) can be replaced by
 \beqnn
 \int_0^{ t} \int_0^\infty \int_0^{L^\xi_\zeta(s)}
 \int_{ (x-s-y)^+ }^{x-s} v_\lambda^g(r)dr \widetilde{N}_\alpha(ds,dy,dz)
 \eeqnn
 and claim (2) holds.
 \qed
 
 Associated  to $v_\lambda^g $, we define a stochastic process $Z_x:=\big\{Z_x(t):t\in[0,x]\big\}$ by
 \beqnn
 Z_x(t)\ar:=\ar \mathbf{E}\big[\lambda L^\xi_\zeta(x)+ g*L^\xi_\zeta(x) \big|\mathscr{G}_t  \big] - \int_t^x \mathcal{V}_\alpha\circ v_\lambda^g(x-s)\mathbf{E}\big[L^\xi_\zeta(s) \big|\mathscr{G}_t  \big] ds.
 \eeqnn
 By Proposition~\ref{MartingaleReProp}(2) and (\ref{eqn.402}), the process $Z_x$ also has the following representation
 \beqnn
 Z_x(t) \ar=\ar Y_x(t)+ \int_0^t \mathcal{V}_\alpha\circ v_\lambda^g (x-s) L^\xi_\zeta(s) ds \cr
 \ar=\ar Y_x(0) + \int_0^t \mathcal{V}_\alpha \circ v_\lambda^g (x-s) L^\xi_\zeta(s) ds \cr
 \ar\ar + \int_0^{ t} \int_0^\infty \int_0^{L^\xi_\zeta(s)} \int_{(x-s-y)^+}^{x-s}v_\lambda^g(r) dr \widetilde{N}_\alpha(ds,dy,dz).
 \eeqnn
 Thus $Z_x$ is a $(\mathscr{G}_t)$-semimartingale.
 Applying It\^o's formula to $\exp\{-Z_x(t) \}$ and then using (\ref{MainThm.Volterra}),  
 \beqlb\label{eqn.600}
 e^{-Z_x(t) }= e^{- Y_x(0) } + M_x(t),\quad t\in[0,x],
 \eeqlb
 where $M_x:=\{ M_x(t):t\in[0,x] \}$ is a $(\mathscr{G}_r)$-local martingale and 
 \beqnn
 M_x(t) := \int_0^t \int_0^\infty \int_0^{L^\xi_\zeta(s)} e^{-Z_x(s)} \Big(  \exp\Big\{  - \int_{(x-s-y)^+}^{x-s}v_\lambda^g(r) dr \Big\}-1 \Big) \widetilde{N}_\alpha(ds,dy,dz).
 \eeqnn
 In the next lemma, we prove the martingality of $e^{-Z_x}:=\{e^{-Z_x(t) }: t\in[0,x]\}$ and the equality (\ref{AffineRepLap}) by using the method developed in \cite[Lemma~6.3]{Jaber2021} and \cite[Lemma~7.3]{JaberLarssonPulido2019}.
 
 \begin{lemma}\label{Lemma.LapFun}
 	For each $x\in[0,T]$, the local-martingale $e^{-Z_x}$ is a true  $(\mathscr{G}_r)$-martingale and (\ref{AffineRepLap}) holds.
 \end{lemma}
 \proof For each $t\geq 0$, define
 \beqnn
 U_x(t):=  \int_0^t \int_0^\infty \int_0^{L^\xi_\zeta(s)}  \Big(  \exp\Big\{  - \int_{(x-s-y)^+}^{x-s}v_\lambda^g(r) dr \Big\}-1 \Big) \widetilde{N}_\alpha(ds,dy,dz),
 \eeqnn
 which is a uniformly square integrable $(\mathscr{G}_r)$-martingale, i.e., by the Burkholder-Davis-Gundy inequality, (\ref{BoundFirstMom}), Proposition~\ref{UpBoundVf01} and the change of variables,
 \beqnn
 \lefteqn{\sup_{t\geq 0}\mathbf{E}\Big[\big|U_x(t)\big|^2\Big] }\ar\ar\cr
 \ar\leq\ar \sup_{t\geq 0} \int_0^t\mathbf{E} \big[L^\xi_\zeta(s)\big]  ds  \int_0^\infty \Big(  \exp\Big\{  - \int_{(x-s-y)^+}^{x-s}v_\lambda^g(r) dr \Big\}-1 \Big)^2\nu_\alpha(dy) \cr
 \ar\leq \ar C   \int_0^x ds\int_0^\infty  \Big|  \int_{(x-s-y)^+}^{x-s}v_\lambda^g(r) dr\Big|^2 \nu_\alpha(dy)\cr
 \ar=\ar C   \int_0^x  s^{\alpha-1}ds \cr
 \ar=\ar \frac{C}{\alpha}\cdot  x^\alpha.
 \eeqnn
 Denote by $\mathcal{E}_{U_x}:= \{\mathcal{E}_{U_x}(t): t\geq 0 \}$ the Dol\'ean-Dade exponential of $U_x$. By It\^o's formula,
 \beqnn
 \mathcal{E}_{U_x}(t)
 \ar=\ar\exp\Big\{ - \int_0^t \mathcal{V}_\alpha\circ v_\lambda^g (x-s) L^\xi_\zeta(s) ds  \cr
 \ar\ar\ \  \quad\quad - \int_0^t \int_0^\infty \int_0^{L^\xi_\zeta(s)}   \int_{(x-s-y)^+}^{x-s}v_\lambda^g(r) dr  \widetilde{N}_\alpha(ds,dy,dz) \Big\}.
 \eeqnn
 Notice that $ \mathcal{E}_{U_x}$ is a non-negative local martingale and hence a supermartingale.
 By Fatou's lemma, we have $\mathbf{E}[\mathcal{E}_{U_x}(t)] \leq 1$ and hence it suffices to identify that $\mathbf{E}[\mathcal{E}_{U_x}(t)]=1$ for any $t\geq 0$.
 For each $t_0\geq 0$ and $n\geq 1$, let 
 \beqnn
 \tau_n:=\inf\big\{ t\geq 0: L_\zeta^\xi(t)\geq n \big\} \wedge t_0
 \quad\mbox{and}\quad 
  \mathcal{E}_{U_x}^n(t):=  \mathcal{E}_{U_x}(\tau_n \wedge t),\quad t\geq 0. 
 \eeqnn 
 By the inequality $|1+(z-1)e^{z} | \leq z^2e^{z} $ for any $z\in\mathbb{R}$ and Proposition~\ref{UpBoundVf01}, there exists a constant $C>0$ such that for any $t\geq 0$,
 \beqnn
 \ar\ar \int_0^t \mathbf{1}_{\{s\leq \tau_n\}} L_\zeta^\xi(s) ds \int_0^\infty \Big| 1-\Big( 1+\int_{(x-s-y)^+}^{x-s}v_\lambda^g(r) dr \Big) \times  \cr
 	\ar\ar\qquad\qquad \qquad \qquad \qquad  \times   \exp\Big\{- \int_{(x-s-y)^+}^{x-s}v_\lambda^g(r) dr\Big\}  \Big| \nu_\alpha(dy) \cr
 \ar\leq\ar n \int_0^x ds \int_0^\infty \Big(  \int_{(x-s-y)^+}^{x-s}v_\lambda^g(r) dr\Big)^2 \exp\Big\{- \int_{(x-s-y)^+}^{x-s}v_\lambda^g(r) dr\Big\} \nu_\alpha(dy) \cr
 \ar\leq\ar C\cdot  \int_0^x s^{\alpha-1}ds  \leq C x^\alpha.
 \eeqnn
 By Theorem~IV.3 in \cite{LepingleMemin1978} with 
 $$y(s,z)= \mathbf{1}_{\{s\leq \tau_n\}}\cdot \bigg(\exp\Big\{- \int_{(x-s-y)^+}^{x-s}v_\lambda^g(r) dr \Big\}-1\bigg),$$
  the process $\mathcal{E}_{U_x}^n$ is a martingale for each $n\geq 1$. Thus
 \beqnn
 1= \mathbf{E}\big[\mathcal{E}_{U_x}^n(t_0)\big]
 \ar=\ar \mathbf{E}\big[\mathcal{E}_{U_x}^n(t_0); \tau_n = t_0\big] + \mathbf{E}\big[\mathcal{E}_{U_x}^n(t_0);\tau_n<t_0\big] \cr
 \ar\ar\cr
 \ar=\ar \mathbf{E}\big[\mathcal{E}_{U_x}(t_0);\tau_n= t_0\big] + \mathbf{E}\big[\mathcal{E}_{U_x}^n(t_0);\tau_n<t_0\big] .
 \eeqnn
 By the monotone convergence theorem and the fact that $\tau_n\overset{\rm a.s.}\to t_0$  as $n\to\infty$, we have 
 $$ \mathbf{E}\big[\mathcal{E}_{U_x}(t_0); \tau_n=t_0\big]\to  \mathbf{E}\big[\mathcal{E}_{U_x}(t_0)\big] .$$
 Thus it suffices to prove that $ \mathbf{E}\big[\mathcal{E}_{U_x}^n(t_0);\tau_n<t_0\big] \to 0$ as $n\to\infty$.
 Associate with the martingale $\mathcal{E}_{U_x}^n$, we define a probability law  $\mathbf{Q}_x^n$ on $(\Omega,\mathscr{G},\mathscr{G}_r)$ by 
 \beqnn
 \frac{d  \mathbf{Q}_x^n}{d \mathbf{P}}= \mathcal{E}_{U_x}^n(\tau_n). 
 \eeqnn
 Since $\mathcal{E}_{U_x}^n(0)\overset{a.s.}=1$, the PRM $N_0(dy,dz)$ is $\mathscr{G}_0$-measurable and has the same distribution under $\mathbf{P}$ and $\mathbf{Q}_x^n$.
 By Girsanov's theorem for random measure; see Theorem~3.17 in \cite[p.170]{JacodShiryaev2003}, the PRM $N_\alpha(ds,dy,dz)$ is a random point measure under $\mathbf{Q}_x^n$ with intensity
 \beqnn
 \mathbf{1}_{\{ s\leq \tau_n \}}\cdot \exp\Big\{- \int_{(x-s-y)^+}^{x-s}v_\lambda^g(r) dr\Big\}  ds \cdot \nu_\alpha(dy) \cdot dz,
 \eeqnn
 and the SVE (\ref{MainThm.SVE}) is equal in distribution to the following SVE under  $\mathbf{Q}_x^n$,
 \beqnn
 L^{\xi}_{\zeta}(t)\ar=\ar   X_{\zeta,0}(t) + B^n(t)+ M^n(t),\quad t\geq 0,  
 \eeqnn
 where $X_{\zeta,0}$ is defined as in (\ref{MartingaleX}), $ M^n$ is defined as in (\ref{MartingaleM}) with $X_{\zeta} $ replaced by $L^{\xi}_{\zeta} $ and
 \beqnn
 B^n(t)\ar:=\ar \int_0^t   \mathbf{1}_{\{ s\leq \tau_n \}} L^\xi_{\zeta}(s)ds  \int_0^\infty  \nabla_y W(t-s) \times \cr
 \ar\ar \times \Big(  \exp\Big\{- \int_{(x-s-y)^+}^{x-s}v_\lambda^g(r) dr\Big\} -1 \Big) \nu_\alpha(dy).
 \eeqnn
 Proposition~\ref{UpBoundVf01} implies that 
 $$\exp \Big\{- \int_{(x-s-y)^+}^{x-s}v_\lambda^g(r) dr \Big\} $$ 
 is uniformly bounded in $s,y\geq 0$.
 Similarly as in the proofs of Lemma~\ref{Lemma.501}-\ref{MomentIncreM}, there exists a constant $C>0$ such  that for any $n\geq 1$, $t>0$, $t_1,t_2\in[0,t]$ and $p> 1$,
 \beqlb\label{eqn.418}
 \mathbf{E}^{\mathbf{Q}_x^n} \Big[\big|X_{\zeta,0}(t)\big|^p + \big|L^{\xi}_{\zeta}(t)\big|^p\Big]< C\cdot (1+t)^{p\alpha} 
 \eeqlb
 and
 \beqnn
 \mathbf{E}^{\mathbf{Q}_x^n} \Big[\big|X_{\zeta,0}(t_1)-X_{\zeta,0}(t_2)\big|^p + \big| M^n (t_1)- M^n (t_2)\big|^p\Big]\leq C \cdot (1+t)^{p\alpha}\cdot |t_1-t_2|^{p\alpha}.
 \eeqnn
 Here $\mathbf{E}^{\mathbf{Q}_x^n}$ is the expectation under $ \mathbf{Q}_x^n$.
 Together with these estimates, an argument similar to that in the proof of  Theorem~\ref{Thm.Regularity} implies that for any $\kappa\in(0,\alpha/2)$, both $ X_{\zeta,0} $ and $ M^n  $ are locally $\kappa$-H\"older continuous under $\mathbf{Q}_n^x$ and the H\"older coefficient has finite moments of all orders.
 Like the argument before Corollary~\ref{Cor.Maxi}, we have
 \beqlb\label{eqn.612}
 \sup_{n\geq 1}\mathbf{E}^{\mathbf{Q}_x^n}\Big[\sup_{s\in[0,t]}| X_{\zeta,0}(s)|^p\Big]+ \sup_{n\geq 1}\mathbf{E}^{\mathbf{Q}_x^n}\Big[\sup_{s\in[0,t]}| M^n(s)|^p \Big]<\infty,
 \eeqlb
 for any $t,p\geq 0$.  
 Noting that $ \int_{(x-s-y)^+}^{x-s}v_\lambda^g(r) dr=0$ when $s\geq x$, we have for any $t\geq 0$,
 \beqnn
 B^n(t) 
 \ar =\ar \int_0^{t\wedge x}  \mathbf{1}_{\{ s\leq \tau_n \}} L^\xi_{\zeta}(s)ds  \int_0^\infty  \nabla_y W(t-s) \times \cr
 \ar\ar\times  \Big(  \exp\Big\{- \int_{(x-s-y)^+}^{x-s}v_\lambda^g(r) dr\Big\} -1 \Big) \nu_\alpha(dy).
 \eeqnn
 By H\"older's inequality and Proposition~\ref{UpBoundVf01}, the preceding inner integral can bounded by $C\cdot((x-s)^+ \cdot (t-s)^+)^{(\alpha-1)/2}$ uniformly in $t,s,y>0$.
 By (\ref{eqn.418}) and H\"older's inequality, we have for $p>1$,
 \beqnn
 \lefteqn{\mathbf{E}^{\mathbf{Q}_x^n}\Big[ \sup_{t\geq 0} \big| B^n(t)\big|^p \Big]}\ar\ar\cr
 \ar\leq\ar  C \int_0^{x\wedge t \wedge t_0} \mathbf{E}^{\mathbf{Q}_x^n}\Big[\big|L^\xi_{\zeta}(s)\big|^p\Big]  (x-s)^{\alpha-1}ds  \Big(\int_0^x   (x-r)^{\alpha-1} dr\Big)^{p-1} \leq C\cdot x^{p\alpha},
 \eeqnn
 for some constant $C>0$ independent of $n\geq 1$. This, together with (\ref{eqn.612}), yields that  
 $$\mathbf{E}^{\mathbf{Q}_x^n}\Big[\sup_{s\in[0,t_0]}\big|L^{\xi}_{\zeta}(s)\big|^p \Big]\leq C ,$$ 
 uniformly in $n\geq 1$. 
 By the definition of $\tau_n$ and Chebyshev's inequality,
 \beqnn
 \mathbf{E}\big[\mathcal{E}_{U_x}^n(\tau_n); \tau_n<t_0 \big]  =\mathbf{Q}_x^n \big(\tau_n<t_0 \big) 
 \ar=\ar \mathbf{Q}_x^n\Big(\sup_{s\in[0,t_0]}L_{\zeta}^{\xi}(s)\geq n\Big)\cr \ar\leq\ar \frac{1}{n}\mathbf{E}^{\mathbf{Q}_x^n}\Big[\sup_{s\in[0,t_0]}L_{\zeta}^{\xi}(s)\Big],
 \eeqnn
 which vanishes as $n\to\infty$. Hence $\mathbf{E}[\mathcal{E}_{U_x}(t_0)]=1$ and $\mathcal{E}_{U_x}$ is a $(\mathscr{G}_r)$-martingale under $\mathbf{P}$.
 By Proposition~\ref{UpBoundVf01},
 \beqnn
 \int_0^\infty \Big| 1- \exp\Big\{  - \int_{(x-y)^+}^{x}v_\lambda^g(r) dr \Big\} \Big|\bar\nu_\alpha(y)dy
 \ar\leq\ar  C \int_0^\infty \Big| \int_{(x-y)^+}^{x}v_\lambda^g(r) dr \Big| \bar\nu_\alpha(y)dy .
 \eeqnn
 Similarly as in the proof for Proposition~\ref{UpBoundVf01}, we can prove that the last quantity is finite.
 From this and the exponential formula for PRMs; see \cite[p.8]{Bertoin1996},
 \beqnn
 \mathbf{E}\big[ \exp\{ - Y_x(0) \} \big] =  \exp \Big\{- \zeta \int_0^\infty \Big(1-\exp\Big\{- \int_{(x-y)^+}^x   v_\lambda^g(s) ds\Big\}\Big)\bar\nu_\alpha(y)dy \Big\}<\infty.
 \eeqnn
 Notice that $e^{ -Z_x }= e^{ -Y_x(0)}\mathcal{E}_{U_x} $ on $[0,x]$.
 The standard conditional expectation argument yields that the local-martingale $e^{ -Z_x} $ is a true $(\mathscr{G}_r)$-martingale under $\mathbf{P}$.
 The equality (\ref{AffineRepLap}) can be obtained immediately from the facts that
 \beqnn
 Z_x(x)= \lambda\cdot  L^\xi_\zeta(x) + g*L^\xi_\zeta(x)\geq 0
 \quad\mbox{and}\quad
 \mathbf{E}\big[e^{-Z_x(x)}\big]= \mathbf{E}\big[e^{-Y_x(0)}\big].
 \eeqnn 
 \qed

 \subsection{Proof for Theorem~\ref{MainThm.03}}
 By Lemma~\ref{Lemma.NonlinearVE} and \ref{Lemma.LapFun}, associated to the unique $\mathcal{A}$-noncontinuable solution $(v_\lambda^g,T_\lambda^g)$ of (\ref{MainThm.Volterra}) we see that Theorem~\ref{MainThm.03} holds for any $x\in [0,T_\lambda^g)$.
 Thus it suffices to prove that $(v_\lambda^g,T_\lambda^g)$  is a $\mathcal{A}$-global solution, i.e., $T_\lambda^g=\infty$.
 In the sequel, we assume for contradiction that $T_\lambda^g<\infty$, which implies that 
 $$\limsup_{t\to T_\lambda^g-}\big|v_\lambda^g(t)\big|= \infty.$$ 
 With the help of the following propositions, we prove that $|v_\lambda^g(t)|\leq C \cdot t^{\alpha-1}$ uniformly in $(0,T_\lambda^g)$; see Lemma~\ref{Lemma.Upperv}, which leads to a contradiction to the preceding assumption,

 \begin{proposition}\label{Prop.01}
 	If $T_\lambda^g<\infty$, there exists a constant $C>0$ such that for any $t\in[0,T_\lambda^g]$ and $y>0$,
 	\beqlb\label{eqn.Prop.01.01}
 	\Big| \int_{(t-y)^+}^t v_\lambda^g(r) dr \Big|\leq C
 	\quad\mbox{and}\quad
 	\big|\mathcal{V}_\alpha\circ  v_\lambda^g(t)\big| \leq C  \int_0^\infty  \Big| \int_{(t-y)^+}^t v_\lambda^g(r) dr \Big|^2 \nu_\alpha(dy).
 	\eeqlb
 	
 \end{proposition}
 \proof  Notice that $e^{-x}-1+x\geq 0$ for any $x\in\mathbb{R}$. By (\ref{MainThm.Volterra}), we have 
 \beqnn
 v_\lambda^g(t) \leq \lambda W'(t)+g*W'(t)\leq C t^{\alpha-1}
 \quad\mbox{and hence}\quad
 \int_{(t-y)^+}^t v_\lambda^g(r) dr\leq C,
 \eeqnn
 uniformly in $t\in [0, T_\lambda^g]$ and $y>0$.
 We now prove 
 $$ \int_{(t-y)^+}^t v_\lambda^g(r) dr \geq - C, $$
  uniformly in $t\in [0, T_\lambda^g]$ and $y>0$. 
 If not, the fact that $ v_\lambda^g \in \mathcal{A}_T$ for any $T\in (0, T_\lambda^g)$ yields that for any $y>0$,
 $$ \int_{(T_\lambda^g-y)^+}^{T_\lambda^g} v_\lambda^g(r) dr  = -\infty. $$
 Moreover, by the inequality $1-e^{-x} \leq x$ for any $x\in\mathbb{R}$, Lemma~\ref{Lemma.LapFun} and (\ref{AffineRepLap}), we have for any $t\in (0,T_\lambda^g)$,
 \beqnn
 0\leq  \int_0^\infty \Big(1-\exp\Big\{- \int_{(t-y)^+}^t   v_\lambda^g(s) ds\Big\}\Big)\bar\nu_\alpha(y)dy \leq \int_0^\infty \bar\nu_\alpha(y)dy\int_{(t-y)^+}^t   v_\lambda^g(s) ds .
 \eeqnn
 The continuity of $ v_\lambda^g$ on $(0, T_\lambda^g)$ implies that the last integral tends to $-\infty$ as $t$ increases to $T_\lambda^g$, which leads to a contradiction and hence the first desired inequality holds.
 The second one can be proved similarly as in the proof of Proposition~\ref{UpBoundVf}.
 \qed
 
 %
 %

 %
 %
 %

 Recall the constant $\theta\in(1,\frac{1}{1-\alpha})$ defined in Proposition~\ref{ContractiveMap}.
 We need the next two constants in the follows
 \beqlb\label{eqn.Con.eta.l}
 \eta \in \Big(\Big(\frac{1}{\theta}-\frac{1+\alpha^2}{2}\Big)^+, \frac{1}{2\theta}\Big)
 \quad\mbox{and}\quad
 \ell \in \Big(\frac{1+\alpha}{2}, \frac{\eta+1-1/\theta}{1-\alpha}\wedge 1\Big).
 \eeqlb
 The first inequality in Proposition~\ref{UpBoundVf}, together with $0<\theta\eta<1/2$ and $\theta(\alpha-1)+1>0$, implies that the following function is well defined on $[0,T_\lambda^g)$,
 \beqnn
 H (t):=  \int_0^t   s^{-\theta\eta}\cdot \big|\mathcal{V}_\alpha \circ v_\lambda^g(t-s)\big|^{\theta} ds.
 \eeqnn
 
 \begin{proposition}
 	If $T_\lambda^g<\infty$, there exists a constant $C>0$ such that for any $t\in(0,T_\lambda^g)$,
 	\beqlb\label{eqn.501}
 	\big|\mathcal{V}_\alpha \circ v_\lambda^g(t)\big|
 	\ar\leq\ar C t^{\alpha-1} + C |H (t) |^{2/\theta}\cdot t^{\alpha+2\eta+1-2/\theta}\cr
 	\ar\ar + C t^{ 2(\alpha +\eta -\ell+1-1/\theta)}   \int_0^t  \frac{|H (t-s) |^{2/\theta}}{s^{\alpha+2-2\ell} }ds  .
 	\eeqlb
 \end{proposition}
 \proof For $y>0$, integrating both sides of (\ref{MainThm.Volterra}) over  $((t-y)^+, t]$ and then using Fubini's theorem,
 \beqnn
 \int_{(t-y)^+}^t v_\lambda^g(s) ds\ar=\ar \lambda \cdot \nabla_y W(t)  +  \big(g - \mathcal{V}_\alpha \circ v_\lambda^g\big) * \nabla_y W(t).
 \eeqnn
 Plugging this into the second inequality in (\ref{eqn.Prop.01.01}) and then using the Cauchy-Schwarz inequality, we have  $|\mathcal{V}_\alpha\circ v_\lambda^g(t)| \leq C \big[J_1(t) +J_2(t)+J_3(t)\big]$ uniformly in $t\in(0,T_\lambda^g)$, 
 where
 \beqnn
 J_1(t) \ar:=\ar   \lambda^2 \int_0^\infty \big|\nabla_y W(t)\big|^2 \nu_\alpha(dy),\cr
 J_2(t)\ar:=\ar \int_0^\infty  |g*\nabla_y W(t) |^2\nu_\alpha(dy), \cr
 J_3(t)\ar:=\ar \int_0^\infty \big| (\mathcal{V}_\alpha\circ v_\lambda^g)* \nabla_y W(t) \big|^2\nu_\alpha(dy).
 \eeqnn
 Similarly as in (\ref{eqn.101})-(\ref{eqn.102}),  we have $J_1(t)\leq C \cdot t^{\alpha-1}$ uniformly in $t\in(0, T_\lambda^g)$.
 By H\"older's inequality, we have 
 $$\big|g*\nabla_y W(t) \big| \leq  \|g\|_{L^\infty}\cdot \big\|\nabla_yW\big\|_{L^1_t} \leq \|g\|_{L^\infty} \cdot t^{1/2}\cdot \big\|\nabla_yW\big\|_{L^2_t} .$$
 Plugging this into $J_2(t)$ and then using Proposition~\ref{WeakConvergenceProp01} with $p=2$ as well as Fubini's theorem, we have uniformly in $t\in[0, T_\lambda^g)$,
 \beqnn
 J_2(t) \ar\leq\ar C   \cdot t\cdot\int_0^\infty \big\|\nabla_yW\big\|_{L^2_t}^2 \nu_\alpha(dy)\cr
 \ar=\ar  C \cdot t\cdot\int_0^tdr \int_0^\infty\big|\nabla_yW(r)\big|^2 \nu_\alpha(dy) \leq C\cdot t^{\alpha+1}.
 \eeqnn
 We now turn to analyze $J_3(t)$.
 Splitting  the interval of integration and then using the Cauchy-Schwarz inequality, we have $J_3(t) \leq J_{31}(t)+J_{32}(t)+J_{33}(t) $ with 
 \beqnn
 J_{31}(t) \ar:=\ar  \big| |\mathcal{V}_\alpha\circ v_\lambda^g| *W(t)  \big|^2 \cdot  \bar\nu_\alpha(t) ,\cr
 \ar\ar\cr
 J_{32}(t)\ar:=\ar 2\int_0^t\Big|\int_{t-y}^t \mathcal{V}_\alpha \circ v_\lambda^g(s)  W(t-s) ds \Big|^2\nu_\alpha(dy),\cr
 J_{33}(t) \ar:=\ar 2\int_0^{t}\Big|\int_0^{t-y}\mathcal{V}_\alpha\circ  v_\lambda^g(s) \nabla_y W(t-s) ds \Big|^2\nu_\alpha(dy).
 \eeqnn
 By H\"older's inequality and (\ref{UpperboundW}),
 \beqnn
 \big| |\mathcal{V}_\alpha\circ v_\lambda^g| *W(t)  \big|
 \ar=\ar \int_0^t s^{-\eta}\big|\mathcal{V}_\alpha\circ  v_\lambda^g(t-s)\big|\cdot s^\eta W(s)ds\cr
 \ar\leq\ar \big|H (t)\big|^{1/\theta}\cdot  \Big(\int_0^t \big|s^\eta W(s)\big|^{\theta/(\theta-1)}ds \Big)^{1-1/\theta} \cr
 \ar\leq\ar  C \cdot \big|H (t)\big|^{1/\theta}\cdot t^{\alpha+\eta +1-1/\theta}
 \eeqnn
 and hence 
 $$J_{31}(t) \leq C\cdot |H(t) |^{2/\theta}\cdot  t^{\alpha+2\eta+1-2/\theta}, $$ 
 uniformly in $t\in (0, T_\lambda^g)$.
 Similarly, we also have
 \beqnn
 \int_{t-y}^t \big|\mathcal{V}_\alpha\circ v_\lambda^g(s)\big|\cdot W(t-s)ds 
 \ar=\ar  \int_0^y  s^{-\eta}\big|\mathcal{V}_\alpha\circ v_\lambda^g(t-s)\big| \cdot s^{\eta} W(s)ds \cr
 \ar\leq\ar \Big( \int_0^y  s^{-\eta\theta} \cdot \big| \mathcal{V}_\alpha\circ v_\lambda^g(t-s) \big|^\theta ds\Big)^{1/\theta}\cr
 \ar\ar\quad \times \Big( \int_0^y \big|s^{\eta} W(s)\big|^{\theta/(\theta-1)} ds\Big)^{1-1/\theta},
 \eeqnn
 which can be bounded by $C \cdot |H(t) |^{1/\theta}\cdot y^{\alpha+\eta+1-1/\theta}$.
 Taking this back into $J_{32}(t)$, we have  
 $$J_{32}(t) \leq C |H(t) |^{2/\theta}\cdot t^{\alpha+2\eta+1-1/\theta},$$ 
 uniformly in $t\in[0, T_\lambda^g)$.
 For $J_{33}(t)$, we first consider its inner integral.
 Like the preceding argument, by H\"older's inequality we have
 \beqlb\label{eqn.103}
  \int_0^{t-y} \big|\mathcal{V}_\alpha\circ v_\lambda^g(s)\nabla_x W(t-s)\big| ds  
 \ar=\ar \int_y^{t} (s-y)^{-\eta}\cdot\big|\mathcal{V}_\alpha\circ v_\lambda^g(t-s)\big|\cdot (s-y)^{\eta} \nabla_y W(s)ds \cr
 \ar\leq\ar \Big(\int_y^{t} (s-y)^{-\eta\theta} \cdot \big|\mathcal{V}_\alpha \circ v_\lambda^g(t-s) \big|^\theta ds \Big)^{1/\theta}\cr
 \ar\ar\quad \times \Big( \int_y^{t} \big|(s-y)^{\eta}\nabla_y W(s)\big|^{\frac{\theta}{\theta-1}}ds \Big)^{1-1/\theta} .
 \eeqlb
 By the change of variables, the first term on the right side of this inequality equals to $ |H(t-y) |^{1/\theta}$.
 Recall the constant $\ell $ in (\ref{eqn.Con.eta.l}).
 By (\ref{UpperboundDiffW}), we have
 $$\big|\nabla_y W(s)\big| \leq C \cdot s^{\alpha(1-\ell)} \cdot (s-y)^{\ell(\alpha-1)} \cdot y^\ell, $$ 
 uniformly in $s>y>0$. 
 Plugging this into the second term on the right side of the inequality in (\ref{eqn.103}), it can be bounded by
 \beqnn
 Ct^{\alpha(1-\ell) }\cdot\Big(\int_y^{t}   (s-y)^{(\eta+(\alpha-1)\ell)\frac{\theta}{\theta-1}}ds\Big)^{\frac{\theta-1}{\theta}} \cdot y^{  \ell }
 \leq C \cdot t^{ \alpha +\eta -\ell+1-1/\theta} \cdot y^{\ell} .
 \eeqnn
 uniformly in $t\in[0, T_\lambda^g)$ and $y\in(0,t)$.
 Taking these two estimates back into (\ref{eqn.103}) and then $J_{33}(t)$, we have
 \beqnn
 J_{33}(t)
 \ar\leq\ar C \cdot t^{ 2(\alpha +\eta -\ell+1-1/\theta)}  \cdot\int_0^t \big|H (t-y)\big|^{2/\theta}\cdot y^{2 \ell -\alpha-2}dy.
 \eeqnn
 Here the constant $C>0$ is independent of $t$.
 Then (\ref{eqn.501}) follows by putting all estimates above together.
 \qed
 
 \begin{proposition}
 	If $T_\lambda^g<\infty$, there exists a constant $C_*>0$ such that for any $t\in(0,T_\lambda^g)$,
 	\beqlb \label{UpperH}
 	H (t)
 	\ar\leq\ar   C_*  t^{ -\eta\theta } +  C_* \int_0^t (t-s)^{-\eta\theta}\cdot|H (s)|^2dt.
 	\eeqlb
 \end{proposition}
 \proof  Raising both sides of the inequality (\ref{eqn.501}) to the $\theta$ power and then using the  power mean inequality, we have for some constant $C>0$ independent of $t$,
 \beqnn
 |\mathcal{V}_\alpha \circ v_\lambda^g(t)|^{\theta}
 \ar\leq\ar Ct^{\theta(\alpha-1)} +  C |H (t) |^{2 }\cdot t^{\theta(\alpha+2\eta+1)-2 }\cr
 \ar\ar\cr
 \ar\ar +C t^{ 2\theta(\alpha +\eta -\ell+1)-2} \Big(  \int_0^t  \frac{|H (t-s) |^{2/\theta}}{s^{\alpha+2-2\ell} }ds\Big)^\theta.
 \eeqnn
 Convolving both sides of this inequality by the power function $s^{-\theta\eta}$,
 \beqlb\label{eqn.502}
 H (t)
 \ar\leq\ar  C\int_0^t (t-s)^{-\eta\theta}\cdot  s^{\theta(\alpha-1)} ds  +  C \int_0^t (t-s)^{-\eta\theta}\cdot|H (s)|^2\cdot s^{\theta(\alpha+2\eta+1)-2}ds\cr
 \ar \ar +C\int_0^t (t-s)^{-\eta\theta}\cdot s^{2\theta(\alpha +\eta -\ell +1)-2 }  \Big( \int_0^s \frac{|H (s-r) |^{2/\theta}}{r^{\alpha+2-2\ell} }dr \Big)^{\theta} ds .
 \eeqlb
 Notice that $\eta\theta<1/2$, $\theta(\alpha-1)+1>0$ and $\theta(\alpha+2\eta+1)-2 >0$.
 A simple calculation shows that  uniformly in $t\in (0,T_\lambda^g)$, the first term on the right side of the above inequality can be bounded by $C \cdot t^{ -\eta\theta }$ and  the second term can be bounded by
 \beqlb\label{eqn.106}
 C\cdot t^{\theta(\alpha+2\eta+1)-2}\int_0^t (t-s)^{-\eta\theta}\cdot|H (s)|^2 ds.
 \eeqlb
 Using H\"older's inequality and the fact that $2 \ell -\alpha-1>0$, we have
 \beqnn
 \Big( \int_0^s \frac{|H (s-r) |^{2/\theta}}{r^{\alpha+2-2\ell} }dr \Big)^{\theta}
 \ar=\ar \Big( \int_0^s \frac{|H (s-r) |^{2/\theta}}{r^{(\alpha+2-2\ell)/\theta} }\cdot r^{(1-1/\theta))(2 \ell -\alpha-2)} dr \Big)^{\theta}  \cr
 \ar\leq\ar  \Big( \int_0^s z^{2 \ell -\alpha-2}dz \Big)^{\theta-1}   \int_0^s \frac{|H (s-r) |^{2 }}{r^{\alpha+2-2\ell} }  dr \cr
 \ar = \ar \frac{s^{(\theta-1)(2 \ell -\alpha-1)}}{|2 \ell -\alpha-1|^{\theta-1}} \int_0^s \frac{|H (s-r) |^{2 }}{r^{\alpha+2-2\ell} }  dr.
 \eeqnn
 Plugging this into the third term on the right side of (\ref{eqn.502}), it can be bounded uniformly in $t\in (0,T_\lambda^g)$ by
 \beqnn
 \lefteqn{   C  \int_0^t (t-s)^{-\eta\theta}\cdot s^{2\theta(\alpha +\eta -\ell +1)-2 +(\theta-1)(2 \ell -\alpha-1)}\cdot \int_0^s \frac{|H (s-r) |^{2 }}{r^{\alpha+2-2\ell} }  dr  ds }\ar\ar \cr
 \ar\leq\ar  C t^{\theta(\alpha +2\eta  +1)-2 -(2 \ell -\alpha-1)}\cdot \int_0^t (t-s)^{-\eta\theta}\cdot \int_0^s \frac{|H (s-r) |^{2 }}{r^{\alpha+2-2\ell} }  dr  ds \cr
 \ar=\ar C \cdot t^{\theta(\alpha +2\eta  +1)-2 -(2 \ell -\alpha-1)} \cdot \int_0^t |H (s)|^{2 }\int_0^{t-s} \frac{(t-s-r)^{-\eta\theta}}{r^{\alpha+2 -2 \ell}}dr ds.
 \eeqnn
 Here the inequality comes from the fact that $2\theta(\alpha +\eta -\ell +1)-2>0$ as well as $(\theta-1)(2 \ell -\alpha-1)>0$, and the equality comes from Fubini's theorem.
 Noting that $2 \ell -\alpha-1> 0$ and $\eta\theta<1/2$,
 a simple calculation induces that uniformly in $t\geq s>0$,
 \beqnn
 \int_0^{t-s} \frac{(t-s-r)^{-\eta\theta}}{r^{\alpha+2 -2 \ell}}dr
 \leq C\cdot (t-s)^{-\eta\theta+2 \ell -\alpha-1 } \leq  C \cdot t^{2 \ell -\alpha-1} \cdot (t-s)^{-\eta\theta} . 
 \eeqnn
 Consequently, the third term on the right side of  (\ref{eqn.502}) can be bounded uniformly in $t\in(0,T_\lambda^g)$ by
 \beqnn
 C\cdot t^{\theta(\alpha+2\eta+1)-2 } \cdot \int_0^t (t-s)^{-\eta\theta  }  |H (s)|^{2}ds,
 \eeqnn
 which can be merged with (\ref{eqn.106}). Then (\ref{UpperH}) follows by  putting these estimates together.
 \qed
 
 \begin{proposition}\label{Prop.FracRicc}
 	For any $C^*>0$, there exists a unique continuous and non-negative solution of
 	\beqlb\label{FracRiccati}
 	\psi(t)=C^* t^{ -\eta\theta } + C^*\int_0^t (t-s)^{-\eta\theta}\cdot|\psi(s)|^2ds,\quad t>0.
 	\eeqlb
 	Moreover, for any $T>0$, there exists a constant $C>0$ such that for any $t\in(0,T]$,
 	\beqnn
 	\psi(t)\leq C \cdot t^{-\eta\theta}.
 	\eeqnn
 \end{proposition}
 \proof By Theorem~6.1(ii) in \cite{JaberLarssonPulido2019} and $0<\eta\theta<1/2$, there exists a unique solution $\tilde\psi\in L^2_{\rm loc}(\mathbb{R}_+;\mathbb{R}_+)$ of (\ref{FracRiccati}).
 Let 
 $$
 \psi(t):=C^* t^{ -\eta\theta } + C^* \int_0^t (t-s)^{-\eta\theta}\cdot|\tilde\psi(s)|^2ds,
 \quad t>0. $$
 By the properties of convolution, it is easy to identify that $ \psi$ is continuous on $(0,\infty)$ and equal to $\tilde\psi$ almost everywhere. Thus $ \psi$ is the unique continuous and non-negative solution of (\ref{FracRiccati}).
 By Theorem~2.a in \cite{CallegaroGrasselliPages2021}, there exist two constants $C_\psi,r_\psi>0$ such that 
 $$\psi(t)\leq C_\psi t^{-\eta\theta},$$
 uniformly on $t\in(0, r_\psi]$. 
 For $T>r_\psi$, the continuity of $\psi$ yields that  $\psi(t)\leq C  \cdot t^{-\eta\theta} $ for any $t\in(0,T]$ and some $C>0$.
 \qed
 
 \begin{proposition}\label{Prop.05}
 	If $T_\lambda^g<\infty$, there exists a constant $C>0$ such that $H(t)\leq C  \cdot t^{-\eta\theta} $ for any $t\in(0,T_\lambda^g)$.
 \end{proposition}
 \proof Choosing the two constants $C^*> C_*>0$ such that the inequality (\ref{UpperH}) turns to be strict.
 It is easy to identify that both $t^{\eta\theta}H(t)$ and $t^{\eta\theta}\psi(t)$ are continuous on $[0,T_\lambda^g)$.
 By Theorem~2.1 in \cite{DentonVatsala2010}\,\footnote{For $T>0$ and $\rho\in(0,1)$, let $f_1,f_2 $ be two functions on $(0,T]$ satisfying  that $t^{\rho}f_i(t) \in C([0,T];\mathbb{R})$ with $i=1,2$. If
 	$
 	f_1(t)< C_1 t^{-\rho} + C_1 \int_0^t (t-s)^{-\rho} f_1(s)ds$
 	and
 	$f_2(t)= C_2 t^{-\rho} + C_2 \int_0^t (t-s)^{-\rho} f_2(s)ds$ with $C_1<C_2$,
 	then $f_1<f_2$ on $(0,T]$.}, the function $H$ can be uniformly bounded by $\psi$ on $(0,T_\lambda^g)$ and the desired result follows directly from Proposition~\ref{Prop.FracRicc}.
 \qed
 
 \begin{lemma}\label{Lemma.Upperv}
 	If $T_\lambda^g<\infty$, there exists a constant $C>0$ such that $|v_\lambda^g(t)|\leq C \cdot t^{\alpha-1} $ for any $t\in(0,T_\lambda^g)$.
 \end{lemma}
 \proof Plugging Proposition~\ref{Prop.05} into (\ref{eqn.501}), we have
 \beqnn
 |\mathcal{V}_\alpha \circ v_\lambda^g(t)|
 \ar\leq\ar C t^{\alpha-1} + C  t^{\alpha +1-2/\theta}
 + C t^{ 2(\alpha +\eta -\ell+1-1/\theta)}   \int_0^t  \frac{ (t-s)^{-2\eta}}{s^{\alpha+2-2\ell} }ds  .
 \eeqnn
 Notice that $\alpha +1-2/\theta \in (\alpha-1,1-\alpha)$,  $\eta <1/2$ and $\alpha+2-2\ell<1$, we have  $|\mathcal{V}_\alpha \circ v_\lambda^g(t)| \leq C\cdot t^{\alpha-1}$ uniformly in  $t\in[0,T_\lambda^g)$.
 Taking this back into (\ref{MainThm.Volterra}), we can get the desired result immediately.
 \qed

 \subsection{Proof for Theorem~\ref{MainThm.01}(3)} Assume that $L^\xi_{\zeta,1}$ and $L^\xi_{\zeta,2}$ are two weak solutions of (\ref{MainThm.SVE}). For any $x,z\geq 0$ and $g\in L^\infty(\mathbb{R}_+;\mathbb{R}_+)$, let $v_{0}^{g} $ be the unique $\mathcal{A}$-global solution of (\ref{MainThm.Volterra}) with $\lambda =0$.
 By Theorem~\ref{MainThm.03},
 \beqnn
 \mathbf{E}\big[e^{-z \cdot g* L^\xi_{\zeta,1}(x)  }\big]
 \ar=\ar \exp \Big\{- \zeta \int_0^\infty \Big(1-\exp\Big\{- \int_{(x-y)^+}^x   v_0^g(s) ds\Big\}\Big)\bar\nu_\alpha(y)dy \Big\}\cr
 \ar=\ar \mathbf{E}\Big[\exp\Big\{-z \cdot g* L^\xi_{\zeta,2}(x)  \Big\}\Big].
 \eeqnn
 The one-to-one correspondence between non-negative random variables and their Laplace transforms yields that the two non-negative random variables
 $g* L^\xi_{\zeta,1}(x)$ and $g* L^\xi_{\zeta,2}(x)$
 are equal in distribution.
 Hence the two solutions $L^\xi_{\zeta,1}$ and $L^\xi_{\zeta,2}$ have the same probability law on $L^1([0,x];\mathbb{R}_+)$ and also on $C([0,x];\mathbb{R}_+)$. By the arbitrariness of $x$, the weak uniqueness of non-negative solutions holds for (\ref{MainThm.SVE}).
 
  \section{Fractional integral representations}\label{FIR}

 In this section we prove the two equivalences in Theorem~\ref{MainThm.FIR}. When $b=0$, they follow directly from (\ref{ScaleF0}).
 We now prove them with $b>0$.
 The equivalence between (\ref{MainThm.Volterra}) and (\ref{FractionalRiccatti}) follows from the resolvent equation (\ref{KenelResolvent01}).
 Indeed, by Theorem~4.6\,\footnote{For two function $f,\mathrm{k}\in L^1_{\rm loc}(\mathbb{R}_+;\mathbb{R})$, we have $x= f+\mathrm{k}*x$ if and only if $x=f+R_\mathrm{k}*f$, where $R_\mathrm{k}$ is the unique solution of $R_\mathrm{k}=\mathrm{k}+\mathrm{k}*R_\mathrm{k}$.} in \cite[p.48]{GripenbergLondenStaffans1990} and (\ref{KenelResolvent01}) we have  $v_\lambda^g$ solves (\ref{FractionalRiccatti}) if and only if
 \beqnn
 v_\lambda^g\ar=\ar \lambda K+   (g- \mathcal{V}_\alpha\circ v_\lambda^g )*K
 -bW'*\big(\lambda K+  ( g- \mathcal{V}_\alpha\circ v_\lambda^g )*K \big) \cr
 \ar\ar\cr
 \ar=\ar \lambda (K-bW'* K ) + (g-\mathcal{V}_\alpha\circ v_\lambda^g )*(K-bW'* K ).
 \eeqnn
 Multiplying both sides by $b$ and then using (\ref{KenelResolvent01}) again, we have
 \beqnn
 b v_\lambda^g
 \ar=\ar \lambda (bK-bW'* bK ) + (g-\mathcal{V}_\alpha\circ v_\lambda^g )*(bK-bW'*b K )\cr
 \ar\ar\cr
 \ar =\ar \lambda  bW'  + (g-\mathcal{V}_\alpha\circ v_\lambda^g )* bW',
 \eeqnn
 which is equivalent to (\ref{MainThm.Volterra}).
 
 The equivalence between (\ref{MainThm.SVE}) and (\ref{FractionalSVR}) can be proved in the same way.
 For convenience, we assume $K(x)=0$ if $x\leq 0$.
 Notice that 
 $$\int_{s-y}^{s} K(r)dr = \int_0^s\nabla_y K(r)dr$$ 
 for any $s,y\geq 0$.
 By Theorem~4.6 in \cite[p.48]{GripenbergLondenStaffans1990} and (\ref{KenelResolvent01}) again, the process $L^\xi_\zeta$ is a solution of (\ref{FractionalSVR}) if and only if it solves
 \beqlb\label{eqn.701}
 L^\xi_\zeta(x)\ar=\ar \zeta - \zeta\cdot bW(x)  +  \int_0^\infty  \int_0^\zeta  \int_{0}^{x}  \nabla_y K(r) dr \widetilde{N}_0( dy,dz)\cr
 \ar\ar -\int_0^x bW'(x-t) dt \int_0^\infty  \int_0^\zeta  \int_{0}^{t}  \nabla_yK(r) dr \widetilde{N}_0( dy,dz)\cr
 \ar\ar + \int_0^x \int_0^\infty \int_0^{L^\xi_\zeta(s)} \int_0^{x-s} \nabla_y K(r)dr \widetilde{N}_\alpha(ds,dy,dz)\cr
 \ar\ar -\int_0^x bW'(x-t) dt\int_0^t \int_0^\infty \int_0^{L^\xi_\zeta(s)} \int_0^{t-s} \nabla_y K(r)dr \widetilde{N}_\alpha(ds,dy,dz).
 \eeqlb
 By the change of variables and Proposition~\ref{Prop.D4},
 \beqlb\label{eqn.108}
 \lefteqn{\int_0^x bW'(x-t)\Big| \int_0^t ds \int_0^\epsilon \Big|  \int_0^{t-s} \nabla_y K(r)dr  \Big|^2 \nu_\alpha(dy) \Big|^{1/2} dt }\qquad \ar\ar\cr
 \ar\leq\ar C\cdot W(x)\cdot \Big| \int_0^x ds \int_0^\epsilon \frac{ |  \int_0^{s} \nabla_y K(r)dr |^2}{y^{\alpha+2}}dy \Big|^{1/2},
 \eeqlb
 which goes to $0$ as $\epsilon\to 0+$.
 By the change of variables and Fubini's theorem,
 \beqnn
 \lefteqn{\int_0^x bW'(x-t) \int_{t-s-y}^{t-s} K(r)dr\cdot \mathbf{1}_{\{ 0\leq s<t\}}dt}\ar\ar\cr
 \ar=\ar \int_0^{x-s} bW'(x-s-t) \int_0^{s} \nabla_y K(r)dr  dt\cr
 \ar=\ar \int_0^{x-s} \nabla_yK(t)\int_t^{x-s} bW'(x-s-r) dr dt \cr
 \ar=\ar \int_0^{x-s} \nabla_yK(t) bW(x-s-t) dr dt \cr
 \ar\ar\cr
 \ar=\ar  b  W*\nabla_yK(x-s).
 \eeqnn
 By the change of variables and Proposition~\ref{Prop.B6}, we have for $\epsilon>0$,
 \beqnn
 \int_0^x ds\int_0^\epsilon |b  W*\nabla_yK(x-s)|^2\nu_\alpha(dy) \leq C \int_0^x ds\int_0^\epsilon \frac{|  W*\nabla_yK(s)|^2}{y^{\alpha+2}}dy ,
 \eeqnn
 which is finite and goes to $0$ as $\epsilon\to 0+$.
 From this, (\ref{eqn.108}) and the stochastic Fubini theorem; see Theorem~\ref{StoFubiniThm}, the stochastic integral
 \beqnn
 \int_0^x \int_0^\infty \int_0^{L^\xi_\zeta(s)} b  W*\nabla_yK(x-s) \widetilde{N}_\alpha(ds,dy,dz)
 \eeqnn
 is well defined and equal almost surely to the last stochastic integral on  the right side of (\ref{eqn.701}).
 Moreover, by Fubini's theorem and (\ref{KenelResolvent02}),
 \beqnn
 \int_0^{x-s} \nabla_y K(r)dr -  b  W*\nabla_yK(x-s)
 \ar=\ar  \big( 1- bW \big) * \nabla_y K(x-s)\cr
 \ar=\ar  W'*L_K*  \nabla_y K(x-s)=   \nabla_y W(x-s).
 \eeqnn
 Thus the subtraction of  the last two terms on the right side of (\ref{eqn.701}) is equal almost surely to
 \beqnn
 \int_0^x \int_0^\infty \int_0^{L^\xi_\zeta(s)}   \nabla_y W(x-s) \widetilde{N}_\alpha(ds,dy,dz).
 \eeqnn
 Similarly, the following stochastic integral
 \beqnn
 \int_0^\infty  \int_0^\zeta   \int_0^x bW'(x-t)\int_{0}^{t}\nabla_yK(r) dr dt  \widetilde{N}_0( dy,dz)
 \eeqnn
 is also well defined and equal almost surely to the second stochastic integral on the right side of (\ref{eqn.701}). Moreover,
 \beqnn
 \int_{0}^{x}  \nabla_y K(r) dr-  \int_0^x bW'(x-t)\int_{0}^{t}\nabla_yK(r) dr dt = \int_0^x \nabla_y W'(r)dy =  \nabla_y W(x).
 \eeqnn
 The subtraction of the first two integrals on the right side of (\ref{eqn.701}) is equal almost surely to
 \beqnn
 \int_0^\infty  \int_0^\zeta  \nabla_y W(x) \widetilde{N}_0( dy,dz).
 \eeqnn
 Putting these results  together,  we see that  (\ref{eqn.701}) turns into (\ref{Equ.SVE}) and hence (\ref{FractionalSVR}) is equivalent to (\ref{MainThm.SVE}).
 \qed
 
  \section{Application to M/G/1 processor-sharing queues} \label{Sec.PS}

 As the continuation of \cite{LambertSimatos2015}, we use the preceding results to establish a SVE for the heavy-traffic limit of heavy-tailed M/G/1 processor sharing queues.
 Recall the sequence $\{\gamma_n\}_{n\geq 1}$ and  the Pareto distribution $\Lambda$ defined in Section~\ref{Sec.CP}.
 In the $n$-th  processor-sharing queue, the arrival of customers to the
 system is described by a Poisson process with rate $\gamma_n>0$ and the amount of processing time that each  customer requires from the server is distributed as $\Lambda$.
 Additionally, there are $z_n$ initial customers in the system at time $0$, whose residual service times are independent and identically distributed with common distribution $\Lambda^*$.
 Here we are interest in the heavy-traffic limit of the queue-length process before the queue becoming empty.
 More precisely, let $q^{(n)}:=\{ q^{(n)}(t):t\geq 0 \}$ be the queue-length process and $\tau^{(n)}$ the first time that the queue becomes empty, i.e. $\tau^{(n)}:= \inf\{t>0:q^{(n)}(t)=0\}$.
 We write  $Q^{(n)}$ for the rescaled queue-length process $ \{n^{-\alpha/(1+\alpha)}\cdot q^{(n)}(nt) :t\in [0,\tau^{(n)}] \}$ under $\mathbf{P}(\,\cdot\,|\tau^{(n)}<\infty)$.
 
 Denote by $\mathcal{E}$ the set of all positive excursions with finite length.
 For each $f\in \mathcal{E}$, let $e_f$ be the right end point of $f$ and
 \beqlb\label{eqn.If}
 \mathcal{I}_f(t):=  \int_0^{t } f(s)ds,\quad t\geq 0.
 \eeqlb
 It is obvious that $\mathcal{I}_f$ is a continuous and non-decreasing function on $\mathbb{R}_+$, which allows us to define its right-inverse function  $ \mathcal{I}_f^{-1}$ by  $\mathcal{I}_f^{-1}(t)= e_f$ if $t> \mathcal{I}_f(\infty)$ and
 $ \mathcal{I}_f^{-1} (t):= \inf\{ s\geq 0 : \mathcal{I}_f(s)\geq t \} $ if $t\in[0, \mathcal{I}_f(\infty)]$.
 Let $\mathscr{L}$ be the Lamperti transformation on $\mathcal{E}$, which is a map acting on an excursion $f\in \mathcal{E}$ by $\mathscr{L}\circ f(t):= f( \mathcal{I}_f^{-1}(t))$ for $t\geq 0$.
 The next corollary is a direct consequence of Theorem~6.5 in \cite{LambertSimatos2015} and Theorem~\ref{MainThm.FIR}.
 
 \begin{corollary}
 	If Condition~\ref{Main.Condition} holds and $z_n/n^{\alpha/(1+\alpha)}\to \zeta>0$ as $n\to\infty$, we have  $Q^{(n)} \to Q_\zeta$ weakly in $D([0,\infty),\mathbb{R}_+)$, where the limit process $ Q_\zeta \in \mathcal{E}$ is the unique weak solution of
 	\beqlb\label{SVE.PS}
 	Q_\zeta(t) \ar =\ar \zeta - b \int_0^t \frac{\big(\mathcal{I}_{Q_\zeta}^{-1}(t)-\mathcal{I}_{Q_\zeta}^{-1}(s)\big)^{\alpha-1}}{\Gamma(\alpha)\Gamma(1-\alpha)}ds \cr
 	\ar\ar + \int_0^\infty  \int_0^\zeta \Big(\int_{(\mathcal{I}_{Q_\zeta}^{-1}(t)-y)^+}^{\mathcal{I}_{Q_\zeta}^{-1}(t)}\frac{r^{\alpha-1}dr }{\Gamma(\alpha)\Gamma(1-\alpha)} \Big) \widetilde{N}_{Q,0}( dy,dz) \cr
 	\ar\ar + \int_0^t \int_0^\infty \Big(\int_{(\mathcal{I}_{Q_\zeta}^{-1}(t)-\mathcal{I}_{Q_\zeta}^{-1}(s)-y)^+}^{\mathcal{I}_{Q_\zeta}^{-1}(t)-\mathcal{I}_{Q_\zeta}^{-1}(s)}\frac{r^{\alpha-1}dr }{\Gamma(\alpha)\Gamma(1-\alpha)} \Big)\widetilde{N}_Q(ds,dy),
 	\quad t\geq0,
 	\eeqlb
 	where $\nu_\alpha(dy)$ is given by (\ref{LevyMeasure}) with $c=\Gamma(1-\alpha)$, $\widetilde{N}_{Q,0}( dy,dz)$ and $\widetilde{N}_Q(ds,dy)$ are two compensated PRMs on $(0,\infty)^2$ with intensity  $\bar\nu_\alpha(y)dydz$ and $ds\nu_\alpha(dy)$, respectively.
 \end{corollary}
 \proof Let $L^\xi_\zeta$ be the unique weak solution of (\ref{FractionalSVR}) with $c=\Gamma(1-\alpha)$.
 From Theorem~6.5 in \cite{LambertSimatos2015}, we have $Q^{(n)} \to \mathscr{L}\circ L^\xi_\zeta \in\mathcal{E}$ weakly in $D([0,\infty),\mathbb{R}_+)$. By (\ref{FractionalSVR})  and the change of variables,
 \beqlb\label{eqn.100}
 \lefteqn{\mathscr{L}\circ L^\xi_\zeta(t)}\ar\ar\cr
 \ar=\ar  \zeta- b  \int_0^{\mathcal{I}_{L^\xi_\zeta}^{-1}(t)} \frac{(\mathcal{I}_{L^\xi_\zeta}^{-1}(t)-s)^{\alpha-1}  }{\Gamma(\alpha)\Gamma(1-\alpha)} L^\xi_\zeta(s)ds  \cr
 \ar\ar +  \int_0^\infty  \int_0^\zeta \Big(\int_{(\mathcal{I}_{L^\xi_\zeta}^{-1}(t)-y)^+}^{\mathcal{I}_{L^\xi_\zeta}^{-1}(t)}\frac{r^{\alpha-1}dr }{\Gamma(\alpha)\Gamma(1-\alpha)} \Big) \widetilde{N}_0( dy,dz)\cr
 \ar\ar + \int_0^{\mathcal{I}_{L^\xi_\zeta}^{-1}(t)} \int_0^\infty \int_0^{L^\xi_\zeta(s)} \Big(\int_{(\mathcal{I}_{L^\xi_\zeta}^{-1}(t)-s-y)^+}^{\mathcal{I}_{L^\xi_\zeta}^{-1}(t)-s}\frac{r^{\alpha-1}dr }{\Gamma(\alpha)\Gamma(1-\alpha)} \Big)
 \widetilde{N}_\alpha(ds,dy,dz)\cr
 \ar=\ar   \zeta- b  \int_0^{t} \frac{(\mathcal{I}_{L^\xi_\zeta}^{-1}(t)-\mathcal{I}_{L^\xi_\zeta}^{-1}(s))^{\alpha-1}  }{\Gamma(\alpha)\Gamma(1-\alpha)} L^\xi_\zeta(\mathcal{I}_{L^\xi_\zeta}^{-1}(s))d\mathcal{I}_{L^\xi_\zeta}^{-1}(s) \cr
 \ar\ar +  \int_0^\infty  \int_0^\zeta \Big(\int_{(\mathcal{I}_{L^\xi_\zeta}^{-1}(t)-y)^+}^{\mathcal{I}_{L^\xi_\zeta}^{-1}(t)}\frac{r^{\alpha-1}dr }{\Gamma(\alpha)\Gamma(1-\alpha)} \Big) \widetilde{N}_0( dy,dz)\cr
 \ar\ar + \int_0^{t} \int_0^\infty \int_0^{L^\xi_\zeta(\mathcal{I}_{L^\xi_\zeta}^{-1}(s))} \Big(\int_{(\mathcal{I}_{L^\xi_\zeta}^{-1}(t)-\mathcal{I}_{L^\xi_\zeta}^{-1}(s)-y)^+}^{\mathcal{I}_{L^\xi_\zeta}^{-1}(t)-\mathcal{I}_{L^\xi_\zeta}^{-1}(s)}\frac{r^{\alpha-1}dr }{\Gamma(\alpha)\Gamma(1-\alpha)} \Big)
 \widetilde{N}_\alpha(d\mathcal{I}_{L^\xi_\zeta}^{-1}(s),dy,dz)\cr
 \ar=\ar  \zeta- b  \int_0^{t} \frac{(\mathcal{I}_{L^\xi_\zeta}^{-1}(t)-\mathcal{I}_{L^\xi_\zeta}^{-1}(s))^{\alpha-1}  }{\Gamma(\alpha)\Gamma(1-\alpha)}   \mathscr{L}\circ L^\xi_\zeta(s)d\mathcal{I}_{L^\xi_\zeta}^{-1}(s) \cr
 \ar\ar +  \int_0^\infty  \int_0^\zeta \Big(\int_{(\mathcal{I}_{L^\xi_\zeta}^{-1}(t)-y)^+}^{\mathcal{I}_{L^\xi_\zeta}^{-1}(t)}\frac{r^{\alpha-1}dr }{\Gamma(\alpha)\Gamma(1-\alpha)} \Big) \widetilde{N}_0( dy,dz)\cr
 \ar\ar + \int_0^{t} \int_0^\infty   \Big(\int_{(\mathcal{I}_{L^\xi_\zeta}^{-1}(t)-\mathcal{I}_{L^\xi_\zeta}^{-1}(s)-y)^+}^{\mathcal{I}_{L^\xi_\zeta}^{-1}(t)-\mathcal{I}_{L^\xi_\zeta}^{-1}(s)}\frac{r^{\alpha-1}dr }{\Gamma(\alpha)\Gamma(1-\alpha)} \Big)
 \widetilde{N}_\alpha(d\mathcal{I}_{L^\xi_\zeta}^{-1}(s),dy,(0,\mathscr{L}\circ L^\xi_\zeta(s)]) .
 \eeqlb
 For any $s\in [0, \mathcal{I}_{L^\xi_\zeta}(\infty)]$, by  (\ref{eqn.If}) we have
 \beqnn
 s= \mathcal{I}_{L^\xi_\zeta} (\mathcal{I}_{L^\xi_\zeta}^{-1}(s)) 
 \ar=\ar \int_0^{\mathcal{I}_{L^\xi_\zeta}^{-1}(s)}  L^\xi_\zeta(r) dr\cr
 \ar =\ar \int_0^{s}  L^\xi_\zeta(\mathcal{I}_{L^\xi_\zeta}^{-1}(r)) d\mathcal{I}_{L^\xi_\zeta}^{-1}(r)
 = \int_0^{s} \mathscr{L}\circ L^\xi_\zeta(r)d\mathcal{I}_{L^\xi_\zeta}^{-1}(r).
 \eeqnn
 Differentiating both sides of this equality and then moving $ \mathscr{L}\circ L^\xi_\zeta$ to the left side of the first equality,
 \beqnn
 \frac{ds}{ \mathscr{L}\circ L^\xi_\zeta(s)} = d\mathcal{I}_{L^\xi_\zeta}^{-1}(s) .
 \eeqnn
 Integrating both side of this equality over $[0,t]$, we have
 \beqnn
 \mathcal{I}_{L^\xi_\zeta}^{-1}(t) = \int_0^t \frac{ds}{ \mathscr{L}\circ L^\xi_\zeta(s)} = \mathcal{I}_{\mathscr{L}\circ L^\xi_\zeta}^{-1}(t) .
 \eeqnn
 Plugging these back into the terms on  the right side of the last equality in (\ref{eqn.100}) yields that
 \beqnn
 \mathscr{L}\circ L^\xi_\zeta(t)\ar=\ar
 \zeta- b  \int_0^{t} \frac{\big(\mathcal{I}_{\mathscr{L}\circ L^\xi_\zeta}^{-1}(t) -\mathcal{I}_{\mathscr{L}\circ L^\xi_\zeta}^{-1}(s) \big)^{\alpha-1}  }{\Gamma(\alpha)\Gamma(1-\alpha)}  ds \cr
 \ar\ar +  \int_0^\infty  \int_0^\zeta \Big(\int_{(\mathcal{I}_{\mathscr{L}\circ L^\xi_\zeta}^{-1}(t)-y)^+}^{\mathcal{I}_{\mathscr{L}\circ L^\xi_\zeta}^{-1}(t)}\frac{r^{\alpha-1}dr }{\Gamma(\alpha)\Gamma(1-\alpha)} \Big) \widetilde{N}_{Q,0}(dy,dz)\cr
 \ar\ar + \int_0^{t} \int_0^\infty   \Big(\int_{(\mathcal{I}_{\mathscr{L}\circ L^\xi_\zeta}^{-1}(t) -\mathcal{I}_{\mathscr{L}\circ L^\xi_\zeta}^{-1}(s)-y)^+}^{\mathcal{I}_{\mathscr{L}\circ L^\xi_\zeta}^{-1}(t) -\mathcal{I}_{\mathscr{L}\circ L^\xi_\zeta}^{-1}(s)}\frac{r^{\alpha-1}dr }{\Gamma(\alpha)\Gamma(1-\alpha)} \Big)
 \widetilde{N}_Q(ds,dy),
 \eeqnn
 where $\widetilde{N}_{Q,0}(dy,dz):=\widetilde{N}_{0}(dy,dz)$ and $\widetilde{N}_Q(ds,dy):= \widetilde{N}_\alpha(d\mathcal{I}_{L^\xi_\zeta}^{-1}(s),dy,(0,\mathscr{L}\circ L^\xi_\zeta(s)])$. It is easy to identify that $\widetilde{N}_Q(ds,dy)$ is a compensated PRM on $(0,\infty)^2$ with intensity
 $ \mathscr{L}\circ L^\xi_\zeta(s) d\mathcal{I}_{L^\xi_\zeta}^{-1}(s) \nu_\alpha(dy) =  ds  \nu_\alpha(dy) $.
 Consequently, the limit process $ \mathscr{L}\circ L^\xi_\zeta$ is a weak solution of (\ref{SVE.PS}).
 The weak uniqueness of solutions of (\ref{SVE.PS}) follows directly from Theorem~\ref{MainThm.01}.
 \qed
 
 \begin{appendix}
   \section{Some technical results about the scale function}\label{TechResW}

 \begin{proposition}\label{WeakConvergenceProp01}
  For $p>1+\alpha$,  there exists a constant $C>0$ such that for any $h\geq 0$,
  \beqnn 
  \int_0^\infty \frac{|\Delta_h  W  (s) |^p }{ (s+h)^{\alpha+1}} ds  + \int_0^h ds\int_0^\infty \frac{|\nabla_y W (s)|^p}{ y^{\alpha+2} } dy  \leq C h^{(p-1)\alpha}.
  \eeqnn
 \end{proposition}
 \proof
 From (\ref{UpperboundDiffW}),  there exists a constant $C>0$ such that for any $h\geq 0$,
 \beqnn
  \int_0^\infty\frac{|\Delta_h  W  (s) |^p }{ (s+h)^{\alpha+1}} ds
  \ar\leq\ar C\int_0^h   (s+h)^{p\alpha-\alpha-1}ds \cr
  \ar\ar + C   \int_h^\infty   \frac{h^ps^{p(\alpha-1)}} {(s+h)^{\alpha+1}}ds \leq C \cdot  h^{(p-1)\alpha}.
 \eeqnn
 Similarly, there exists a constant $C>0$ such that for any $h\geq 0$,
 \beqnn
 \lefteqn{\int_0^h ds\int_0^s |\nabla_y W (s)|^p y^{-\alpha-2}dy}\ar\ar\cr
 \ar \leq\ar C \int_0^h ds \int_0^{s/2} \frac{(s-y)^{p(\alpha-1)}}{ y^{\alpha+2-p}} dy   +C  \int_0^h  ds \int_{s/2}^\infty\frac{s^{p\alpha}}{y^{\alpha+2}}dy  \cr
 \ar\leq\ar C\int_0^h s^{p(\alpha-1)}  ds \int_0^{s/2}y^{p-\alpha-2}dy +C\int_0^h s^{p\alpha-\alpha-1}ds
 \leq C\cdot h^{(p-1)\alpha}.
 \eeqnn
 \qed

 \begin{proposition}\label{WeakConvergenceProp02}
 For  $p\geq 2$, there exists a constant $C>0$ such that for any $h\geq 0$,
 \beqlb\label{WeakConvergence.eqn.02}
 \int_0^\infty ds \int_0^{s} \frac{|\nabla_y\Delta_h W (s) |^p }{ y^{\alpha+2} } dy  \leq C  h^{(p-1)\alpha} .
 \eeqlb
 \end{proposition}
 \proof
 We first split the double integral in (\ref{WeakConvergence.eqn.02}) into the following four parts:
 \beqnn
 J_1(h)\ar:=\ar \int_0^{4h} ds\int_0^{s}  \frac{|\nabla_y\Delta_h W (s) |^p }{ y^{\alpha+2} } dy  ,\cr
 J_2(h)\ar:=\ar \int_{4h}^\infty ds\int_{s-h}^{s}  \frac{|\nabla_y\Delta_h W (s) |^p }{ y^{\alpha+2} } dy, \cr
 J_3(h)\ar:=\ar  \int_{4h}^\infty ds\int_{s/2}^{s-h} \frac{|\nabla_y\Delta_h W (s) |^p }{ y^{\alpha+2} } dy, \cr
 J_4(h)\ar:=\ar \int_{4h}^\infty ds \int_0^{s/2}  \frac{|\nabla_y\Delta_h W (s) |^p }{ y^{\alpha+2} } dy.
 \eeqnn
 The power mean inequality, along with the equality $\nabla_y\Delta_h W(s)=  \nabla_yW(s+h)-\nabla_yW(s)$, implies that uniformly in $h\geq 0$,
 \beqnn
 J_1(h)\ar\leq\ar C\int_0^{4h}ds \int_0^{s}  \frac{|\nabla_y W (s+h)  |^p}{y^{\alpha+2}} dy
 + C\int_0^{4h}ds\int_0^{s} \frac{|\nabla_y W (s)|^p}{ y^{\alpha+2} } dy.
 \eeqnn
 By Proposition~\ref{WeakConvergenceProp01}, the second term on the right side of this inequality can be bounded by $C h^{(p-1)\alpha}$ uniformly in $h\geq 0$.
 Applying the change of variables to the first term and then using Proposition~\ref{WeakConvergenceProp01} again, it can be bounded uniformly in $h\geq 0$ by
 \beqnn
 C \int_h^{5h}ds \int_0^{s-h}  \frac{|\nabla_y W (s) |^p}{y^{\alpha+2}} dy
 \leq  C \int_0^{5h}ds \int_0^{s}  \frac{|\nabla_y W (s) |^p}{y^{\alpha+2}} dy \leq C h^{(p-1)\alpha}.
 \eeqnn
 Hence $J_1(h) \leq  C\cdot h^{(p-1)\alpha}$ uniformly in $h\geq 0$.
 Similarly, notice that 
 $$\nabla_y\Delta_h W(s)= \Delta_h\nabla_y W(s)=\Delta_hW(s) -\Delta_h W(s-y), $$
 by the power mean inequality we have uniformly in $h\geq 0$,
 \beqnn
 J_{2}(h)\ar\leq\ar C\cdot\int_{4h}^\infty ds \int_{s-h}^{s} \frac{|\Delta_h W(s) |^p}{ y^{\alpha+2} } dy+ C\cdot \int_{4h}^\infty ds \int_0^h   \frac{|\Delta_h W (y) |^p }{ (s-y)^{\alpha+2} } dy.
 \eeqnn
 Here the change of variables is also used to get the second integral.
 From (\ref{UpperboundDiffW}) we have $|\Delta_h W(s)| \leq C s^{\alpha-1}h$ for any $s\geq h\geq 0$ and $|\Delta_h W (y) | \leq  Ch^\alpha$ for any $y\in(0,h]$.
 Thus
 \beqnn
 J_{2}(h)
 \ar \leq\ar  Ch^p \int_{4h}^\infty  s^{p(\alpha-1)}ds \int_{s-h}^{s}   \frac{dy}{y^{\alpha+2}}+ C h^{p\alpha+1}   \int_{4h}^\infty  \frac{ds}{ (s-h)^{\alpha+2}}\cr
 \ar\leq\ar Ch^{p+1} \int_{4h}^\infty \frac{s^{p(\alpha-1)}}{(s-h)^{\alpha+2} }ds     + C h^{p\alpha+1}   \int_{3h}^\infty   s^{-\alpha-2}ds \leq C h^{(p-1)\alpha}.
 \eeqnn
 We now turn to consider $J_3(h)$ and $J_4(h)$.
 By (\ref{UpperboundW}), we have
 \beqnn
 |\nabla_y\Delta_h W(s)|
 \ar=\ar \Big|\int_0^yd\tilde{y}\int_0^h W''(s+\tilde{h}-\tilde{y})d\tilde{h}\Big|\cr
 \ar\leq\ar  C \int_0^yd\tilde{y}\int_0^h |s+\tilde{h}-\tilde{y}|^{\alpha-2}d\tilde{h}
 \leq C \cdot h\cdot \frac{ |s-y|\wedge y}{(s-y)^{2-\alpha}},
 \eeqnn
  uniformly in $h\geq 0$ and $s\geq  y> 0$. 
 Plugging this into $J_3(h)$ and $J_4(h)$ yields that uniformly in $h\geq 0$,
 \beqnn
 J_3(h) \ar\leq\ar Ch^p\int_{4h}^\infty ds \int_{s/2}^{s-h}  \frac{(s-y)^{p(\alpha-1)}}{  y^{\alpha+2} } dy
 \leq C h^{ p\alpha}\int_{4h}^\infty s^{-\alpha-1} ds
 \leq C h^{(p-1)\alpha}
 \eeqnn
 and
 \beqnn
 J_4(h) \ar\leq\ar Ch^p\int_{4h}^\infty \big| s/2 \big|^{p(\alpha-2)} ds \int_0^{s/2}   y^{p-\alpha-2}dy
 \leq Ch^{(p-1)\alpha}.
 \eeqnn
 The desired result follows immediately by putting these estimates together.
 \qed

 \begin{proposition} \label{WeakConvergenceProp03}
 For $p\geq 2$,  there exists a constant $C>0$ such that for any $h\geq 0$,
 \beqnn
 \int_0^\infty ds \int_{s}^{s+h} \frac{|\nabla_y W (s+h) -W (s) |^p}{ y^{\alpha+2}} dy \leq C h^{(p-1)\alpha}.
 \eeqnn
 \end{proposition}
 \proof  We split the preceding double integral into the following three parts:
 \beqnn
 J_1(h)\ar:=\ar  \int_0^h ds \int_{s}^{h} \frac{|\nabla_y W (s+h) -W (s) |^p}{ y^{\alpha+2}} dy,\cr
 J_2(h)\ar:=\ar \int_0^h ds \int_{h}^{s+h} \frac{|\nabla_y W (s+h) -W (s) |^p}{ y^{\alpha+2}} dy ,\cr
 J_3(h)\ar:=\ar \int_h^\infty ds\int_{s}^{s+h} \frac{|\nabla_y W (s+h) -W (s) |^p}{ y^{\alpha+2}} dy.
 \eeqnn
 By the power mean inequality, we have uniformly in $h\geq 0$,
 \beqnn
 J_1(h)\ar\leq \ar
 C\int_0^h ds \int_{s}^{h} \frac{|\nabla_y W(s+h) |^p}{ y^{\alpha+2} }dy+  C\int_0^h  ds \int_{s}^{h} \frac{|W(s)|^p}{y^{\alpha+2}}dy.
 \eeqnn
 Since $W(s)\leq Cs^\alpha$ uniformly in $s\geq 0$; see (\ref{UpperboundW}), the second term on the right side of this inequality can be bounded by $ C\int_0^h s^{p\alpha-\alpha-1}  ds  \leq  C\cdot h^{(p-1)\alpha}$ uniformly in $h\in[0,1]$.
 For the first term, choosing a positive constant $\theta $ satisfying that  $1+\alpha<p\theta<(1-\alpha)^{-1} \wedge p$, by (\ref{UpperboundDiffW}) we have
 		\beqnn
 		\nabla_y W(s+h) 
 		\ar=\ar |\nabla_y W(s+h)|^{1-\theta}\cdot |\nabla_y W(s+h)|^{\theta} \cr
 		\ar\ar\cr
 		\ar\leq\ar  C (s+h)^{(1-\theta)\alpha } (s+h-y)^{\theta(\alpha-1)} y^{\theta},
 		\eeqnn
 		uniformly in $h\geq s\geq 0$ and $y\in[h,s+h]$.
 		Then
 		\beqnn
 		\lefteqn{\int_0^h ds \int_{s}^{h} \frac{|\nabla_y W(s+h) |^p }{y^{\alpha+2}} dy}\ar\ar\cr
 		\ar\leq\ar C \int_0^h ds \int_{s}^{h}(s+h)^{(1-\theta)p\alpha } (s+h-y)^{p\theta(\alpha-1)} y^{p\theta-\alpha-2}dy\cr
 		\ar\leq\ar C\cdot h^{(1-\theta)p\alpha }\int_0^h ds \int_{s}^{h} (s+h-y)^{p\theta(\alpha-1)} y^{p\theta-\alpha-2}dy\cr
 		\ar\leq\ar C\cdot h^{(1-\theta)p\alpha }\int_0^h s^{p\theta-\alpha-2} ds \int_0^s (y+h-s)^{p\theta(\alpha-1)} dy\cr
 		\ar\leq\ar C\cdot h^{(1-\theta)p\alpha +p\theta(\alpha-1)+1}\int_0^h s^{p\theta-\alpha-2} ds \cr\ar\ar\cr
 		\ar\leq\ar  Ch^{(p-1)\alpha}.
 		\eeqnn
 		Hence  $J_1(h) \leq  C\cdot h^{(p-1)\alpha}$ uniformly in $h \geq 0$.
 		We turn to consider $J_2(h)$.
 		Notice that 
 		$$\nabla_y W(s+h) - W(s)=\Delta_h W(s) - W(s+h-y)$$ 
 		for any $y\in[h,s+h]$.
 		By the power mean inequality,
 		\beqnn
 		J_2(h)\ar\leq\ar  C \int_0^h ds \int_{h}^{s+h} \frac{|\Delta_h W(s) |^p}{ y^{\alpha+2} }dy + C\int_0^h ds \int_{h}^{s+h} \frac{|  W(s+h-y)|^p}{ y^{\alpha+2}}dy,
 		\eeqnn
 		uniformly in $h\geq 0$.
 		Since $|\Delta_hW(s)| \leq C(s+h)^\alpha$ and $W(s+h-y) \leq C (s+h-y)^\alpha$ uniformly in $s,h\geq 0$ and $y\in[0,s+h]$, there exists a constant $C>0$ such that for any $h\geq 0$,
 		\beqnn
 		J_2(h)\ar\leq\ar  C\int_0^h  ds \int_{h}^{s+h} \frac{(s+h)^{p\alpha} }{y^{\alpha+2} } dy
 		+ C\int_0^hds \int_{h}^{s+h} \frac{ (s+h-y)^{p\alpha} }{y^{\alpha+2} } dy .
 		\eeqnn
 		The first double integral on the  right side of this inequality can be bounded by $Ch^{-\alpha-1}\int_0^h  (s+h)^{p\alpha} ds \leq Ch^{(p-1)\alpha}.$
 		Using the change of variables and then Fubini's theorem to the second double integral,
 		\beqnn
 		\int_0^hds \int_{h}^{s+h}  \frac{ (s+h-y)^{p\alpha} }{y^{\alpha+2} }  dy
 		\ar=\ar \int_0^hds \int_{0}^{s} \frac{(s-y)^{p\alpha} }{(y+h)^{\alpha+2}} dy \cr
 		\ar=\ar  \int_0^h  ds \int_s^h \frac{(y-s)^{p\alpha}}{(s+h)^{\alpha+2}} dy
 		\leq C\cdot h^{(p-1)\alpha}
 		\eeqnn
 		and hence $J_2(h)\leq C\cdot h^{(p-1)\alpha}$ uniformly in $h\geq 0$.
 		Similarly, we also have
 		\beqnn
 		J_3(h)\ar\leq \ar  \int_h^\infty ds\int_{s}^{s+h} \frac{|\Delta_h W(s) |^p}{ y^{\alpha+2} } dy
 		+ \int_h^\infty ds \int_{s}^{s+h}  \frac{|  W(s+h-y) |^p}{ y^{\alpha+2} } dy \cr
 		\ar\leq\ar Ch^p \int_h^\infty  ds \int_{s}^{s+h}   \frac{s^{p(\alpha-1)}}{y^{\alpha+2} } dy
 		+C \int_h^\infty ds \int_{s}^{s+h}\frac{ (s+h-y)^{p\alpha}}{ y^{\alpha+2} } dy .
 		\eeqnn
 		The first term on the right side of the second inequality can be bounded by
 		$$C\cdot h^p \int_h^\infty  s^{p(\alpha-1)-\alpha-1}ds  \leq C\cdot h^{(p-1)\alpha},$$ 
 		uniformly in $h\geq 0$.
 		By the change of variables, the second term equals to
 		\beqnn
 		C \int_h^\infty ds\int_{0}^{h}  \frac{(h-y)^{p\alpha} }{(y+s)^{\alpha+2}} dy
 		\leq C \int_h^\infty   ds\int_{0}^{h} \frac{y^{p\alpha} }{s^{\alpha+2}}  dy \leq  C\cdot h^{(p-1)\alpha}
 		\eeqnn
 		and $J_3(h)\leq Ch^{(p-1)\alpha}$ uniformly in $h\geq 0$. The desired result follows by putting these estimates together.
 		\qed

 		\begin{proposition}\label{Prop.D4}
 			For any $x\geq 0$, we have as $\epsilon\to 0+$,
 			\beqnn
 			\int_0^x ds\int_0^\epsilon \frac{|\nabla_y W(s)|^2}{y^{\alpha+2}} dy +\int_0^x ds\int_0^\epsilon \frac{|\int_0^s \nabla_y K(r)dr|^2}{y^{\alpha+2}}dy  \to 0.
 			\eeqnn
 		\end{proposition}
 		\proof Here we just prove the convergence of the first integral to $0$.
 		The second one can be prove in the same way.
 		For convenience, we assume $x\geq \epsilon\geq 0$. By Fubini's theorem and the fact that $W(s-y)=0$ for $y\geq s$, we can split the targeted integral into two parts
 		\beqnn
 		I_1(\epsilon):= \int_0^\epsilon dy \int_0^y  \frac{|W(s)|^2}{y^{\alpha+2}}ds
 		\quad \mbox{and}\quad
 		I_2(\epsilon):= \int_0^\epsilon dy \int_y^x  \frac{|\nabla_y W(s)|^2}{y^{\alpha+2}} ds.
 		\eeqnn
 		By (\ref{UpperboundW}), we have 
 		$\int_0^y |W(s)|^2 ds  \leq C \cdot y^{2\alpha+1}$ uniformly in $y \geq 0$ and hence
 		\beqnn
 		\lim_{\epsilon \to 0+} I_1(\epsilon) \leq C\cdot \epsilon^{\alpha} = 0. 
 		\eeqnn 
 		Let $\vartheta\in (\alpha+1,(1-\alpha)^{-1}\wedge 2)$.  Using (\ref{UpperboundDiffW}), we  have uniformly in $s\geq y\geq 0$,
 		\beqnn
 		\big|\nabla_y W(s) \big|^2 =  \big|\nabla_y W(s) \big|^{\vartheta} \cdot \big|\nabla_y W(s) \big|^{2-\vartheta}\leq C\cdot y^\theta(s-y)^{\vartheta(\alpha-1)}\cdot s^{(2-\vartheta)\alpha}.
 		\eeqnn
 		Plugging this into $I_2(\epsilon)$ implies that
 		\beqnn
 		I_2(\epsilon) \ar\leq\ar C\cdot  \int_0^\epsilon dy  \int_y^x  \frac{(s-y)^{\vartheta(\alpha-1)}}{y^{\alpha+2-\vartheta}}\cdot s^{(2-\vartheta)\alpha} ds \cr
 		\ar \leq\ar C\cdot x^{(2-\vartheta)\alpha} \int_0^\epsilon \frac{(x-y)^{\vartheta(\alpha-1)+1}}{y^{\alpha+2-\vartheta}} dy
 		\leq C \cdot \epsilon^{\vartheta-\alpha-1},
 		\eeqnn
 		which goes to $0$ as $\epsilon\to 0 + $. The desired result follows by putting these estimates together.
 		\qed

 		\begin{proposition}\label{Prop.B5}
 			For any $x\geq 0$, we have as $\epsilon\to 0+$,
 			\beqnn
 			\int_0^x ds\int_0^\epsilon\frac{|K*\nabla_y W(s)|^2}{y^{\alpha+2}} dy \to 0.
 			\eeqnn
 		\end{proposition}
 		\proof We still assume $x\geq \epsilon\geq 0$ and then split the targeted double integral into the following two parts
 		\beqnn
 		I_3(\epsilon)\ar:=\ar \int_0^\epsilon dy \int_0^y\frac{|K* W(s)|^2}{y^{\alpha+2}} ds
 		\quad \mbox{and}\quad
 		I_4(\epsilon):= \int_0^\epsilon dy \int_y^x \frac{|K*\nabla_y W(s)|^2}{y^{\alpha+2}} ds.
 		\eeqnn
 		By (\ref{UpperboundW}), we have uniformly in $y\geq 0$,
 		\beqnn
 		\int_0^y\big|K* W(s)\big|^2 ds
 		\ar\leq\ar C\cdot  \int_0^y\Big|\int_0^s (s-r)^{\alpha-1}r^\alpha dr\Big|^2 ds \leq C\cdot  \int_0^y s^{4\alpha} ds  \leq C \cdot y^{4\alpha+1} .
 		\eeqnn
 		Taking  this into $ I_3(\epsilon) $, we have $I_3(\epsilon)\leq C\cdot \epsilon^{3\alpha}\to 0$ as $\epsilon\to 0+$.
 		By the Cauchy-Schwarz inequality, we have $I_4(\epsilon) \leq C\big[I_{41}(\epsilon)+I_{42}(\epsilon)\big]$ uniformly in $x\geq \epsilon\geq 0$, where
 		\beqnn
 		I_{41}(\epsilon)\ar :=\ar  \int_0^\epsilon \frac{dy}{y^{\alpha+2}} \int_y^x\Big|\int_0^y (s-r)^{\alpha-1}  W(r)dr\Big|^2 ds ,\cr
 		I_{42}(\epsilon)\ar :=\ar  \int_0^\epsilon \frac{dy}{y^{\alpha+2}} \int_y^x\Big|\int_y^s (s-r)^{\alpha-1} \nabla_y W(r)dr\Big|^2 ds.
 		\eeqnn
 		By (\ref{UpperboundW}) and the fact that $(s-r)^{\alpha-1} \leq (y-r)^{\frac{\alpha-1}{2}}\cdot (s-y)^{\frac{\alpha-1}{2}}$ for $s>y>r>0$, there exists a constant $C>0$ such that for any $x\geq y\geq 0$,
 		\beqnn
 		\int_y^x\Big|\int_0^y (s-r)^{\alpha-1}  W(r)dr\Big|^2 ds
 		\ar\leq\ar C\cdot \int_y^x\Big|\int_0^y (y-r)^{\frac{\alpha-1}{2}}r^\alpha dr\Big|^2 (s-y)^{ \alpha-1 } ds\cr
 		\ar\ar\cr
 		\ar\leq\ar C\cdot x^\alpha \cdot y^{3\alpha+1}.
 		\eeqnn
 		Plugging this back into $I_{41}(\epsilon)$, we have $I_{41}(\epsilon) \leq C\cdot \epsilon^{2\alpha}\to 0$ as $\epsilon \to 0+$.
 		Let $\vartheta\in(\frac{\alpha+1}{2}, \frac{\alpha+1/2}{1-\alpha}\wedge 1)$.
 		By (\ref{UpperboundDiffW}), we have 
 		$$| \nabla_y W(r)| \leq C\cdot (r-y)^{(\alpha-1)\vartheta}y^\vartheta \cdot x^{\alpha(1-\vartheta)} ,$$ 
 		uniformly in $0\leq y\leq r\leq x$ and hence
 		\beqnn
 		\int_y^x\Big|\int_y^s (s-r)^{\alpha-1} \nabla_y W(r)dr\Big|^2 ds
 		\ar\leq\ar C\cdot  y^{2\vartheta} \int_y^x\Big|\int_y^s (s-r)^{\alpha-1} (r-y)^{(\alpha-1)\vartheta}  dr\Big|^2 ds \cr
 		\ar\leq\ar C \cdot y^{2\vartheta} \int_y^x (s-y)^{2\alpha+2(\alpha-1)\vartheta}  ds
 		\leq  C \cdot   y^{2\vartheta} .
 		\eeqnn
 		Taking this back into $I_{42}(\epsilon)$, we have $I_{42}(\epsilon)\leq  C\cdot \epsilon^{2\theta-\alpha-1} \to 0$ as $\epsilon\to 0+$.
 		The desired result follows by putting all results above together.
 		\qed

 		\begin{proposition}\label{Prop.B6}
 			For any $x\geq 0$, we have as $\epsilon\to 0+$,
 			\beqnn
 			\int_0^x ds\int_0^\epsilon \frac{|W* \nabla_y K(s)|^2}{y^{\alpha+2}}dy \to 0 .
 			\eeqnn
 		\end{proposition}
 		\proof Notice that  $|W* \nabla_y K(s)| \leq C\cdot x^\alpha\cdot \int_0^s \nabla_y K(r)dr$ uniformly in $x\geq s\geq 0$.
 		Thus
 		\beqnn
 		\int_0^x ds\int_0^\epsilon \frac{|W* \nabla_y K(s)|^2}{y^{\alpha+2}}dy
 		\leq C\cdot \int_0^x ds\int_0^\epsilon \frac{|\int_0^s \nabla_y K(r)dr|^2}{y^{\alpha+2}}dy ,
 		\eeqnn
 		which goes to $0$ as $\epsilon\to 0+$; see Proposition~\ref{Prop.D4}.
 		\qed

 	\section{Marked Hawkes point measures}\label{Appendix.HP}

 	Let $\mathbb{U}$ be a Lusin topological space endowed with the Borel $\sigma$-algebra $\mathscr{U}$.
 	Let $\{\sigma_k:k=1,2\cdots\}$ be a sequence of increasing, $(\mathscr{F}_t)$-adapted random times and $\{ \eta_k:k=1,2,\cdots \}$ be a sequence of i.i.d. $\mathbb{U}$-valued random variables with distribution $\nu_H(du)$.
 	We assume that $\eta_k$ is independent of  $\{\sigma_j:j=1,\cdots,k \}$ for any $k\geq 0$.
 	In terms of these two sequences we define the $(\mathscr{F}_t)$-random point measure
 	\beqnn
 	N_H(ds,du):= \sum_{k=1}^\infty \mathbf{1}_{\{\sigma_k\in ds, \eta_k\in du  \}}
 	\eeqnn
 	on $(0,\infty)\times \mathbb{U}$.
 	We say $N_H(ds,du)$  is a \textit{marked Hawkes point measure} (MHP) on $(0,\infty)\times\mathbb{U}$ if it has a $(\mathscr{F}_t)$-intensity $Z(s-)\cdot ds\cdot \nu_{H}(du)$ with the intensity process $Z:=\{ Z(t):t\geq 0 \}$ given by
 	\beqnn
 	Z(t)\ar=\ar \mu(t) + \sum_{k=1}^{N_H(t)} \phi(t-\sigma_k, \eta_k), \quad t\geq 0,
 	\eeqnn
 	for some \textsl{kernel} $\phi: \mathbb{R}_+\times \mathbb{U} \to [0,\infty)$ and some $\mathscr{F}_0$-measurable, non-negative functional-valued random variable $\{  \mu(t) :t\geq 0 \}$.
 	We usually interpret $\phi(\cdot,u)$ and $\mu$ as the impacts of an event with mark $u$ and all events prior to time $0$ on the arrival of future events respectively. 
 	Following the argument in \cite[p.93]{IkedaWatanabe1989}; see also Section~2 in \cite{HorstXu2021}, on an extension of the original probability space we can define a time-homogeneous PRM $N(ds,du,dz)$ on $(0,\infty)\times \mathbb{U}\times \mathbb{R}_+$ with intensity $ ds\cdot\nu_H(du)\cdot dz$ such that
 	\beqnn
 	N_H(ds,du)\ar=\ar  \int_0^{Z(s-)}  N(ds,du,dz)
 	\eeqnn
 	and hence the intensity process at time $t$ can be rewritten into
 	\beqnn
 	Z(t) \ar=\ar \mu (t) +\int_0^t \int_\mathbb{U}\int_0^{Z(s-)} \phi(t-s, u)N(ds,du,dz) .
 	\eeqnn
 	
 	Denote by $\phi_H:=\{\phi_H(t):t\geq 0\}$ the mean impacts of an event on the arrival of future events with
 	$$ \phi_H(t):= \int_\mathbb{U} \phi(t,u) \nu_H(du).$$
 	We assume $\phi_H$ is locally integrable.
 	Let $R_H:=\{R_H(t):t\geq 0\}$ be the \textit{resolvent} of $\phi_H$ defined as the unique solution to
 	$$
 	R_H(t) = \phi_H(t)  +  \phi_H* R_H(t) .
 	$$
 	It is usual to interpret $R_H$ as the mean impacts of an event and its triggered events on the arrivals of future events.
 	In addition, we introduce a two-parameter function 
 	$$R(t, u) = \phi(t,u) + R_H*\phi(t,u),\quad (t,u) \in \mathbb{R}_+\times \mathbb{U}$$
 	to describe the mean impacts of an event with mark $u$ on the arrivals of future events.
 	An argument similar to the one used in Section~2 in \cite{HorstXu2021} induces the following proposition immediately.
 	
 	\begin{proposition}\label{SVR} \it
 		The intensity process $Z$ satisfies the following SVE
 		\beqnn
 		Z(t)
 		\ar=\ar  \mu(t) + \int_0^t R_H(t-s)\mu(s)ds\cr
 		\ar\ar  +\int_0^t \int_\mathbb{U}\int_0^{Z(s-)}  R(t-s, u)\widetilde{N}(ds,du,dz), \quad t\geq 0, 
 		\eeqnn
 		where $\widetilde{N}(ds,du,dz):=N(ds,du,dz)- ds\cdot\nu_H(du)\cdot dz$.
 	\end{proposition}

 	\section{Stochastic integrals with respect to $\mathbb{H}^\#$-semimartingale} \label{AppendixB} 
 	
 	In this section we give a brief introduction to the stochastic integrals with respect to infinite-dimensional semimartingales; readers may refer to \cite{KurtzProtter1996} for more details.
 	Let $\mathbb{H}$ be a separable Banach space endowed with a norm $\|\cdot\|_{\mathbb{H}}$.
 	We now give the definition of  $\mathbb{H}^\#$-semimartingales.
 	
 	\begin{definition}
 		We say	$Y$ is a $(\mathscr{F}_t)$-adapted \textit{$\mathbb{H}^\#$-semimartingale}, if it is a stochastic process indexed by $\mathbb{H} \times \mathbb{R}_+$ such that
 		\begin{enumerate}
 			\item[$\bullet$] For each $f\in \mathbb{H}$,  $Y(f):=\{ Y(f,t):t\geq 0 \}$ is a c\'adl\'ag $(\mathscr{F}_t)$-semimartingale with $Y(f,0)\overset{\rm a.s.}=0$;
 			
 			\item[$\bullet$] For each $t\geq 0$, $\alpha_1,\cdots,\alpha_m\in\mathbb{R}$ and $f_1,\cdots,f_m\in\mathbb{H}$, 
 			\beqnn
 			Y\bigg(\sum_{k=1}^m \alpha_k f_k,t\bigg)\overset{\rm a.s.}=\sum_{k=1}^m \alpha_k Y(f_k,t).
 			\eeqnn
 		\end{enumerate}
 	\end{definition}
 	
 	Let $\mathbb{H}_0$ be a dense subset of $\mathbb{H}$ and $\mathcal{S}_0$ the collection of $\mathbb{H}$-valued stochastic processes of the form
 	\beqnn
 	X(t):= \sum_{k=1}^m \xi_k(t)\varphi_k \quad \mbox{with}\quad \xi_k(t):=\sum_{i=0}^{\infty} \eta_i^k\cdot\mathbf{1}_{[\tau_i^k,\tau_{i+1}^k)}(t),
 	\eeqnn
 	where $m\geq 1$, $\varphi_1,\cdots,\varphi_m\in\mathbb{H}_0$, $\{\tau_i^k\}_{i\geq 0}$ is a sequence of non-decreasing $(\mathscr{F}_t)$-stopping times and $\eta_i^k \in\mathbb{R}^d$ is  $\mathscr{F}_{\tau_i^k}$-measurable.
 	For any $X\in\mathcal{S}_0$,  define
 	\beqnn
 	X_-\cdot Y(t) =\sum_{k=1}^m \int_0^t \xi_k(s-)d Y(\varphi_k,t), \quad t\geq 0.
 	\eeqnn
 	\begin{definition}
 		The $\mathbb{H}^\#$-semimartingale $Y$ is \textit{standard} if
 		\beqlb\label{eqn.Ht}
 		\mathcal{H}_t:=  \Big\{ \sup_{s\leq t}|X_-\cdot Y(s)| : X \in\mathcal{S}_0,\, \sup_{s\leq t}\|X(s)\|_{\mathbb{H}}\leq 1  \Big\}
 		\eeqlb
 		is stochastically bounded for each $t\geq 0$.
 	\end{definition}
 	For any $\mathbb{H}$-valued c\'adl\'ag process $X$ and standard $\mathbb{H}^\#$-semimartingale $Y$,
 	we can find a sequence $\{X^\epsilon\}_{\epsilon>0}\subset \mathcal{S}_0$ such that as $\epsilon \to0$,
 	\beqnn
 	\sup_{t\in[0,T]}\|X^\epsilon(t)-X(t)\|_\mathbb{H}\overset{\rm  a.s. }\to0
 	\quad \mbox{and} \quad
 	X_-\cdot Y := \lim_{\epsilon\to 0+} X^\epsilon_-\cdot Y
 	\eeqnn
 	exists a.s. in the sense that $\sup_{t\in[0,T]} |X_-\cdot Y(t)- X^\epsilon_-\cdot Y(t) | \overset{\rm p}\to 0$.
 	Moreover, the limit process $X_-\cdot Y$ is c\'adl\'ag,  independent of $\{X^\epsilon\}_{\epsilon>0}$ and called the \textit{stochastic integral} of $X$ with respect to $Y$. 
 	For any $(\mathscr{F}_t)$-stopping time $\sigma$, we have
 	\beqnn
 	X_-\cdot Y(t\wedge \sigma)= X_-^\sigma\cdot Y(t)
 	\quad\mbox{with}\quad
 	X_-^\sigma(t):= X_-(t) \cdot \mathbf{1}_{[0,\sigma)}(t),
 	\quad t\geq 0. 
 	\eeqnn 
 	
 	\begin{definition}\label{Definition.A1}
 		A sequence of $\mathbb{H}^\#$-semimartingales $\{Y_n\}_{n\geq 1}$ is \textit{uniformly tight} if $\{\mathcal{H}_{n,t}\}_{n\geq 1}$ is uniformly stochastically bounded for any $t\geq 0$, where $\mathcal{H}_{n,t}$ is defined by (\ref{eqn.Ht}) with $Y$ replaced by $Y_n$.
 		
 		We say $Y_n$ {\rm converges weakly} to $Y$ and write $Y_n \Rightarrow Y$ if for any $m\geq 1$ and $f_1,\cdots,f_m\in\mathbb{H}$,
 		\beqnn
 		(Y_n(f_1),\cdots, Y_n(f_m))\overset{\rm d}\to(Y(f_1),\cdots, Y(f_m))\quad \mbox{in }D([0,\infty),\mathbb{R}^m).
 		\eeqnn
 		In addition,
 		we also write $(X_n,Y_n) \Rightarrow (X,Y)$ if
 		\beqnn
 		(X_n, Y_n(f_1),\cdots,Y_n(f_m))\overset{\rm d}\to(X,Y(f_1),\cdots, Y(f_m))\quad \mbox{in }D([0,\infty),\mathbb{H}\times\mathbb{R}^m).
 		\eeqnn
 	\end{definition}

 	\section{Stochastic integrals with respect to Poisson random measure}

 	Let $\widetilde{N}_1(ds,dy,dz)$ be a compensated $(\mathscr{F}_t)$-PRM on $(0,\infty)^3$  with intensity $ds\cdot\nu_1(dy)\cdot dz$, where $\nu_1(dy)$ is a $\sigma$-finite measure on $\mathbb{R}_+$ such that $\nu_1(x,\infty)<\infty$ for any $x>0$.
 	Let $\{X(t):t\geq 0\}$ be a $(\mathscr{F}_t)$-predictable and non-negative process.
 	
 	\begin{theorem}[Maximal inequality] \label{Thm.BDG}
 		For $p\geq 1$ and $T>0$, let $f $ be a measurable function on $\mathbb{R}_+^2$ satisfying
 		$$\int_0^Tds \int_0^\infty \big|f(s,y)\big|^{2p}\nu_1(dy) <\infty. $$
 		If $\sup_{t\in[0,T]}\mathbf{E}\big[\big|X(t)\big|^{p}\big]<\infty$,  there exists a constant $C>0$ depending only on $p$ such that
 		\beqlb \label{BDG}
 		\lefteqn{\mathbf{E}\Big[\sup_{t\in[0,T]}\Big| \int_0^t \int_0^\infty \int_0^{X(s)} f(s,y)\widetilde{N}_1(ds,dy,dz) \Big|^{2p}  \Big]}\ar\ar\cr
 		\ar\leq\ar C  \sup_{t\in[0,T]} \mathbf{E} \Big[\big|X(t)\big|^p \Big] \cdot \Big| \int_0^T  \int_0^\infty\big|f(s,y)\big|^{2} \nu_1(dy)  ds\Big|^{p} \cr
 		\ar\ar + C  \sup_{t\in[0,T]}\mathbf{E}\Big[\big|X(t)\big|\Big]\cdot \int_0^T  \int_0^\infty\big|f(s,y)\big|^{2p} \nu_1(dy)  ds .
 		\eeqlb
 	\end{theorem}
 	\proof  
 	By the maximal inequality established in \cite[Theorem 1, p.297]{MarinelliRockner2014} for purely discontinuous martingales, there exists a constant $C>0$ depending only on $p$  such that
 	\beqnn
 	\lefteqn{\mathbf{E}\Big[ \sup_{t\in[0,T]}\Big| \int_0^t \int_0^\infty \int_0^{X(s)} f(s,y)\widetilde{N}_1(ds,dy,dz) \Big|^{2p}  \Big]}\ar\ar\cr
 	\ar\leq\ar C \cdot \mathbf{E}\Big[ \Big| \int_0^T X(s) ds \int_0^\infty|f(s,y)|^{2} \nu_1(dy) \Big|^{p}\Big]\cr
 	\ar\ar +  C \cdot \mathbf{E}\Big[ \int_0^T X(s) ds \int_0^\infty|f(s,y)|^{2p} \nu_1(dy) \Big] \cr
 	\ar\leq\ar C \cdot \mathbf{E}\Big[ \Big| \int_0^T X(s) ds \int_0^\infty|f(s,y)|^{2} \nu_1(dy) \Big|^{p}\Big] \cr
 	\ar\ar + C \sup_{t\in[0,T]}\mathbf{E}\Big[\big|X(t)\big|\Big]\cdot\int_0^T  ds \int_0^\infty|f(s,y)|^{2p} \nu_1(dy) .
 	\eeqnn
 	By H\"older's inequality; see footnote~\ref{Footnote.Holder}, the first expectation on  the right side of the second inequality can be bounded by
 	\beqnn
 	\ar\ar \int_0^T \mathbf{E} \Big[\big|X(r)\big|^p \Big] \int_0^\infty|f(r,y)|^{2} \nu_1(dy) dr
 	\cdot \Big| \int_0^T \int_0^\infty|f(s,z)|^{2} \nu_1(dz)ds  \Big|^{p-1}  ,
 	\eeqnn
 	which can be bounded by 
 	$$ \sup_{t\in[0,T]} \mathbf{E} \Big[\big|X(t)\big|^p  \Big] \cdot\big| \int_0^T \int_0^\infty|f(s,y)|^{2} \nu_1(dy)ds  \big|^{p}.$$
 	The desired result holds.
 	\qed
 	
 	\begin{theorem}[Stochastic Fubini theorem]\label{StoFubiniThm}
 		Let $(\mathbb{V},\mathscr{V},m)$ be a measurable space.
 		For $T\geq 0$, let $f$ be a measurable function on $\mathbb{V}$ and $g,h$  two measurable functions on $ \mathbb{V}\times [0,T]\times \mathbb{R}_+$ satisfying that as $\epsilon \to 0+$,
 		\beqnn
 		\int_\mathbb{V} |f(v)|m(dv) \Big| \int_0^T ds \int_0^\epsilon |g(v,s,y)|^2\nu_1(dy)\Big|^{1/2} \to  0, \cr
 		\int_0^T ds \int_0^\epsilon \Big|\int_\mathbb{V} f(v)g(v,s,y)m(dv)\Big|^2\nu_1(dy) \to  0
 		\eeqnn
 		and
 		\beqnn
 		\int_\mathbb{V} |f(v)|m(dv) \int_0^T ds \int_0^\epsilon |h(v,s,y)|\nu_1(dy)   \ar\to\ar 0, \cr
 		\int_0^T ds \int_0^\epsilon \Big| \int_\mathbb{V} f(v)h(v,s,y)m(dv)\Big| \nu_1(dy)  \ar\to\ar 0.
 		\eeqnn
 		If $ \|X\|_{L^\infty_T}<\infty$ a.s., we have
 		\beqlb\label{SFubni}
 		\ar\ar\int_\mathbb{V} f(v)m(dv) \int_0^T \int_0^\infty \int_0^{X(s)}g(v,s,y) \widetilde{N}_1(ds,dy,dz) \cr
 		\ar\ar \quad  \overset{\rm a.s.}= \int_0^T \int_0^\infty \int_0^{X(s)} \int_\mathbb{V} f(v) g(v,s,y) m(dv) \widetilde{N}_1(ds,dy,dz)
 		\eeqlb
 		and
 		\beqlb\label{SFubni01}
 		\ar\ar\int_\mathbb{V} f(v)m(dv) \int_0^T \int_0^\infty \int_0^{X(s)}h(v,s,y) N_1(ds,dy,dz) \cr
 		\ar\ar \quad  \overset{\rm a.s.}= \int_0^T \int_0^\infty \int_0^{X(s)} \int_\mathbb{V} f(v) h(v,s,y) m(dv) N_1(ds,dy,dz).
 		\eeqlb
 	\end{theorem}
 	\proof
 	Here we just prove the first the equality (\ref{SFubni}) and the second one can be proved in the same way.
 	It is easy to identify that the two integrals  in (\ref{SFubni}) are well-defined. We now show they are equal almost surely. For any $\epsilon\in(0,1]$, by the assumption that $\nu_1(\epsilon,\infty)<\infty$ and  Fubini's theorem,
 	\beqnn
 	\lefteqn{\int_\mathbb{V} f(v)m(dv)\int_0^T \int_\epsilon^\infty \int_0^{X(s)} g(v,s,y) \widetilde{N}_1(ds,dy,dz)}\ar\ar\cr
 	\ar=\ar \int_\mathbb{V} f(v)m(dv)\int_0^T \int_\epsilon^\infty \int_0^{X(s)} g(v,s,y) N_1(ds,dy,dz)\cr
 	\ar\ar -\int_\mathbb{V} f(v) m(dv)\int_0^T \int_\epsilon^\infty  X(s) g(v,s,y) ds\nu(dy)\cr
 	\ar=\ar \int_0^T \int_\epsilon^\infty \int_0^{X(s)} \int_\mathbb{V} f(v) g(v,s,y) m(dv) N_1(ds,dy,dz)\cr
 	\ar\ar -\int_0^T \int_\epsilon^\infty X(s) \int_\mathbb{V} f(v)g(v,s,y) m(dv) ds\nu(dy)\cr
 	\ar=\ar\int_0^T \int_\epsilon^\infty \int_0^{X(s)}\int_\mathbb{V} f(v) g(v,s,y) m(dv)\widetilde{N}_1(ds,dy,dz).
 	\eeqnn
 	Here the two stochastic integrals with respect to $N_1(ds,dy,dz)$ on the right side of the first and second equality are finite sums.
 	Thus the difference between the two integrals in (\ref{SFubni}) can be bounded by $|A_1(\epsilon) | +|A_2(\epsilon)|$ uniformly in $\epsilon\in(0,1]$ with
 	\beqnn
 	A_1(\epsilon)\ar:=\ar \int_\mathbb{V} f(v) m(dv)\int_0^T \int_0^\epsilon \int_0^{X(s)}g(v,s,y) \widetilde{N}_1(ds,dy,dz),\cr
 	A_2(\epsilon)\ar:=\ar \int_0^T \int_0^\epsilon \int_0^{X(s)} \int_\mathbb{V} f(v) g(v,s,y) m(dv) \widetilde{N}_1(ds,dy,dz).
 	\eeqnn
 	It suffices to prove  that $|A_1(\epsilon) | +|A_2(\epsilon)|\overset{\rm p}\to 0$ as $\epsilon\to 0+$.
 	For any $\eta>0$ and $J> 0$, we have 
 	$$\mathbf{P}(| A_1(\epsilon)|\geq \eta) \leq \mathbf{P}(| A_1(\epsilon)|\geq \eta, \|X\|_{L^\infty_T}\leq J)+ \mathbf{P}(  \|X\|_{L^\infty_T}> J).$$
  Since $ \|X\|_{L^\infty_T}<\infty$ a.s., then $\mathbf{P}(  \|X\|_{L^\infty_T}> J)\to0$ as $J\to\infty$. Notice that
 	\beqnn
 	\lefteqn{\mathbf{P}\big(| A_1(\epsilon)|\geq \eta, \|X\|_{L^\infty_T}\leq J\big)}\ar\ar\cr
 	\ar\ar\cr
 	\ar\leq\ar \mathbf{P}\Big(\Big|\int_\mathbb{V} f(v)m(dv) \int_0^T \int_0^\epsilon \int_0^{X(s)\wedge J}g(v,s,y) \widetilde{N}_1(ds,dy,dz) \Big|\geq \eta  \Big) .
 	\eeqnn
 	By Chebyshev's inequality and then Fubini's theorem,
 	\beqnn
 	\lefteqn{\mathbf{P}\big(| A_1(\epsilon)|\geq \eta, \|X\|_{L^\infty_T}\leq J\big)}\ar\ar\cr
 	\ar\ar\cr
 	\ar\leq\ar \frac{1}{\eta}\int_\mathbb{V} |f(v)| m(dv) \mathbf{E}\Big[\Big|\int_0^T \int_0^\epsilon \int_0^{X(s)\wedge J}g(v,s,y) \widetilde{N}_1(ds,dy,dz)\Big| \Big].
 	\eeqnn
 	By (\ref{BDG}) with $p=1/2$ and then Jensen's inequality, the last expectation  can be bounded by
 	\beqnn
 	\lefteqn{\mathbf{E}\Big[\Big|\int_0^T \int_0^\epsilon \int_0^{X(s)\wedge J}|g(v,s,y)|^2 N_1(ds,dy,dz)\Big|^{1/2} \Big]}\ar\ar\cr
 	\ar\leq\ar  \Big|\mathbf{E}\Big[\int_0^T \int_0^\epsilon \int_0^{X(s)\wedge J}|g(v,s,y)|^2 N_1(ds,dy,dz) \Big]\Big|^{1/2}\cr
 	\ar\leq \ar \Big|J \int_0^T ds\int_0^\epsilon  |g(v,s,y)|^2 \nu_1(dy)\Big|^{1/2}
 	\eeqnn
 	and hence
 	\beqnn
 	\lefteqn{\mathbf{P}\big(| A_1(\epsilon)|\geq \eta, \|X\|_{L^\infty_T}\leq J\big)}\ar\ar\cr
 	\ar\ar\cr
 	\ar\leq\ar \frac{\sqrt{J}}{\eta}\int_\mathbb{V} |f(v)|\cdot\Big| \int_0^T ds \int_0^\epsilon  |g(v,s,y)|^2 \nu_1(dy)\Big|^{1/2} m(dv),
 	\eeqnn
 	which goes to $0$ as $\epsilon\to 0+$. Putting these estimates together, we have 
 	\beqnn
 	\mathbf{P}\big(| A_1(\epsilon)|\geq \eta\big)\to 0
 	\quad\mbox{and hence}\quad
 	| A_1(\epsilon)| \overset{\rm p}\to 0,
 	\eeqnn
 	as  $\epsilon\to 0+$.
 	Similarly, we also can prove that $| A_2(\epsilon)| \overset{\rm p}\to 0$ as  $\epsilon\to 0+$.
 	\qed

 \end{appendix}

  \begin{acks}[Acknowledgments]
 	The author is grateful to Matthias Winkel who noticed the inaccuracy on the H\"older continuity and recommended several helpful references.
 	The author also like to thank the three professional referees for their careful and insightful reading of the paper, and for comments, which led to many improvements.  
 \end{acks}

 \end{document}